\documentclass[12pt,oneside,french]{article}
\usepackage{amssymb,amsmath,babel}
\ifx\optionkeymacros\undefined\else \fi

\catcode`\Œ=\active\defŒ{{\aa}}       
\catcode`\º=\active\defº{\int}        
\catcode`\=\active\def{\c c}        
\catcode`\¶=\active\def¶{\partial}    
\catcode`\Ä=\active\defÄ{\oint}       
\catcode`\Æ=\active\defÆ{\triangle}   
\catcode`\Â=\active\defÂ{\neg}        
\catcode`\µ=\active\defµ{\mu}         
\catcode`\¿=\active\def¿{{\o}}        
\catcode`\¹=\active\def¹{\pi}         
\catcode`\Ï=\active\defÏ{{\oe}}       
\catcode`\§=\active\def§{{\ss}}       
\catcode`\ =\active\def {\dagger}     
\catcode`\Ã=\active\defÃ{\sqrt}       
\catcode`\·=\active\def·{\Sigma}      
\catcode`\Å=\active\defÅ{\approx}     
\catcode`\½=\active\def½{\Omega}      
\catcode`\£=\active\def£{{\it\$}}     
\catcode`\°=\active\def°{\infty}      
\catcode`\¤=\active\def¤{{\S}}        
\catcode`\¦=\active\def¦{{\P}}        
\catcode`\¥=\active\def¥{\bullet}     
\catcode`\»=\active\def»{\leavevmode\raise.585ex\hbox{\b a}}      
\catcode`\¼=\active\def¼{\leavevmode\raise.6ex\hbox{\b o}}        
\catcode`\­=\active\def­{\not=}       
\catcode`\²=\active\def²{\leq}        
\catcode`\³=\active\def³{\geq}        
\catcode`\Ö=\active\defÖ{\div}        
\catcode`\É=\active\defÉ{{\dots}}     
\catcode`\¾=\active\def¾{{\ae}}       
\catcode`\Ç=\active\defÇ{\ll}         
\catcode`\Ò=\active\defÒ{``}          
\catcode`\Á=\active\defÁ{!`}          
\catcode`\¢=\active\def¢{\rlap/c}     
\catcode`\Ô=\active\defÔ{`}           
\catcode`\Õ=\active\defÕ{'}           

\catcode`\=\active\def{{\AA}}       
\catcode`\'=\active\def'{\c C}        
\catcode`\¯=\active\def¯{{\O}}        
\catcode`\¸=\active\def¸{\Pi}         
\catcode`\Î=\active\defÎ{{\OE}}       
\catcode`\®=\active\def®{{\AE}}       
\catcode`\×=\active\def×{\diamond}    
\catcode`\¡=\active\def¡{\accent'27}  
\catcode`\Ó=\active\defÓ{''}          
\catcode`\±=\active\def±{\pm}         
\catcode`\È=\active\defÈ{\gg}         
\catcode`\À=\active\defÀ{?`}          
\catcode`\Ð=\active\defÐ{--}          
\catcode`\Ñ=\active\defÑ{---}         

\catcode`\Š=\active\defŠ{\"a}        
\catcode`\'=\active\def'{\"e}        
\catcode`\•=\active\def•{\"{\i}}     
\catcode`\š=\active\defš{\"o}        
\catcode`\Ÿ=\active\defŸ{\"u}        
\catcode`\Ø=\active\defØ{\"y}        
\catcode`\€=\active\def€{\"A}        
\catcode`\…=\active\def…{\"O}        
\catcode`\†=\active\def†{\"U}        
\catcode`\‡=\active\def‡{\'a}        
\catcode`\Ž=\active\defŽ{\'e}        
\catcode`\'=\active\def'{\'{\i}}     
\catcode`\—=\active\def—{\'o}        
\catcode`\œ=\active\defœ{\'u}        
\catcode`\ƒ=\active\defƒ{\'E}        
\catcode`\ˆ=\active\defˆ{\`a}        
\catcode`\=\active\def{\`e}        
\catcode`\"=\active\def"{\`{\i}}     
\catcode`\˜=\active\def˜{\`o}        
\catcode`\=\active\def{\`u}        
\catcode`\Ë=\active\defË{\`A}        
\catcode`\‹=\active\def‹{\~a}        
\catcode`\–=\active\def–{\~n}        
\catcode`\›=\active\def›{\~o}        
\catcode`\Ì=\active\defÌ{\~A}        
\catcode`\"=\active\def"{\~N}        
\catcode`\Í=\active\defÍ{\~O}        
\catcode`\‰=\active\def‰{\^a}        
\catcode`\=\active\def{\^e}        
\catcode`\"=\active\def"{\^{\i}}     
\catcode`\™=\active\def™{\^o}        
\catcode`\ž=\active\defž{\^u}        

\let\optionkeymacros\null

\newtheorem{Theoreme}{Th\'eor\`eme}
\newtheorem{Lemme}{Lemme}
\newtheorem{Proposition}{Proposition}
\newtheorem{Corollaire}{Corollaire}

\newcommand{\NN}{\mathbb N}
\newcommand{\PP}{\mathbb P}
\newcommand{\ZZ}{\mathbb Z}
\newcommand{\QQ}{\mathbb Q}
\newcommand{\RR}{\mathbb R}
\newcommand{\CC}{\mathbb C}
\newcommand{\GG}{\mathbb G}
\newcommand{\FF}{\mathbb F}
\newcommand{\sqm}[4]{\displaystyle{
\left({#1 \atop #3}{#2 \atop #4}\right)}}
\newcommand{\binomial}[2] {{\binom{#1}{#2}}} 
\def\tr{{\bf t}}
\def\no{{\bf n}}

\def\und{\underline}
\def\SL{{\bf SL}}
\def\GL{{\bf GL}}
\def\aa{a}
\def\SP{{\bf SP}}

\begin{document}

\begin{center}

\LARGE{Introduction aux formes \\modulaires de Hilbert et ˆ leur\\  propriŽtŽs diffŽrentielles.}\\

\vspace{15pt}

\Large{Federico Pellarin}

\end{center}

\vspace{15pt}

\tableofcontents

\section{Introduction.}

Pour expliquer les motivations de cette sŽrie de six exposŽs, nous dŽcrivons brivement
un thŽorme de Nesterenko dans \cite{Nesterenko:Introduction} (voir le texte \cite{Bosser:Independance} de Bosser dans ce volume).

Posons $q=e^{2\pi{\rm i}\tau}$ avec $\tau$ un nombre complexe de partie imaginaire positive
(o ${\rm i}$ dŽsigne le nombre complexe $\sqrt{-1}$), soient $E_2(q),E_4(q),E_6(q)$ les s\'eries
d'Eisenstein usuelles, de poids $2,4,6$ respectivement. 

\begin{Theoreme} Si $q\not=0$ alors le sous-corps 
de $\CC$ engendr\'e par $q,E_2(q),E_4(q),\\ E_6(q)$ a un degr\'e de transcendance sur $\QQ$
au moins $3$.\label{lemme:nesterenko}\end{Theoreme}

La preuve de Nesterenko s'appuye sur la propri\'et\'e suivante des s\'eries $E_2,E_4,E_6$.
L'op\'erateur de d\'erivation $D=\displaystyle{q\frac{d}{dq}}$ fait de l'anneau $\QQ[E_2,E_4,E_6]$
un anneau {\em dif\-f\'erentiel}. Plus en particulier, nous avons~:
\begin{eqnarray}
DE_2 & = & \frac{1}{12}(E_2^2-E_4)\nonumber\\
DE_4 & = & \frac{1}{3}(E_2E_4-E_6)\label{eq:systeme}\\
DE_6 & = & \frac{1}{2}(E_2E_6-E_4^2)\nonumber.
\end{eqnarray}
On remarque que~:
\begin{eqnarray}
\QQ[E_2,E_4,E_6] &= &\QQ[E_2,DE_2,D^2E_2]\label{eq:monogene}.
\end{eqnarray}

La forme $E_2(e^{2\pi{\rm i}\tau})$, quasi-modulaire de poids $2$, contient
donc toutes les informations qui permettent d'ac\-cŽder
ˆ une preuve du thŽorme.

Dans ces exposŽs, nous voulons faire une introduction ˆ des problmes 
analogues, mais en deux ou plusieurs variables complexes. 

\subsection{Motivation et structure de ce texte.}
 
Nous commenons par dŽcrire les groupes modulaires de Hilbert et leur action sur les
produits de copies de demi-plans supŽrieurs complexes. 
Nous faisons ensuite une courte introduction aux formes modulaires de Hilbert, qui sont des 
formes modulaires de plusieurs variables complexes, associŽes ˆ des corps de nombres totalement rŽels~: le nombre de variables 
est Žgal au degrŽ du corps.
Nous construisons des formes modulaires de Hilbert~: des sŽries d'Eisenstein, et des sŽries thta. 

Nous regardons ensuite les zŽros de certaines formes modulaires de Hilbert, par analogie avec 
le discriminant elliptique $\Delta=(E_4^3-E_6^2)/1728$ (de
Jacobi), qui est sans zŽros pour $|q|<1$. La fonction
$E_2$ est Žgale ˆ la dŽrivŽe logarithmique $D\Delta/\Delta$.
Gr‰ce ˆ (\ref{eq:monogene}), ceci permet de retrouver le systme (\ref{eq:systeme}).

En deux ou plusieurs variables, nous
observons qu'il n'existe pas de  forme modulaire de Hilbert sans zŽros dans son domaine de dŽfinition, 
ce qui constitue une difficultŽ 
dans la construction d'une gŽnŽralisation du systme (\ref{eq:systeme}).

De plus, en deux ou plusieurs variables, il existe une division importante entre 
formes modulaires dites de poids {\em parallle} et celles dites de poids {\em non parallle}, qui
n'existe pas en une variable complexe.

Il s'avre que des gŽnŽralisations partielles de la forme modulaire $\Delta$ de Jacobi en 
plusieurs variables existent. Nous donnons un exemple
explicite en deux variables complexes~: la forme modulaire ainsi construite permettra de calculer explicitement
un certain anneau de formes modulaires de Hilbert de deux variables complexes.

En gŽnŽral, les structures d'anneaux (graduŽs) des formes modulaires de Hilbert de poids
parallle peuvent tre dŽterminŽes de
plusieurs faons. La mŽthode la plus ancienne consiste ˆ Žtudier des sŽries d'Eisenstein de petit poids
(le plus souvent de poids $(1,\ldots,1)$), pour des sous-groupes de congruence, tordues par des caractres
de Dirichlet, ou des sŽries thtas~: c'est la mŽthode introduite par Hecke, 
avec des contributions de Kloosterman, Gštzky, Gundlach, Maass, Resnikoff.
Les thŽormes de structure les plus prŽcis concernent uniquement 
le cas des formes modulaires de Hilbert de deux variables
complexes.

Cette mŽthode a ŽtŽ remplacŽe ou complŽtŽe dans 
plusieurs cas par l'utilisation de la gŽomŽtrie des surfaces lorsque Hirzebruch
a dŽcouvert comment construire explicitement un modle dŽsingularisŽ de compactifications de surfaces modulaires de Hilbert.
On a alors pu utiliser la thŽorie de l'intersection dans les surfaces, et dŽterminer explicitement 
des bases d'espaces vectoriels de formes modulaires (i.e. sections de faisceaux inversibles).
Le livre \cite{Geer:Hilbert} contient les fondations de cette technique.

Dans ces exposŽs nous proposons une mŽthode encore diffŽrente, qui consiste
ˆ exploiter les propriŽtŽs diffŽrentielles des formes modulaires (Jacobi, Rankin, Resnikoff).
En Žtudiant ces propriŽtŽs, nous dŽcrivons compl\-tement 
la structure de l'anneau des formes modulaires de Hilbert de poids parallle, associŽes au corps $\QQ(\sqrt{5})$.
Aprs avoir construit deux formes modulaires $\varphi_2,\chi_5$ (de poids $(2,2)$ et $(5,5)$), nous montrons que
l'anneau des formes modulaires de Hilbert de poids parallle 
\[{\cal T}(\Gamma)=\bigoplus_{r\in\NN}M_{r\und{1}}(\Gamma)\]
pour le corps $\QQ(\sqrt{5})$
est un anneau de type fini, dont on peut dŽterminer
explicitement la structure, et est engendrŽ par les images de $\varphi_2,\chi_5$ par certains opŽrateurs
diffŽrentiels que nous dŽcrirons.

Par contre, nous montrons que l'anneau des formes modulaires de Hilbert
(de poids quelconque) associŽ ˆ un corps quadratique rŽel $K$ quelconque
\[{\cal L}(\Gamma)=\bigoplus_{\und{r}\in\NN^2}M_{\und{r}}(\Gamma)\]
n'est pas un anneau de type fini, 
et ne peut mme pas tre engendrŽ par les images d'un ensemble
fini de formes  modulaires, par des opŽrateurs diffŽren\-tiels.

Suivant Resnikoff, nous Žtudions le corps diffŽrentiel engendrŽ par les dŽrivŽes partielles
d'une seule forme modulaire de Hilbert. Dans le cas de $K=\QQ(\sqrt{5})$, on verra que l'anneau 
engendrŽ par toutes les dŽrivŽes partielles de la forme modulaire $\varphi_2$ est contenu dans un anneau de type fini,
que l'on dŽcrira explicitement.

\subsection{PrŽlude~: le cas elliptique.\label{section:prelude}}
Une forme modulaire {\em elliptique} est par dŽfinition une forme modulaire pour le groupe $\SL_2(\ZZ)$ (cf. \cite{Serre:Modulaires} p.
135 ou \cite{Martin:Transcendance}).
\begin{Proposition}
Soit $F$ une forme modulaire elliptique non constante de poids $f$. Alors~:
\begin{enumerate}
\item La fonction $fF(d^2F/dz^2)-(f+1)(dF/dz)^2$ est une forme modulaire de poids $2f+4$.
\item Les fonctions $F,(dF/dz),(d^2F/dz^2),(d^3F/dz^3)$ sont algŽbriquement dŽ\-pendantes sur $\CC$.
\item Les fonctions $F,(dF/dz),(d^2F/dz^2)$ sont algŽbriquement indŽpendantes sur $\QQ(e^{2\pi{\rm i}z})$.
\end{enumerate}
\label{proposition:rankin_mahler}\end{Proposition}
Ceci est la proposition 1.1 p. 2 de \cite{Nesterenko:Introduction}~: ce rŽsultat est dž ˆ Mahler, mais 
s'appuie fortement sur des contributions de Hurwitz et Rankin. L'opŽrateur
\[F\mapsto (2\pi{\rm i})^{-2}\left(fF\frac{d^2F}{dz^2}-(f+1)\left(\frac{dF}{dz}\right)^2\right)\] est une spŽcialisation des {crochets de Rankin-Cohen} ainsi
dŽfinis.  Pour deux formes modulaires
$F,G$ de poids respectifs $f,g$~:
\[[F,G]_{n}:=\frac{1}{(2\pi{\rm i})^{n}}\sum_{r=0}^n(-1)^r\binomial{f+n-1}{n-r}\binomial{g+n-1}{r}F^{(r)}G^{(n-r)},\] o $F^{(r)}$ dŽsigne 
la dŽrivŽe $r$-ime de $F$ par rapport ˆ la variable complexe $z$.
La forme modulaire $[F,G]_{2n}$ est de poids $f+g+2n$ si elle est non nulle. Par exemple, on a~:
\[[F,F]_2 = (f+1)(2\pi{\rm i})^{-2}\left(fF\frac{d^2F}{dz^2}-(f+1)\left(\frac{dF}{dz}\right)^2\right).\] Les crochets de Rankin-Cohen 
sont ŽtudiŽs en dŽtail dans 
\cite{Martin:Transcendance} et \cite{Zagier:Indian}.

Il existe des rŽlations liant ces crochets, lorsqu'ils sont appliquŽs ˆ une mme forme modulaire. En voici un exemple~:
\[[F,[F,F]_2]_2=\frac{6f(f+1)}{(f+2)(f+3)}F[F,F]_4.\] D'autres exemples peuvent tre construits avec la 
remarque du paragraphe \ref{section:differentiels}.

Le plus simple des crochets est le crochet de Rankin \[[F,G]_1=fF\frac{dG}{dz}-gG\frac{dF}{dz}.\] Ce crochet est antisymŽtrique (\footnote{Plus gŽnŽralement,
les crochets $[\cdot,\cdot]_{2n+1}$ sont
antisymŽtriques, et les crochets $[\cdot,\cdot]_{2n}$ sont symŽ\-triques.}), et il peut
tre utilisŽ pour dŽterminer la structure d'anneau graduŽ de toutes les formes modulaires elliptiques~: 
cet anneau est isomorphe ˆ l'anneau de polyn™mes $\CC[X_0,X_1]$, en deux indŽterminŽes $X_0,X_1$.

On commence par une sŽrie d'Eisenstein
\[E_4^*(z)=\sum_{(m,n)\not=(0,0)}\frac{1}{(mz+n)^4},\]
qui est une forme modulaire elliptique non nulle de poids $4$. La proposition \ref{proposition:rankin_mahler} points 1, 2, et le fait que l'opŽrateur
diffŽrentiel
$F\mapsto[F,F]_2$ est d'ordre $2$, implique que $\Delta^*=[E_4^*,E_4^*]_2$ est une forme modulaire non nulle de poids $12$.
De plus, $E_4^*,\Delta^*$ sont algŽbriquement indŽpendantes.

En utilisant le systme (\ref{eq:systeme}), on trouve $D\Delta^*=E_2\Delta^*$. Ainsi, $\Delta^*(e^{2\pi{\rm i}z})$
ne s'annule pas dans 
${\cal H}=\{z\in\CC\mbox{ tel que }\Im(z)>0\}$ (on peut Žgalement
appliquer la formule p. 143 de \cite{Serre:Modulaires} (voir aussi \cite{Martin:Transcendance}).

On calcule ensuite $[E_4^*,\Delta^*]_1$. On voit facilement que cette forme modulaire est non nulle~: en effet, si cette
forme Žtait nulle on aurait l'existence d'une constante $c\in\CC^\times$ et d'entiers $x,y$ non nuls, avec $E_4^*{}^x=c\Delta^*{}^y$.
Or, $E_4^*$ est non parabolique, et $\Delta^*$ est parabolique~: donc $[E_4^*,\Delta^*]_1\not=0$.

La forme modulaire $[E_4^*,\Delta^*]_1$ est parabolique de poids $18$. Ceci veut dire que la fonction
$E_6^*:=[E_4^*,\Delta^*]_1/\Delta^*$ est une forme modulaire non nulle de poids $6$, car $\Delta^*(e^{2\pi{\rm i}z})$ ne s'annule pas
dans ${\cal H}$. Nous devons maintenant prouver que
$E_4^*,E_6^*$ sont algŽbriquement indŽpendantes, mais ceci ne dŽcoule pas de la proposition \ref{proposition:rankin_mahler}.

Pour ceci, on note qu'il existe une combinaison linŽaire $\lambda E_4^*{}^3+\mu E_6^*{}^2$ qui est une forme parabolique. 
C'est donc un multiple $c\Delta^*$. Si cette combinaison linŽaire n'est pas triviale, la constante de proportionnalitŽ
$c$ est non nulle et $E_4^*,E_6^*$ sont algŽbriquement indŽpendantes car $E_4^*$ et $\Delta^*$ le sont,
par la proposition \ref{proposition:rankin_mahler}.
On doit calculer quelques coefficients de Fourier. On note $c_4,c_6,c_{12}$ des nombres complexes non nuls
tels que les sŽries de Fourier de $E_4(q)=c_4E^*_4(z),E_6(q)=c_6E^*_6(z),\Delta(q)=c_{12}\Delta^*(z)$ 
(avec $q=e^{2\pi{\rm i}z}$) soient telles que le coefficient du terme
de plus petit degrŽ en $q$ soit Žgal ˆ $1$.
On a donc~:
\begin{eqnarray*}
E_4 & = & 1+240q+\cdots\\
\Delta & = & q+\cdots\\
\Delta^{-1}[E_4,\Delta]_1 & = & (4(1+\cdots)(q+\cdots)-12(q+\cdots)(240q+\cdots))/(q+\cdots)\\
& = & 4(1-504q+\cdots),
\end{eqnarray*}
d'o $E_6=1-504q+\cdots$. On calcule~:
\begin{eqnarray*}
[E_4,E_6]_1 & = & 4(1+\cdots)(-504q+\cdots)-6(1+\cdots)(240q+\cdots)\\
& = & -2\cdot 1728\Delta.
\end{eqnarray*}
Ceci implique qu'il n'existe pas de relation $E_4^n=\lambda E_6^m$ avec $\lambda\in\CC,n,m\in\ZZ$,
et la constante de proportionnalitŽ $c$ est non nulle~: donc les formes modulaires $E_4,E_6$ sont algŽbriquement indŽpendantes
(naturellement, on peut aussi voir ceci en calculant quelques autres coefficients de Fourier, mais notre
souci sera par la suite de gŽnŽraliser {\em cette} dŽmonstration).

On montre maintenant que toute forme modulaire $F$ s'Žcrit de manire unique comme {\em polyn™me isobare} 
(\footnote{Un polyn™me isobare est un polyn™me en des formes modulaires qui est une forme modulaire.})
de $\CC[E_4,E_6]$~: on fait une rŽcurrence sur le poids $k$ de $F$. Si $F$ est parabolique, $F/\Delta$ est une forme modulaire
de poids $k-12$. Si $F$ n'est pas parabolique et est de poids $k\geq 4$ alors il existe un polyn™me isobare $P\in\CC[E_4,E_6]$
de mme poids que $F$, tel que $F-P$ soit une forme parabolique, et on se ramne au cas prŽcŽdent.

Il reste ˆ traiter le cas o $F$ est de poids $0$ ou $2$. Si $F$ est de poids $0$, alors $F$ est une constante (cf. \cite{Serre:Modulaires} p. 143).
Si $F$ a poids $2$ et est non nulle, alors elle ne peut pas tre une forme parabolique (il n'existe pas de 
formes modulaires de poids nŽgatif~: cf. \cite{Serre:Modulaires} p. 143). Il faut de plus que $F^2=\lambda_1E_4$ et $F^3=\lambda_2E_6$, avec
$\lambda_1\lambda_2\not=0$ (car sinon, on trouve des formes modulaires de poids nŽgatif qui doivent tre nulles), mais ceci entra"ne une relation 
de dŽpendance algŽbrique entre $E_4$ et $E_6$, qui ne peut pas exister~: il n'existe pas de forme modulaire non nulle
de poids $2$.

Les autres dŽtails que nous avons omis dans cette dŽmonstration, peuvent tre facilement trouvŽs par le lecteur.
Nous avons commencŽ ce raisonnement avec la sŽrie d'Eisenstein $E_4^*$~: on peut voir qu'on peut commencer Žgalement avec $E_6^*$
ou $\Delta^*$. 

Dans la suite, nous voulons gŽnŽraliser ces arguments aux formes modulaires de Hilbert.

\section{Groupes modulaires de Hilbert.}

Nous faisons ici une courte introduction aux groupes modulaires de Hilbert.
Le groupe Bih$({\cal H}^n)$ des automorphismes
biholomorphes de ${\cal H}^n$ est l'ex\-tension de sa composante neutre Bih${}_0({\cal H}^n)$ par le groupe fini des permutations
des  coordonnŽes $z_i$. D'autre part, Bih${}_0({\cal H}^n)$ est le produit direct de $n$ copies de Bih$({\cal H}^1)\cong$PSL${}_2(\RR)$~:
\[1\rightarrow\mbox{PSL}{}_2(\RR)^n\rightarrow\mbox{Bih}({\cal H}^n)\rightarrow{\mathfrak S}_n\rightarrow 1,\] o ${\mathfrak S}_n$ dŽsigne le groupe
symŽtrique agissant sur un ensemble fini ˆ $n$ ŽlŽments.
Le groupe Bih${}_0({\cal H}^n)$ agit sur ${\cal H}^n$ de la manire suivante. 
Posons $\und{z}=(z_1,\ldots,\\ z_n)$, soit $\gamma\in$Bih${}_0({\cal H}^n)$~:
\[\gamma=\left(\sqm{a_1}{b_1}{c_1}{d_1},\ldots,\sqm{a_n}{b_n}{c_n}{d_n}\right).\] Alors~:
\begin{eqnarray*}
(\gamma,\und{z})\in\mbox{Bih}_0({\cal
H}^n)\times{\cal H}^n & \mapsto \gamma(\und{z}):=\displaystyle{\left(\frac{a_1z_1+b_1}{c_1z_1+d_1},\ldots,
\frac{a_nz_n+b_n}{c_nz_n+d_n}\right)}\in{\cal H}^n.
\end{eqnarray*}
Soit $K$ un corps totalement rŽel de degrŽ $n$, d'anneau d'entiers ${\cal O}_K$. Nous avons $n$ plongements~:
\[\sigma_i:K\rightarrow\RR.\] Posons $\Sigma=(\sigma_i)_{i=1,\ldots,n}$.
Le {\em groupe modulaire de Hilbert} (associŽ ˆ ${\cal O}_K$) est le sous-groupe 
$\Gamma_K=\SL_2({\cal O}_K)$
de $\SL_2(K)$~:
\begin{eqnarray*}\Gamma_K & = &\left\{\sqm{a}{b}{c}{d}\mbox{ avec }a,b,c,d\in{\cal O}_K\mbox{ et
}ad-bc=1\right\}.\end{eqnarray*} 
Le groupe $\Gamma_K$ se plonge dans $\SL_2(\RR)^n$ par le biais du plongement $\Sigma$.
Si $n=1$, on retrouve le groupe modulaire $\SL_2(\ZZ)$ usuel. Si $n=2$, nous noterons $\sigma_1(\mu)=\mu$ et $\sigma_2(\mu)=\mu'$.

Si ${\mathfrak A}$ est un idŽal non nul de ${\cal O}_K$, le noyau $\Gamma_K({\mathfrak A})$ 
de l'homomorphisme $\Gamma_K\rightarrow\SL_2({\cal O}_K/{\mathfrak A})$
est un sous-groupe normal d'indice fini de $\Gamma_K$. On l'appelle parfois {\em sous-groupe de congruence principal}
associŽ ˆ ${\mathfrak A}$.

Dans toute la suite, lorsque le corps de nombres $K$ sera clairement dŽter\-minŽ, nous Žcrirons $\Gamma=\Gamma_K$ et $\Gamma({\mathfrak A})=
\Gamma_K({\mathfrak A})$.

\subsection{PropriŽtŽs de base.\label{section:sans_demonstration}}

Comme dans le cas $n=1$, l'action de $\Gamma$ (ou de $\Gamma({\mathfrak A})$) possde des bonnes propriŽtŽs topologiques,
que nous ne vŽrifions pas ici (voir \cite{Freitag:Hilbert}, chapitre 1).

On montre ainsi que $\Gamma$ agit {\em proprement} sur ${\cal H}^n$
(c'est-ˆ-dire que si $A$ est un compact de ${\cal H}^n$, alors $\gamma(A)\cap A\not=\emptyset$ pour
au plus une finitude de $\gamma\in\Gamma$). 

Pour montrer cette propriŽtŽ, on utilise le fait que $\Gamma$ est {\em discret} dans $\SL_2(\RR)^n$ (\footnote{Un
sous-ensemble
$E\subset\RR^N$ est discret si pour tout compact $U$ de $\RR^N$, $U\cap E$ est fini.}),
ce qui est garanti par le fait que ${\cal O}_K$ est ˆ son tour discret dans $\RR^n$ via le
plongement $\Sigma$ (voir \cite{Freitag:Hilbert} proposition 1.2 p. 7 et proposition 2.1 p. 21).

\medskip

\noindent {\bf Exemple pour $n=2$}.
Si $K=\QQ(\sqrt{5})$, le groupe $\Gamma_K$ est engendrŽ par les
trois matrices~:
\[S=\sqm{1}{1}{0}{1},T=\sqm{0}{-1}{1}{0},U=\sqm{\epsilon}{0}{0}{\epsilon'}\]
avec $\epsilon=\displaystyle{\frac{1+\sqrt{5}}{2}}$,
qui satisfont
\begin{equation}U^{-1}=TUT,S^{-1}=TSTST.\label{eq:relations_goetzky}\end{equation}
Pour vŽrifier ceci, on remarque tout d'abord que ${\cal O}_K=\ZZ[\epsilon]$, et donc que
le sous-groupe $\displaystyle{\sqm{*}{*}{0}{*}}\subset\Gamma$ est engendrŽ par $S,U$. 

Soit $\displaystyle{\sqm{\alpha}{\beta}{\gamma}{\delta}}\in\Gamma$, avec $\gamma\not=0$. Pour tout $\nu\in {\cal O}_K$,
\[\sqm{\alpha}{\beta}{\gamma}{\delta}=\sqm{\alpha\nu-\beta}{\;\;\alpha}{\gamma\nu-\delta}{\;\;\gamma}\cdot T\cdot\sqm{1}{\nu}{0}{1}.\]
Nous pouvons choisir $\nu$ avec $|\no(\gamma\nu-\delta)|<|\no(\gamma)|$, o $\no(\gamma)$ dŽsigne la norme
absolue de $\gamma$, et appliquer une hypothse de 
rŽcurrence sur l'entier $|\no(\gamma)|$ (voir \cite{Goetzky:Anwendung} pp. 414-416).

\medskip

\noindent {\bf Notation}. Dans ces exposŽs, nous faisons jouer ˆ $K=\QQ(\sqrt{5})$ un r™le particulier. Dans ce cas, et uniquement
dans ce cas, nous Žcrirons $\Upsilon=\SL_2({\cal O}_K)$.

\medskip

Soit $\Omega$ un sous-groupe discret de $\SL_2(\RR)^n$.
L'action de $\Omega$ s'Žtend de manire naturelle en une action sur $\overline{{\cal H}}^n$, o
\[\overline{{\cal H}}:={\cal H}\cup\PP_1(\RR).\]

\noindent {\bf DŽfinition.} On dit que $\Omega$ a une {\em pointe} en $\infty:=({\rm i}\infty,\ldots,{\rm i}\infty)
\in\overline{{\cal H}}^n$ si $\Omega$ contient un sous-groupe $\Omega_\infty$~:
\[\Omega_\infty=\left\{\left(\sqm{\alpha_1}{\nu_1}{0}{\alpha_1^{-1}},\ldots,
\sqm{\alpha_n}{\nu_n}{0}{\alpha_n^{-1}}\right)\mbox{ avec }\alpha_i\in\RR_{>0}\mbox{ et }\nu_i\in\RR\right\},\]
ayant la propriŽtŽ que l'ensemble $\{(\nu_1,\ldots,\nu_n)\}\subset\RR^n$ est un rŽseau, et
l'ensem\-ble $\{(\alpha_1,\ldots,\alpha_n)\}\subset\RR_{>0}^n$ est un sous-groupe multiplicatif libre de rang $n-1$ sur
$\ZZ$. Par exemple, $\infty\in\overline{{\cal H}}$ est une pointe de $\Omega=\SL_2(\ZZ)$, et $\Omega_\infty$ est le groupe
cyclique engendrŽ par toutes les puissances de la matrice $\sqm{1}{1}{0}{1}$.

Soit $a\in\overline{{\cal H}}^n$. On dit que $\Omega$ a une {\em pointe} en $a$ s'il existe un ŽlŽment
$A\in\SL_2(\RR)^n$ tel que $A(a)=\infty$ et tel que le groupe $A\cdot\Omega\cdot A^{-1}$ a une pointe en $\infty$.

\medskip

Le sous-groupe $\Omega_\infty$ s'appelle le {\em stabilisateur} de $\infty$.
Pour $K$ corps de nombres totalement rŽel de degrŽ $n$, 
$\Gamma_\infty$ est isomorphe ˆ un produit semi-direct de ${\cal O}_K\cong\ZZ^n$ par $({\cal O}_K^\times)^2
\cong\ZZ^{n-1}$~: $\infty$ est toujours une pointe des groupes modulaires de Hilbert $\Gamma$.
Par exemple, Le groupe $\Upsilon_\infty$ est engendrŽ par $S,U$.

\medskip

On peut Žtendre ces dŽfinitions ˆ un quelconque sous-groupe de $\Gamma$ d'indice fini, par exemple $\Gamma({\mathfrak A})$, pour un idŽal
non nul ${\mathfrak A}$ de ${\cal O}_K$. Plus l'indice est grand, plus on trouve dans $\overline{{\cal H}}^n$ un grand nombre
de  pointes non Žquivalentes modulo l'action de $\Gamma({\mathfrak A})$.
Par un calcul ŽlŽmentaire, on voit que les pointes de $\Gamma$ se situent dans $\PP_1(\RR)$.

Le groupe $\Gamma$ peut possŽder plusieurs pointes non Žquivalentes. 

\begin{Lemme}
L'ensemble des classes d'Žquivalen\-ce des
pointes de $\overline{\cal H}^n$ sous l'action de $\Gamma$ 
est en bijection avec les classes d'idŽaux modulo les idŽaux principaux de $K$.\label{lemme:Maass}
\end{Lemme}

\noindent {\bf DŽmonstration.}
Tout d'abord, montrons que les ŽlŽ\-ments de $\PP_1(K)$ sont des pointes. 

Nous avons dŽjˆ vu que $\infty$ est une pointe.
Soit $a\in K=\PP_1(K)-\{\infty\}$ et considŽrons \[A:=\displaystyle{\sqm{0}{1}{-1}{a}}\in\SL_2(K).\]
On a $A(a)=\infty$, donc le groupe~:
\[A\cdot\Gamma\cdot A^{-1}\]
possde une pointe en $\infty$, et $\Gamma$ a une pointe en $a$. Nous venons de montrer que
$\PP_1(K)\subset\{\mbox{pointes}\}$.

RŽciproquement, soit $p=(p_1,\ldots,p_n)\in\PP_1(\RR)$ une pointe. Par dŽfinition il existe une matrice $P$ de $\SL_2(\RR)$
telle que $P(p)=\infty$ et telle que le groupe $P\cdot\Gamma\cdot P^{-1}$ contient un
groupe de translations isomorphe ˆ ${\cal O}_K$. 
Toutes ces translations sont des ŽlŽments {\em para\-boliques}, i.e.
de trace rŽduite $\pm 2$. Soit $t\not=\displaystyle{\sqm{1}{0}{0}{1}}$ une de ces translations.
En vertu de la condition sur la trace, les Žquations~:
\[\sigma_i(t)(p_i)=p_i\mbox{ pour }i=1\ldots,n\]
ont une unique solution $p$, et de plus, $p\in\Sigma(K)$, donc $\PP_1(K)\supset\{\mbox{pointes}\}$.

Le fait que $\PP_1(K)/\Gamma$ soit en bijection avec le groupe des classes d'idŽaux de $K$ est bien connu et se voit
en sachant que
tout idŽal fractionnaire non principal ${\cal I}$ de $K$ est engendrŽ par deux ŽlŽments $a,b\in K$ distincts.

\subsection{Une description de $X_\Gamma$.\label{section:description}}

Un {\em ensemble fondamental} $E\subset{\cal H}^2$ pour $\Gamma$ est un sous-ensemble 
tel que pour tout $\und{z}\in{\cal H}^2$ il existe $\gamma\in\Gamma$ avec $\gamma(\und{z})\in E$.

Un {\em domaine fondamental} $D\subset{\cal H}^2$ pour $\Gamma$ est un ensemble fondamental mesurable
pour la mesure de Lebesgue, tel que le sous-ensemble dont les ŽlŽments sont les 
$\und{z}\in D$ tels que $\gamma(\und{z})\in D$ pour quelques
$\gamma\in\Gamma-\left\{\sqm{1}{0}{0}{1}\right\}$, ait une mesure nulle.

Le groupe $\Gamma$ Žtant discret, tout ensemble fondamental qui est mesurable contient un
domaine fondamental (voir \cite{Freitag:Hilbert}). 

\medskip

\noindent {\bf 1}. L'espace $X_\Gamma:=({\cal H}^n\cup\PP_1(K))/\Gamma$, muni d'une certaine topologie
naturelle engendrŽe par des systmes de voisinages de points de ${\cal H}^n$ et de pointes, 
est compact. Nous ne ne donnerons pas de dŽmonstra\-tion gŽnŽrale de ce
fait, pour laquelle nous renvoyons ˆ \cite{Freitag:Hilbert}, chapitre 1, ou ˆ \cite{Geer:Hilbert}, chapitre 1.
Dans le paragraphe
\ref{section:goetzky} nous construirons explicitement un domaine fondamental dans le cas $K=\QQ(\sqrt{5})$, et cette
propriŽtŽ de compacitŽ peut tre vue facilement ˆ partir de cette construction.

\medskip

\noindent {\bf 2.} En gŽnŽral, l'espace $X_\Gamma$ possde des singularitŽs. Ces singularitŽs ont leur
support dans l'ensemble des poin\-tes
et des points {\em elliptiques} (points fixes de transformations d'ordre fini). Voici un exemple explicite 
de point elliptique dans le cas $K=\QQ(\sqrt{5})$~: le point
$\und{\zeta}:=(\zeta_5,\zeta_5^2)\in{\cal H}^2$, avec $\zeta_5=e^{2\pi i/5}$. 

Ce point est un point elliptique d'ordre $5$ pour $\Upsilon$, car la matrice $\phi=\displaystyle{\sqm{-\epsilon'}{-1}{1}{0}}$, qui est 
d'ordre $5$, satisfait $\phi(\und{\zeta})=\und{\zeta}$.

\medskip

\noindent {\bf 3.} Si $n>1$, {\em toutes} les pointes dŽterminent des singularitŽs. Une dŽmonstra\-tion homologique
de ce fait se trouve dans \cite{Freitag:Hilbert} p. 30. Les points elliptiques dŽterminent parfois
des singularitŽs~: on peut dŽmontrer que l'orbite du point $\und{\zeta}$ dŽfini ci-dessus est un point rŽgulier
dans $X_\Upsilon$.

Pour simplifier l'exposition, {\em nous supposons dorŽnavant que l'anneau des entiers de $K$ est principal}. 
C'est le cas si $K=\QQ$ ou $K=\QQ(\sqrt{5})$.

\section{Formes et fonctions modulaires de Hilbert.}

Il existe des formes modulaires pour $\SL_2(\ZZ)$ sans zŽros 
dans ${\cal H}$~: la forme $\Delta$ en est un exemple. Nous introduisons ici les formes modulaires de Hilbert associŽes ˆ un
corps de nombres totalement rŽel de degrŽ $n$, et nous faisons une Žtude analytique~: il en rŽsultera qu'il n'existe
pas de forme modulaire de Hilbert sans zŽros dans ${\cal H}^n$.

\medskip

Fixons un corps totalement rŽel $K$.
Soit $\Omega$ un sous-groupe d'indice fini de $\Gamma:=\SL_2({\cal O}_K)$ (dans la plupart des exemples, $\Omega$ est un sous-groupe
de congruence principal), soit $\xi:\Omega\rightarrow\CC^\times$ un caractre.
Nous dŽfinissons les formes modulaires de Hilbert pour $\Omega$ avec caractre, de deux manires diffŽrentes, suivant que $K=\QQ$ ou $K\not=\QQ$.

Si $K=\QQ$, une {\em forme modulaire de Hilbert $f$ de poids $n\in\ZZ$ pour $\Omega$ avec caractre $\xi$} est par dŽfinition une forme modulaire 
pour le groupe $\Omega\subset\SL_2(\ZZ)$ de poids $n$, avec caractre $\xi$.
En particulier, toute forme modulaire de Hilbert pour $K=\QQ$ est
holomorphe aux pointes de $\Omega$.

\medskip

\noindent {\bf DŽfinition}.
Si $K\not=\QQ$,
une {\em forme modulaire de Hilbert de poids $\und{r}=(r_1,\ldots,\\ r_n)\in\ZZ^n$ avec caractre $\xi$ pour $\Omega$}
est une fonction holomorphe $f:{\cal H}^n\rightarrow\CC$ satisfaisant~:
\begin{equation}
f(\gamma(\und{z}))=\xi(\gamma)\prod_{i=1}^n(\sigma_i(c)z_i+\sigma_i(d))^{r_i}f(\und{z}),\label{eq:f_modulaire}\end{equation}
pour tous $\und{z}=(z_1,\ldots,z_n)$ et $\gamma=\displaystyle{\sqm{a}{b}{c}{d}}\in\Omega$.

Si $\xi=1$, on parle plus simplement de formes modulaires de Hilbert (sans mentionner \og avec caractre trivial\fg).

\medskip

\noindent {\bf DŽfinition}. Une {\em fonction modulaire de Hilbert} de poids $\und{r}=(r_1,\ldots,r_n)\in\ZZ^n$ pour $\Omega$
est une fonction {\em mŽromorphe} $f:{\cal H}^n\rightarrow\CC$ satisfaisant (\ref{eq:f_modulaire})
pour tous $\und{z}=(z_1,\ldots,z_n)$ et $\gamma=\displaystyle{\sqm{a}{b}{c}{d}}\in\Omega$, avec $\xi=1$.

\subsection{ Fonctions rŽgulires aux pointes.}

Notons
$\tr:K\rightarrow\QQ$ la trace de $K$ sur $\QQ$. Pour simplifier les notations, nous Žtendons cette fonction de la manire suivante.
Si $\und{z}=(z_1,\ldots,z_n)\in\CC^n$, nous Žcrivons $\tr(\und{z})=z_1+\cdots+z_n$~;
ainsi par exemple, si $\nu\in K$, nous Žcrivons $\tr(\nu\und{z})$ plut™t que $\sigma_1(\nu) z_1+\cdots+\sigma_n(\nu)z_n$. 

\medskip

Les conditions analytiques imposŽes dans la dŽfinition d'une forme mo\-dulaire de Hilbert 
pour $K\not=\QQ$ sont plus faibles que pour $K=\QQ$.
On {\em ne demande pas} que les formes modulaires de Hilbert soient holomorphes ˆ l'infini. Nous expliquons ici ce
phŽnomne.

Soit $K$ un corps de nombres totalement rŽel, de degrŽ $n>1$. Le thŽorme des unitŽs de Dirichlet implique~:
\[{\cal O}_K^\times\cong\ZZ^{n-1}\times\frac{\ZZ}{2\ZZ}.\] Soit $\und{c}=(c_2,\ldots,c_{n})\in\ZZ^{n-1}$.

\medskip

\noindent {\bf DŽfinition.}
Soit $M$ un $\ZZ$-module complet de $K$,
soit $V$ un sous-groupe d'indice fini de ${\cal O}_K^\times$ tel que $\eta M=M$ pour tout $\eta\in V$.
L'{\em espace vectoriel ${\cal A}_{\und{c}}(V,M)$ des fonctions rŽ\-gu\-lires} de poids
$\und{c}$ au voisinage de $\infty$ est
constituŽ des sŽries de Fourier ˆ coefficients 
complexes, absolument convergentes dans un sous-ensemble non vide de ${\cal H}^n$~:
\begin{equation}f(\und{z}):=\sum_{\nu\in M^*}d_\nu\exp\{2\pi i\tr(\und{z}\nu)\},\label{eq:fourier}\end{equation}
o $M^*$ dŽsigne le dual de $M$ pour la trace $\tr$, avec la condition que
\[d_{\nu\epsilon}=d_{\nu}\prod_{i=2}^n\sigma_i(\epsilon)^{c_i},\]
pour tout $\epsilon\in V^2$ et $\nu\in M^*$.

\medskip
 
Si $f\in{\cal A}_{\und{c}}(V,M)$ a la sŽrie de Fourier (\ref{eq:fourier}), 
alors $d_\nu\not=0$ implique $\sigma_i(\nu)\geq 0$
pour tout $i=1,\ldots,n$, d'o le fait que $f$ converge au voisinage de la pointe $\infty$. 
La vŽrification est laissŽe en exercice~: il faut utiliser l'exis\-tence 
d'unitŽs $\epsilon\in{\cal O}_K^\times$ avec
\[\sigma_1(\epsilon)>1,0<\sigma_2(\epsilon)<1,\ldots,0<\sigma_n(\epsilon)<1, \] et l'hypothse de convergence. On en dŽduit
le principe de Koecher-Gštzky~:

\medskip

\begin{Lemme} Les sŽries (\ref{eq:fourier}) qui sont dans ${\cal A}_{\und{c}}(V,M)$ 
ont leur ensemble d'indi\-ces dans $M_{K,+}^*\cup\{0\}$, avec~:
\[M_{K,+}^*=\{\nu\in M_K^*\mbox{ avec }\sigma_i(\nu)\geq 0,i=1,\ldots,n\}.\]
De plus, si $d_0\not=0$ dans (\ref{eq:fourier}), alors $\und{c}=\und{0}$.\label{lemme:Koecher}\end{Lemme}

\medskip

\noindent {\bf Exemple.}
Pour $K=\QQ(\sqrt{5})$, $M={\cal O}_K,V={\cal O}_K^\times$, l'ensemble ${\cal A}_{\und{0}}(V,M)$ est l'anneau des fonctions $f$ analytiques
sur 
${\cal H}^n$, bornŽes au voisinage de $\infty$,
satisfaisant~:
\[f(S(\und{z}))=f(\und{z})\mbox{ et }f(U(\und{z}))=f(\und{z}).\]

\medskip

Pour $K\not=\QQ$,
toute forme modulaire de Hilbert $f$ de poids $(r_1,\ldots,r_n)$ pour $\Omega$ appartient ˆ ${\cal A}_{\und{c}}(V,M)$
pour certains $V,M$,  
avec $\und{c}=(r_2-r_1,\ldots,r_n-r_1)$.
On voit que
\[f(\und{z})=\sum_{\nu\in{\cal O}_K^*}d_\nu\exp\{2\pi i\tr(\nu \und{z})\},\] 
car $f$ est invariante par les translations de ${\cal O}_K$. Puis on applique le principe de Koecher
pour vŽrifier que l'ensemble des indices a son support dans ${\cal O}_{K,+}^*\cup\{0\}$.

Soit $a\in K$~; on considre une matrice $A\in\SL_2(K)$ telle que $A(a)=\infty$, on pose
$A^{-1}=\sqm{u}{v}{w}{x}$. Par exemple~;
\[A=\sqm{1}{\beta}{0}{1}\cdot\sqm{\alpha}{0}{0}{\alpha^{-1}}\cdot\sqm{0}{1}{-1}{a},\alpha\in K^*,\beta\in K\] convient. 
La fonction~:
\[f|_A(\und{z}):=\prod_{i=1}^n(\sigma_i(w)z_i+\sigma_i(x))^{-r_i}f(A^{-1}(\und{z}))\]
est une forme modulaire de poids $\und{r}$ pour le groupe $A\cdot\Omega\cdot A^{-1}$. Ce groupe a une pointe en $\infty$,
et $f|_A\in{\cal A}_{\und{c}}(\tilde{V},\tilde{M})$ pour $\tilde{V}\subset{\cal O}_K^\times$ et $\tilde{M}$ un $\ZZ$-module 
complet de $K$, que l'on peut expliciter en calculant le stabilisateur de $\infty$ dans $A\cdot\Omega\cdot A^{-1}$~:
donc $f$ est {\em rŽgulire en toute pointe}. En particulier~:
\[f|_A(\und{z})=\sum_{\nu\in \tilde{M}^*_+\cup\{0\}}d_\nu(A)\exp\{2\pi i\tr(\nu z)\}.\]

\medskip

\begin{small}
\noindent {\bf Exemples d'autres fonctions rŽgulires.}
Nous ajoutons un appendice ˆ ce paragraphe, bien qu'il n'y ait pas un lien direct avec les formes 
modulaires de Hilbert. On appelle {\em sŽrie gŽomŽtrique} de Hecke, la sŽrie~:
\[g_K(\und{z}):=\sum_{\nu\in {\cal O}^*_{K,+}\cup\{0\}}\exp\{\tr(\nu \und{z})\}\in{\cal A}_{\und{0}}(V,M),\]
avec $V={\cal O}_K^{\times}$ et $M={\cal O}_K$.

Les propriŽtŽs asymptotiques de cette sŽrie ont ŽtŽ ŽtudiŽes par Hecke dans \cite{Hecke:Analytische}.
Hecke a aussi ŽtudiŽ la sŽrie~:
\[b_K(\und{z}):=\sum_{\epsilon\in ({\cal O}_K^\times)^2}\exp\{\tr(\epsilon \und{z})\}\in{\cal A}_{\und{0}}(V,M).\]
Si $K\not=\QQ$, ces sŽries
sont transcendantes sur $\QQ(\exp_K(z_1,\ldots,z_n))$, pour toute fonction exponentielle $\exp_K:\CC^n\rightarrow\GG_m^n(\CC)$
ayant ses pŽriodes dans $\Sigma({\cal O}_K)\subset\CC^n$. Si $K=\QQ$ alors $g_K,b_K$ sont algŽbriques sur $\QQ(e^{2\pi i z})$.
On peut montrer que $g_K,b_K$ ne sont pas des formes modulaires de Hilbert. Ces sŽries ont des 
propriŽtŽs arithmŽtiques trs intŽressantes, liŽes ˆ la thŽorie de Mahler d'Žquations fonctionnelles 
associŽes ˆ certaines transformations monomiales.\end{small}

\subsection{PropriŽtŽs de base des formes modulaires de Hil\-bert.}

\noindent {\bf DŽfinition (poids parallle).} Soit $K$ un corps de nombres totalement rŽel de degrŽ $n$.
Nous dirons qu'une
forme modulaire de Hilbert associŽe ˆ $K$ est {\em de poids parallle $r\in\ZZ$,} si son poids est $(r,\ldots,r)\in\ZZ^n$.

\medskip

\noindent {\bf DŽfinition (formes paraboliques)}. Une {\em forme parabolique} pour $\Omega$ est une forme modulaire de Hilbert $f:{\cal H}^n\rightarrow\CC$ 
pour le groupe $\Omega$ telle que, pour tout $a\in\PP_1(K)$ et tout $A\in\SL_2(K)$ avec $A(a)=\infty$, on ait
\[f|_A(\und{z})=\sum_{\nu\in \tilde{M}^*_+\cup\{0\}}d_\nu(A)\exp\{2\pi i\tr(\nu z)\},\] avec $d_0(A)=0$.

\medskip

Pour $\und{r}\in\ZZ^n$, nous notons $M_{\und{r}}(\Omega)$ et $S_{\und{r}}(\Omega)$ respectivement les espa\-ces vectoriels de formes
modulaires de poids $\und{r}$ pour $\Omega$, et des formes para\-boliques pour $\Omega$ de poids $\und{r}$.

\begin{Lemme} Soit $\und{r}=(r_1,\ldots,r_n)\in\ZZ^n$. Si $r_1+\cdots+r_n<0$, alors $M_{\und{r}}(\Omega)=\{0\}$.
Si $r_1+\cdots+r_n=0$, alors $M_{\und{r}}(\Omega)\subset\CC$. En particulier, $M_{\und{0}}(\Omega)=\CC$. 
\label{lemme:M0=CC}\end{Lemme}

\noindent {\bf DŽmonstration.} Nous donnons ˆ la fin une deuxime preuve de l'ŽgalitŽ $M_{\und{0}}(\Omega)=\CC$, mais commenons 
avec une remarque importante.

Le groupe $\SL_2(\ZZ)$ se plonge dans tout groupe modulaire de Hilbert $\Gamma$~: par exemple,
son image dans $\Upsilon$ est engendrŽe par $S,T$.
Il existe donc un lien entre les formes modulaires de Hilbert et les formes modulaires elliptiques (formes
modulaires de Hilbert pour $K=\QQ$)~: soit
\[\Xi_n=\{(z,\ldots,z)\in{\cal H}^n\}\cong {\cal H},\]
soit $f$ une forme modulaire de poids $\und{r}=(r_1,\ldots,r_n)$ pour $\Gamma$.
Alors la restriction $g=f|_{\Xi_n}$ est une forme modulaire elliptique de poids $r_1+\cdots+r_n$.

Soient $f\in M_{\und{r}}(\Omega)$ et $B\in\SL_2(K)$. La restriction $f|_B(z,\ldots,z)$ est une forme modulaire
de poids $r_1+\cdots+r_n$, pour un sous-groupe d'indice fini de $\SL_2(\ZZ)$. Comme la rŽunion des images de $\Xi_n$
par les ŽlŽments $B\in\SL_2(K)$ est dense dans ${\cal H}^n$, si $f$ est non nulle, on peut choisir $B$
telle que $f|_B(z,\ldots,z)$ soit non nulle. Pour cela, il faut que $r_1+\cdots+r_n\geq 0$.

Supposons maintenant que $r_1+\cdots+r_n=0$.
Pour tout $B\in\SL_2(K)$ la forme modulaire $f|_B(z,\ldots,z)$ est constante, et ceci implique $f$ constante,
car dans ce cas, $f$ est une fonction analytique dans ${\cal H}^n$ telle que sur un sous-ensemble dense
de ${\cal H}^n$, il existe $n$ dŽrivŽes directionnelles linŽairement indŽpendantes annulant $f$.
En particulier, $M_{\und{0}}(\Omega)=\CC$.

Voici maintenant une deuxime dŽmonstration de ce dernier fait~: pour simplifier, nous supposons que
$\Omega$ a la seule classe d'Žquivalence de pointe $\infty$.
Soit $f\in M_{\und{0}}(\Omega)$~: 
Si $f$ n'est pas une forme parabolique, alors
\[f-f(\infty)\] est une forme parabolique.

D'aprs le lemme \ref{lemme:Koecher},
$f$ peut tre prolongŽe analytiquement en $\infty$, et en appliquant des
transformations de $\SL_2(K)$, par continuitŽ en tout ŽlŽment de $\PP_1(K)$~: $f(\PP_1(K))=0$. 

La compacitŽ de $X_\Omega$ implique que la fonction continue $|f|$ est
bornŽe sur $X_\Omega$, et le maximum de $|f|$ est atteint dans ${\cal H}^n$.
Le principe du maximum pour les fonctions analytiques de plusieurs variables
complexes implique que $f=0$. Dans tous les cas $f$ est constante.
\begin{Lemme} Supposons que $n\geq 2$.
Si $\und{r}=(r_1,\ldots,r_n)$ est tel que $r_1=0$ et il existe un indice $i$ avec $r_i\not=0$, alors
$M_{\und{r}}(\Omega)=\{0\}$. 
\label{lemme:M0r=CC}\end{Lemme}
\noindent {\bf DŽmonstration.}
Soit $f$ une forme modulaire de Hilbert de poids $\und{r}$, comme dans l'ŽnoncŽ. D'aprs le lemme \ref{lemme:Koecher}, $f$
est une forme parabolique. 

Comme dans la preuve du lemme \ref{lemme:M0=CC}, $f$ se prolonge par continuitŽ sur $\PP_1(K)$~: $f(\PP_1(K))=0$. 
La fonction \[g(\und{z})=|f(\und{z})|y_1^{r_1}y_2^{r_2}\cdots y_n^{r_n}=
|f(\und{z})|y_2^{r_2}\cdots y_n^{r_n}\] est continue, $\Omega$-invariante,
et se prolonge en fonction continue sur $X_\Omega$, nulle aux pointes.

Soit $\und{\zeta}\in{\cal H}^n$ un ŽlŽment tel que $\max g=g(\und{\zeta})$.
La fonction \[z_1\mapsto f(\und{z})y_2^{r_2}\cdots y_n^{r_n}\] est analytique sur ${\cal H}$,
et d'aprs le principe du maximum, la fonction~:
\[z_1\mapsto f(z_1,\zeta_2,\ldots,\zeta_n)\] est constante.

On vŽrifie que les seules fonctions $g\in{\cal A}_{\und{c}}(V,M)$ telles que $\partial g/\partial z_1=0$ sont les
fonctions constantes si $\und{c}=\und{0}$, et les fonctions nulles si $\und{c}\not=0$ (utiliser le fait
que la projection de $\Sigma({\cal O}_K)\subset\RR^n$ sur un facteur quelconque $\RR^{n-1}$ est dense). Ici, $\und{c}=
(r_2,\ldots,r_n)$ est non nul, donc $f$ est nulle.

\begin{Lemme} Si $\und{r}=(r_1,\ldots,r_n)$ est tel que $r_1<0$, alors
$M_{\und{r}}(\Omega)=\{0\}$.
\label{lemme:<0=0}\end{Lemme}
\noindent {\bf DŽmonstration.} Dans le paragraphe \ref{section:exemples}, nous construirons pour tout $r\in\ZZ_{\geq 1}$
une forme modulaire de Hilbert non nulle de poids $r\und{2}$ (une sŽrie d'Eisenstein). 
Soit $f\in M_{\und{r}}(\Omega)$ de poids non parallle $\und{r}=(r_1,\ldots,r_n)$, et supposons que $r_1<0$. Soit $E$ une 
forme modulaire de Hilbert de poids $\und{2}$ non nulle. Un certain
produit $f^aE^b$, pour $a,b$ entiers positifs non nuls, est une forme modulaire de poids $(0,*)$.
Le lemme \ref{lemme:M0r=CC} implique que $f^aE^b=0$, et $f=0$. Si $f$ est de poids parallle, ce raisonnement 
implique seulement que $f$ est une constante. Mais $\sum_ir_i=nr_1<0$, et le lemme \ref{lemme:M0=CC} implique que 
$f=0$.

\medskip

Nous notons $F_{\und{r}}(\Omega)$ l'espace vectoriel des fonctions modulaires $f:{\cal H}^n\rightarrow\CC$
de poids $\und{r}$, et $F(\Omega)=F_{\und{0}}(\Omega)$ le {\em corps} des fonctions modulaires de poids
$\und{r}=\und{0}$. Voici une autre consŽquence du principe de Koecher.

\medskip

\begin{Lemme} Si $K\not=\QQ$, il n'existe pas de fonction modulaire non constante dans $F(\Omega)$, sans p™les ou
sans zŽros dans ${\cal H}^n$.\label{lemme:poles}\end{Lemme}

\medskip

\noindent {\bf DŽmonstration.} Supposons que $f\in F_{\und{0}}(\Omega)$ soit sans p™les. Alors $f$ est holomorphe dans ${\cal H}^n$,
d'aprs le lemme \ref{lemme:Koecher}, donc
$f$ est une forme modulaire de poids $\und{0}$, qui est constante. Si $f$ est sans zŽros, alors $f^{-1}$ est 
sans p™les, donc holomorphe.

\medskip

Pour $K=\QQ$, cette propriŽtŽ est fausse car la fonction
modulaire $j$ est holomorphe dans ${\cal H}$.
Nous pouvons Žtendre ceci aux formes modulaires.

\begin{Lemme} Si $K\not=\QQ$, il n'existe pas de formes modulaires non constantes dans $M_{\und{r}}(\Omega)$, sans zŽros
dans ${\cal H}^n$.\label{lemme:nonexistence}\end{Lemme}

\noindent {\bf DŽmonstration.} Soit $f$ une forme modulaire de poids $\und{r}$, non constante et sans zŽros dans ${\cal H}^n$.
Quitte ˆ se ramener ˆ une puissance entire convenable de $f$, on peut supposer que $\und{r}\in(2\NN_{>0})^n$.

Nous verrons plus loin que~:
\[\lim_{m\rightarrow\infty}\dim_\CC M_{m\und{r}}(\Omega)=\infty.\]
Si $m$ est assez grand, il existe une forme modulaire de Hilbert $G$, de poids $m\und{r}$, telle que
$F^m,G$ soient $\CC$-linŽairement indŽpendantes.

Donc $G/F^m$ est une fonction modulaire de poids $\und{0}$, non constante, et sans p™les~: une contradiction.

\medskip

Une autre manire de procŽder est  d'appliquer les opŽrateurs \og slash\fg \\ (composition de formes modulaires
avec des transformations de ${\cal H}^2$ associŽes ˆ des matrices de $\GL_2(K)$~:
voir p. 251 de Zagier \cite{Zagier:Introduction}). La propriŽtŽ du lemme \ref{lemme:nonexistence} est fausse
pour $K=\QQ$, car la forme parabolique $\Delta$ de Jacobi ne s'annule pas dans ${\cal H}$.

\section{Exemples de formes modulaires de Hilbert.\label{section:exemples}}

Les techniques qui permettent de construire des formes modulaires de Hilbert, et les opŽrations qui permettent
de construire des nouvelles formes modulaires ˆ partir de 
certaines formes modulaires donnŽes, sont peu nombreuses mais puissantes.

\subsection{SŽries d'Eisenstein.\label{section:recette}}

Soit ${\mathfrak A}$ un idŽal non nul de ${\cal O}_K$, soit $r$ un entier rationnel
$\geq 1$, soit $s$ un nombre complexe de partie rŽelle $>2-r$. 
Notons ${\cal S}$ un sous-ensemble maximal de ${\mathfrak A}\times{\mathfrak A}$ ne contenant pas $(0,0)$, avec la propriŽtŽ que 
si $(m_1,l_1),(m_2,l_2)\in{\cal S}$ sont deux ŽlŽments distincts, alors pour toute unitŽ $\eta$ de $K$ totalement positive,
on a $\eta m_1\not=m_2$ et $\eta l_1\not=l_2$ .
La sŽrie~:
\[E_{r,{\mathfrak A}}(\und{z},s):=\sum_{(m,l)\in{\cal S}}\prod_{i=1}^n
(\sigma_i(m)z_i+\sigma_i(l))^{-r}|(\sigma_i(m)z_i+\sigma_i(l))|^{-s}\]
converge uniformŽ\-ment dans tout compact de ${\cal H}^n$ et dŽfinit une fonction holomorphe de la variable complexe $s$,
pour $\und{z}\in{\cal H}^n$ fixŽ.  En effet, pour $\und{z}=(z_1,\ldots,z_n)\in(\CC-\RR)^n$, 
il existe une constante $c>0$ telle que~:
\[|x z_i+y|^2\geq c(x^2+y^2),\] pour tout $(x,y)\in\RR^2$. Donc, pour $(m,l)\in K^2$~:
\[\left|\prod_{i=1}^n(\sigma_i(m)z_i+\sigma_i(l))^{-r-s}\right|\leq c^{-1}
\prod_{i=1}^n(\sigma_i(m)^2+\sigma_i(l)^2)^{-r/2-\Re(s)/2}.\]
Soit $K_1$ le corps de nombres $K({\rm i})$. Le nombre rationnel $\prod_i(\sigma_i(m)^2+\sigma_i(l)^2)$ est la 
valeur absolue de la norme
de $K_1$ sur $\QQ$ du nombre $m{\rm i}+l\in K_1$. On a donc majorŽ la valeur 
absolue du terme gŽnŽral de notre sŽrie
$E_{r,{\mathfrak A}}(\und{z},s)$ par le terme gŽnŽral d'une sŽrie extraite de la sŽrie 
dŽfinissant la fonction zta du corps $K_1$, multipliŽe
par une constante positive, d'o la convergence.

Soit ${\cal S}'$ un autre sous-ensemble maximal de ${\mathfrak A}\times{\mathfrak A}$ ayant la mme propriŽtŽ que ${\cal S}$.
Il est clair que~:
\begin{eqnarray*}\lefteqn{
\sum_{(m,l)\in{\cal S}}\prod_{i=1}^n(\sigma_i(m)z_i+\sigma_i(l))^{-r}|(\sigma_i(m)z_i+\sigma_i(l))|^{-s}=}\\
& = &\pm\sum_{(m',l')\in{\cal S}'}\prod_{i=1}^n(\sigma_i(m')z_i+\sigma_i(l'))^{-r}|(\sigma_i(m')z_i+\sigma_i(l'))|^{-s}.
\end{eqnarray*}
Le signe $\pm$ provient de la possible existence d'unitŽs de norme nŽgati\-ve.

Pour $s$ fixŽ dans le domaine de convergence, la fonction 
$E_{r,{\mathfrak A}}(\cdot,s):{\cal H}^n\rightarrow\CC$ ainsi dŽfinie
est \og automorphe de poids $(r+s,\ldots,r+s)$\fg pour $\Gamma$. Par exemple, si $r\geq 3$ et $s=0$, 
on construit de cette faon des
formes modulaires de Hilbert de poids $(r,\ldots,r)$. Pour $\gamma=\sqm{a}{b}{c}{d}\in\Gamma$ on a de faon plus dŽtaillŽe~:
\begin{eqnarray*}
E_{r,{\mathfrak A}}(\gamma(\und{z}),s) & = & \prod_{i=1}^n(\sigma_i(c)z_i+\sigma_i(d))^{r}|\sigma_i(c)z_i+\sigma_i(d)|^s\\
& &\times\sum_{(m,n)\in{\cal S}}
\prod_{i=1}^n(\sigma_i(m')z_i+\sigma_i(l'))^{-r}|(\sigma_i(m')z_i+\sigma_i(l'))|^{-s},
\end{eqnarray*}
o $m'=am+cl,l'=bm+dl$. Mais la transformation dŽfinie par ces relations envoie ${\cal S}$ en un autre 
sous-ensemble ${\cal S}'$ de ${\mathfrak A}^2$ ayant les mmes propriŽtŽs que ${\cal S}$.

Si $r$ est impair, $E_{r,{\mathfrak A}}$ est identiquement nulle quand $K$ a une unitŽ
de norme $-1$. Toutefois, $E_{r,{\mathfrak A}}$ n'est pas en gŽnŽral identiquement nulle.

Pour $r\geq 4$ pair, $E_{r,{\mathfrak A}}(\und{z},0)$ n'est jamais nulle. Pour le voir, 
il suffit de choisir un ordre particulier de sommation dans la sŽrie $\sum_{(m,n)\in{\cal
S}}$. Nous pouvons en effet Žcrire, gr‰ce ˆ la convergence absolue~:
\[\sum_{(m,l)\in{\cal S}}=\sum_{m=0,l\in{\cal T}}+\sum_{m\in{\cal T},l\in{\mathfrak A}},\]
o ${\cal T}$ est un sous-ensemble maximal de ${\mathfrak A}$ ayant la propriŽtŽ que $0\not\in{\cal T}$ et si $n_1,n_2\in{\cal T}$
avec $n_1\not=n_2$, alors $n_1/n_2$ n'est pas une unitŽ totalement positive de $K$. La sŽrie~:
\[\sum_{m\in{\cal T},l\in{\mathfrak A}}\prod_{i=1}^n(\sigma_i(m)z_i+\sigma_i(l))^{-r}\]
est nulle en $\und{z}=\infty$. D'autre part,
\[\sum_{m=0,l\in{\cal T}}\prod_{i=1}^n(\sigma_i(m)z_i+\sigma_i(l))^{-r}=\sum_{l\in{\cal T}}\no(l)^{-r},\]
est non nul. Donc la sŽrie d'Eisenstein $E_{r,{\mathfrak A}}(\und{z})=E_{r,{\mathfrak A}}(\und{z},0)$
est une forme modulaire de Hilbert de poids $(r,\ldots,r)$, 
qui n'est pas une forme parabolique (et n'est pas non plus constante).

\medskip

Plus $r\geq 3$ est petit, plus la sŽrie d'Eisenstein $E_{r,{\mathfrak A}}(\und{z})=E_{r,{\mathfrak A}}(\und{z},0)$ 
converge lentement, plus elle est \og intŽressante\fg, comme
nous le verrons plus loin. Hecke (\cite{Hecke:Analytische} pp.
391-395) a introduit la technique du prolongement analytique des fonctions analytiques d'une variable 
complexe, pour construire des fonctions satisfaisant des relations d'automorphie de poids $\und{2}$ et $\und{1}$. 
Il construit ces fonctions
en explicitant d'abord le dŽveloppement en sŽrie de Fourier de $E_{r,{\mathfrak A}}(\und{z},s)$, pour $s$ fixŽ.
En appliquant la formule de Poisson, il remarque que les coefficients de Fourier sont des fonctions analytiques 
de la variable complexe $s$, et qu'on
peut faire un prolongement analytique. Puis il calcule une limite pour $s\rightarrow 0$,
et remarque que trs souvent, les zŽros des facteurs gamma Žliminent les termes
non holomorphes en $\und{z}$ 
Ce procŽdŽ est
devenu classique, c'est pourquoi nous n'en avons donnŽ qu'une esquisse.

Pour $r=2$,
on trouve, lorsque $\Re(s)>0$, que la limite~:
\[E_{2,{\mathfrak A}}(\und{z}):=
\lim_{s\rightarrow 0}\sum_{(m,l)\in{\cal
S}}\prod_{i=1}^n(\sigma_i(m)z_i+\sigma_i(l))^{-2}|(\sigma_i(m)z_i+\sigma_i(l))|^{-s}\] existe. Si $K\not=\QQ$, alors
$E_{2,{\mathfrak A}}(\und{z})$ est holomorphe, et c'est donc  une forme modulaire pour $\Gamma$, de poids $\und{2}$, qui est
non nulle.

Si $K=\QQ$, $E_{2,\ZZ}(z)$ n'est pas holomorphe (mais possde des propriŽtŽs d'automorphie), et
\[E_2(z)=\zeta_\QQ(2)^{-1}E_{2,\ZZ}(z)+\displaystyle{\frac{3}{\pi\Im(z)}}\] 
est holomorphe, mais non modulaire (cf.
l'article de Bertrand dans \cite{Nesterenko:Introduction}, chapitre 1). Le terme $\displaystyle{\frac{3}{\pi\Im(z)}}$
appara"t aprs avoir fait un prolongement analytique.

Les sŽries $E_{r,{\mathfrak A}}(\und{z})$ poss\-dent des dŽveloppements en sŽrie de Fourier qui peuvent
tre calculŽs explicitement~: ce calcul est trs utile si ${\mathfrak A}={\cal O}_K$. 
Soit $\zeta_K(s)$ la fonction zta associŽe ˆ $K$, Žcrivons~: 
\[\zeta_K(r)^{-1}E_{r,{\cal O}_K}(\und{z})=\sum_{\nu\in {\mathfrak A}^*_{+}\cup\{0\}}d_\nu\exp\{2\pi
{\rm i}\tr(z\nu)\}.\] Pour $r$ pair, les coefficients $d_\nu$ sont des nombres rationnels 
que l'on peut calculer explicitement. En particulier, $d_0=1$.  Les dŽtails 
de ces calculs sont dŽcrits dans
\cite{Geer:Hilbert} pp. 19-21~: ils sont classiques, c'est pourquoi 
nous ne les reportons pas ici en toute gŽnŽralitŽ.

\medskip

\noindent {\bf Une recette.}
Voici une faon de calculer les coefficients de Fourier des sŽries $E_{r,{\mathfrak A}}(\und{z})$ 
dans le cas o $[K:\QQ]=2$, $r$ pair, nombre de classes d'idŽaux $1$, 
et ${\mathfrak A}={\cal O}_K$ (nous rendons plus explicites les  formules de la proposition 6.4 pp. 19-20 de
\cite{Geer:Hilbert}, pour les applications que nous avons en vue). On trouve~:
\[E_{r,{\cal O}_K}(\und{z})=\zeta_K(r)(1+\sum_{\nu\in{\cal O}^*_{K,+}}b_r(\nu)\exp\{2\pi{\rm i}\tr(\nu\und{z})\}),\]
o \[b_r(\nu)=\kappa_r\sum_{(\mu)|\nu\sqrt{d}}|\no(\mu)|^{r-1},\] la somme Žtant Žtendue aux 
idŽaux entiers $(\mu)$ de $K$ (donc principaux) qui divisent $\nu\sqrt{d}$. 
Nous avons posŽ~:
\[\kappa_r=\frac{(2\pi)^{2r}\sqrt{d}}{((r-1)!)^2d^r\zeta_K(r)},\] et $d$ dŽsigne le discriminant de $K$.

En particulier, les sŽries d'Eisenstein $E=E_{r,{\cal O}_K}$ pour $K$ quadratique rŽel et $r$ pair satisfont~:
\[E(z,z')=E(z',z).\] On dit que ce sont des formes modulaires {\em symŽtriques}.

\medskip

\noindent {\bf Un calcul explicite.} Soit $K=\QQ(\sqrt{5})$. En utilisant les formules dŽcrites ci-dessus on voit que~:
\[\zeta_K(2)^{-1}E_{2,{\cal O}_K}(\und{z})=(1+120(\exp\{2\pi{\rm i}\tr(\mu \und{z})\}+\exp\{2\pi {\rm i}\tr(\mu'
\und{z})\})+\cdots)\] o $\displaystyle{\mu=\frac{1+\sqrt{5}}{2\sqrt{5}}}$.
Comme~:
\[E_4(t)=1+240\exp\{2\pi i t\}+\cdots\] on voit que
\[\zeta_K(2)^{-1}E_{2,{\cal O}_K}|_{\Xi_2}=E_4.\]

\subsection{Fonctions thtas.\label{section:theta}} La nature des sŽries d'Eisenstein ne permet pas toujours d'obtenir 
la non nullitŽ des formes modulaires ainsi construites (exemples~: sŽries de poids impairs et $K$ a une unitŽ de norme nŽgative, ou 
sŽries d'Eisenstein de poids $\und{1}$ tordues par des caractres). 
Les fonctions thtas en revanche, sont toujours non nulles. 

Il convient de faire un dŽtour passant par les groupes symplectiques et les fonctions thta classiques.

Nous notons ${\cal H}_n\subset\CC^n$ le {\em demi-espace supŽrieur de Siegel} constituŽ des matrices complexes symŽtriques 
$\und{Z}$ ˆ $n$ lignes et $n$ colonnes,
ayant une partie imaginaire dŽfinie positive. Pour 
$(\und{u},\und{Z},\und{r},\und{s})\in\CC^n\times{\cal H}_n\times\RR^n\times\RR^n$, posons~:
\[\vartheta(\und{u},\und{Z};\und{r},\und{s})=\sum_{x\in\ZZ^n+\und{r}}\exp\left\{2\pi{\rm i}
\left(\frac{1}{2}{}^t x\cdot\und{Z}\cdot x+
{}^t x\cdot(\und{u}+\und{s})\right)\right\}.\] 
Cette sŽrie converge uniformŽment dans les compacts de $\CC^n\times{\cal H}_n$ et dŽfinit
la fonction {\em thta de caractŽristique $(\und{r},\und{s})$}.
Posons~:
\[E_n=\sqm{0}{1_n}{-1_n}{0}\in\GL_{2n}(\RR),\] o $1_n$ est la matrice identitŽ d'ordre $n$.
Le groupe {\em symplectique} $\SP_n(\RR)\subset\GL_{2n}(\RR)$ est dŽfini par~:
\[\SP_n(\RR)=\left\{T\mbox{ tels que }T\cdot E_n\cdot{}^tT=E_n\right\}.\] 
Il agit sur ${\cal H}_n$ de la manire suivante. Si $\und{Z}\in{\cal H}_n$ et 
$U=\sqm{A}{B}{C}{D}\in\SP_n(\RR)$ avec $A,B,C,D$ matrices carrŽes d'ordre $n$, alors on pose~:
\[U(\und{Z})=(A\cdot\und{Z}+B)\cdot(C\cdot\und{Z}+D)^{-1}.\]
Noter que $\SP_1(\RR)=\SL_2(\RR)$, et $\SP_n(\RR)\subset\SL_{2n}(\RR)$.

Soit  $S$ une matrice ˆ $n$ lignes et $n$ colonnes, notons $\{S\}\in\CC^n$ le vecteur colonne dont 
les coefficients sont les ŽlŽments de la diagonale de $S$. 

Nous recopions ici un cas particulier d'une proposition trs classique que l'on peut trouver, 
par exemple, dans \cite{Shimura:Theta}, proposition 1.3, pp. 676-677, dont nous omettons la dŽmonstration.
\begin{Proposition}
Soient $\und{r},\und{s}$ des ŽlŽments de $(\ZZ/2)^n$. Pour tout $U=\sqm{A}{B}{C}{D}\in\SP_n(\ZZ)$, on a~:
\[\vartheta(\und{0},U(\und{Z});\und{r},\und{s})=\zeta_U\det(C\cdot\und{Z}+D)^{1/2}\vartheta(\und{0},\und{Z};\und{r}',\und{s}'),\]
o $\zeta_U$ est une racine huitime de l'unitŽ qui dŽpend de $\und{r},\und{s}$ et $U$,
\[\binomial{\und{r}'}{\und{s}'}={}^tU\cdot\binomial{\und{r}}{\und{s}}+\frac{1}{2}\binomial{\{{}^tA\cdot C\}}{\{{}^tB\cdot D\}}.\]
En particulier pour tout $U=\sqm{A}{B}{C}{D}\in\SP_n(\ZZ)(8)$ {\em (\footnote{Sous-groupe de congruence principal de niveau
$8$, engendrŽ par toutes les matrices $U$ de $\SP_n(\ZZ)$ congrues ˆ la matrice unitŽ modulo $8$M${}_{2n,2n}(\ZZ)$.})}~:
\begin{equation}
\vartheta(\und{0},U(\und{Z});\und{r},\und{s})=\zeta_U\det(C\cdot\und{Z}+D)^{1/2}\vartheta(\und{0},\und{Z};\und{r},\und{s}).\label{eq:modularite_vartheta}
\end{equation}
\label{lemme:shimura1}
\end{Proposition}
On dit que $\vartheta(0,U(\und{Z});\und{r},\und{s})$ est une {\em forme modulaire de Siegel} de poids $1/2$, niveau $8$, et 
systme multiplicateur cyclotomique d'ordre divisant $8$.

\medskip 

Nous revenons ˆ notre corps totalement rŽel $K$, et au groupe modulaire de Hilbert $\Gamma$.
Soit ${\cal I}$ un ${\cal O}_K$-module inversible de rang $1$ de $K$, posons~:
\[\SL_2({\cal O}_K\oplus{\cal I}):=\SL_2(K)\cap\sqm{{\cal O}_K}{{\cal I}^{-1}}{{\cal I}}{{\cal O}_K}.\]
Soit $(\beta_1,\ldots,\beta_n)$ une 
base de ${\cal I}$ sur $\ZZ$. ConsidŽrons la matrice~:
\[B=\left(\begin{array}{ccc}
\sigma_1(\beta_1) & \cdots & \sigma_1(\beta_n)\\ \vdots &  & \vdots\\ \sigma_n(\beta_1) & \cdots & \sigma_n(\beta_n)\end{array}\right).\]
L'application~:
\[W_B(\und{z})={}^tB\cdot\mbox{Diag}(z_1,\ldots,z_n)\cdot B,\] dŽfinit une application $W_B:{\cal H}^n\rightarrow {\cal H}_n$.
Soit $a$ un ŽlŽment de $K$ et posons $\Phi(a)=\mbox{Diag}(\sigma_1(a),\ldots,\sigma_n(a))$. Nous avons un homomorphisme
de groupes $I:\SL_2(K)\rightarrow
\SP_n(\QQ)$ dŽfini par~:
\begin{eqnarray*}
I\left(\sqm{a}{b}{c}{d}\right)
& = &\sqm{{}^tB}{0}{0}{B^{-1}}\cdot\sqm{\Phi(a)}{\Phi(b)}{\Phi(c)}{\Phi(d)}\cdot\sqm{{}^tB^{-1}}{0}{0}{B}\\
& = &
\sqm{{}^tB\cdot\Phi(a)\cdot{}^tB^{-1}}{{}^tB\cdot\Phi(b)\cdot B}{B^{-1}\cdot\Phi(c)\cdot{}^tB^{-1}}{B^{-1}\cdot\Phi(d)\cdot B}.
\end{eqnarray*}
On voit facilement que $I(\SL_2({\cal O}_K\oplus{\cal I}))\subset\SP_n(\ZZ)$.

Soit $\gamma\in\SL_2(K)$, Žcrivons $I(\gamma)=\sqm{P}{Q}{R}{S}$. Nous avons les propriŽtŽs de compatibilitŽ
suivantes, qui peuvent tre vŽrifiŽes sans difficultŽ (voir \cite{Shimura:Theta} p. 683, voir aussi \cite{Hammond:Modular}
pp. 499-506)~:
\begin{eqnarray}
I(\gamma)(W_B(\und{z})) & = & W_B(\gamma(\und{z}))\label{eq:comp1}\\
R\cdot W_B(\und{z})+S & = & B^{-1}\cdot\label{eq:comp2}\\ 
& &\cdot\mbox{Diag}(\sigma_1(c)z_1+\sigma_1(d),\ldots,\sigma_n(c)z_n+\sigma_n(d))\cdot B.\nonumber
\end{eqnarray} 

Nous dŽfinissons maintenant les fonctions thta sur ${\cal H}^n$ associŽes ˆ $K$.
Soient $\rho,\delta$ des ŽlŽments de $K$, considŽrons des variables complexes $(\und{z},\und{t})\in{\cal H}^n\times\CC^n$, posons~:
\begin{equation}
\theta(\und{t},\und{z};\rho,\delta)=
\sum_{\nu\in {\cal I}+\rho}\exp\left\{2\pi {\rm i}\;\tr\left(\frac{1}{2}\und{z}\nu^2+\nu(\und{t}+\delta)\right)\right\}.
\label{eq:def_theta}
\end{equation}
Cette sŽrie converge uniformŽment dans les compacts de $\CC^n\times{\cal H}^n$.
Une vŽrification directe nous donne
\[\theta(\und{t},\und{z};\rho,\delta)=\vartheta({}^tB\cdot\und{t},W_B(\und{z});B^{-1}\cdot\rho,{}^tB\cdot\delta).\]
Toutes les propriŽtŽs d'automorphie dŽcrites ci-dessus pour les fonctions $\vartheta$ se transmettent sur les fonctions
$\theta$ via les conditions de compatibilitŽ (\ref{eq:comp1}) et (\ref{eq:comp2}).
Si par exemple ${\cal I}={\cal O}_K^*$, $\rho\in {\cal O}_K^*/2$ et $\delta\in {\cal O}_K/2$, la fonction~:
\[\und{z}\mapsto\theta(\und{0},\und{z};\rho,\delta)\] est une forme modulaire de Hilbert de poids $(1/2,\ldots,1/2)$ pour $\Gamma^\sharp(8)$,
o $\Gamma^\sharp=\SL_2({\cal O}_K\oplus{\cal I})$ et
\[\Gamma^\sharp(8)=\left\{\gamma\in\Gamma^\sharp\mbox{ tel que }\gamma\equiv\sqm{1}{0}{0}{1}\;\mbox{mod} \;(8{\cal
O}_K)\right\},\] avec caractre cyclotomique d'ordre divisant $8$, d'aprs la proposition \ref{lemme:shimura1}, et les
conditions de  compatibilitŽ (\ref{eq:comp1}) et (\ref{eq:comp2}).

\subsection{Le cas $n=2$ dŽtaillŽ.\label{section:cass=2}}
Ici on pose $n=2$, et on s'occupe des caractŽristiques $(\und{r},\und{s})\in(\ZZ/2)^4$.
L'ensemble $\{\vartheta(\und{t},\und{Z};\und{r},\und{s})\}_{(\und{r},\und{s})\in(\ZZ/2)^4}$ est de cardinal $16$.
En fixant $\und{Z}\in{\cal H}_n$, la fonction $\und{t}\mapsto\vartheta(\und{t},\und{Z};\und{r},\und{s})$ peut tre une
fonction paire ou impaire, comme le suggre une vŽrification directe. Plus prŽcisŽment~:
\[\vartheta(-\und{t},\und{Z};\und{r},\und{s})=(-1)^{|(\und{r},\und{s})|}\vartheta(\und{t},\und{Z};\und{r},\und{s}),\]
o $|(\und{r},\und{s})|=4(r_1s_1+r_2s_2)\in\ZZ$. On en dŽduit qu'il existe exactement $10$ fonctions
thta paires, ou de manire Žquivalente, $10$ {\em caractŽristiques paires}. Notons ${\mathfrak K}\subset(\ZZ/2)^4$ un ensemble 
complet de classes de congruence distinctes de $(\ZZ/2)^4$ modulo $\ZZ^4$, et ${\mathfrak K}^*\subset{\mathfrak K}$ un sous-ensemble
complet de reprŽsentants de classes de caractŽristiques paires.

Nous considŽrons une action $\alpha$ de $\SP_2(\ZZ)$ sur $(\ZZ/2)^4$ dŽfinie de la manire suivante 
(suggŽrŽe par la proposition \ref{lemme:shimura1}).
Pour $U=\sqm{A}{B}{C}{D}\in\SP_2(\ZZ)$ et $(\und{r},\und{s})\in(\ZZ/2)^4$,
\[\alpha_U\binomial{\und{r}}{\und{s}}={}^tU\cdot\binomial{\und{r}}{\und{s}}+
\frac{1}{2}\binomial{\{{}^tA\cdot C\}}{\{{}^tB\cdot D\}}.\] On vŽrifie que c'est une action de groupe.

\begin{Lemme} Pour tout $U\in\SP_2(\ZZ)$, l'application $\alpha_U:{\mathfrak K}\rightarrow{\mathfrak K}$ est une bijection.
L'action $\alpha$ agit aussi sur ${\mathfrak K}^*$.
\label{lemme:action2}\end{Lemme}

\noindent {\bf DŽmonstration.} On peut utiliser le fait que $\SP_2(\ZZ)$ est engendrŽ par les matrices~:
\[\sqm{0}{-1_2}{1_2}{0}\mbox{ et }\sqm{1_2}{T}{0}{1_2},\] avec $T={}^tT$. De plus, l'action se factorise
par une action de $\SP_2(\FF_2)\cong{\mathfrak S}_6$.

\medskip

En particulier, la fonction~:
\[\vartheta(\und{Z})=\prod_{(\und{r},\und{s})\in{\mathfrak K}^*}\vartheta(\und{0},\und{Z};\und{r},\und{s})\]
est une forme modulaire de Siegel de poids $5$, avec caractre cyclotomique d'ordre divisant $8$ (lire \cite{Igusa: Siegel}).
On en dŽduit que~:
\begin{equation}
\Theta^\sharp(\und{z})=\prod_{(\rho,\delta)\in{\mathfrak K}^*}\vartheta(\und{0},W_B(\und{z});\rho,\delta)=\vartheta(W_B(\und{z}))
\label{eq:theta_big}\end{equation}
est une forme modulaire de Hilbert de poids $(5,5)$, avec caractre 
pour le groupe $\SL_2({\cal O}_K\oplus{\cal I})$. Cette forme modulaire n'est pas nulle en gŽnŽral, comme expliquŽ 
dans \cite{Hammond:Modular}, p. 507.

\section{La forme parabolique $\Theta$.}

La forme para\-bolique de Jacobi $\Delta:{\cal H}\rightarrow\CC$, de poids $12$ pour $\SL_2(\ZZ)$, peut tre dŽfinie par~:
\[\Delta(z^2)=2^{-8}\left(\theta_2(z)\theta_3(z)\theta_4(z)\right)^{8},\] o les $\theta_i$ sont les fonctions
\begin{eqnarray}
\theta\left(0,z;0,0\right) , & \displaystyle{\theta\left(0,z;\frac{1}{2},0\right)}. & \theta\left(0,z;0,\frac{1}{2}\right).
\label{eq:les_f_theta}
\end{eqnarray}
Nous cherchons des analogues de $\Delta$ en dimension supŽrieure.

Dans ce paragraphe, nous construisons explicitement une forme modulaire $\Theta$
de poids $(5,5)$ pour $\Upsilon=\SL_2({\cal O}_K)$ (avec $K=\QQ(\sqrt{5})$), dont le lieu des zŽros
est l'image de $\Xi_2$ modulo l'action de $\Upsilon$, avec multiplicitŽ $1$. 

\medskip

Soit $\mu$ un ŽlŽment de $K^\times$, 
soit ${\cal I}$ un ${\cal O}_K$-module inversible de rang $1$. On a un
isomorphisme
\[s_\mu:\SL_2({\cal O}_K\oplus{\cal I})\rightarrow\SL_2({\cal O}_K\oplus(\mu){\cal I}),\] dŽfini par~:
\[\gamma=\sqm{a}{b}{c}{d}\mapsto\sqm{a}{\mu^{-1}b}{\mu c}{d}.\]
On a que $\gamma(\sigma(\mu)\und{z})=\sigma(\mu)s_\mu(\gamma)(\und{z})$.

Si $F(z_1,z_2)$ est une forme modulaire de Hilbert de poids $\und{r}$ pour $\Gamma({\cal O}_K\oplus{\cal I})$, et si
$\mu\in K$ est tel que $\mu>0$ et $\mu'>0$, alors $F(\mu^{-1}z_1,\mu'{}^{-1}z_2)$ est une forme modulaire
de Hilbert de poids $\und{r}$ pour $\Gamma({\cal O}_K\oplus(\mu){\cal I})$.

\medskip

Posons 
$\epsilon=\displaystyle{\frac{1+\sqrt{5}}{2}}$.
On prend ${\cal I}={\cal O}_K^*$ dans (\ref{eq:theta_big}), puis on compose $\Theta^\sharp$. Ensuite, 
on utilise les propriŽtŽs de l'application $s_\mu$ ci-dessus, et on vŽrifie que~:
\[\Theta(z_1,z_2)=\Theta^\sharp\left(\epsilon \frac{z_1}{\sqrt{5}},-\epsilon'\frac{ z_2}{\sqrt{5}}\right)\] est une forme modulaire
de Hilbert de poids $(5,5)$ avec caractre, pour le groupe modulaire $\Upsilon$.

Mais $\Upsilon$ n'a pas de caractre
non trivial (cf. \cite{Maass:Modulformen} p. 72). Donc la fonction $\Theta$ est une {\em forme modulaire de Hilbert}
de poids $(5,5)$ pour $\Upsilon$~: c'est de plus une forme parabolique, car certains facteurs du produit 
(\ref{eq:theta_big}) le sont (\footnote{On peut trouver une dŽmonstration plus directe de ces propriŽtŽs 
dans \cite{Gundlach:Bestimmung1} pp. 231-237. Notre construction est plus susceptible de gŽnŽralisation.}).

On considre~:
\begin{eqnarray*}\lefteqn{{\cal E}:=\{(\alpha_i,\beta_i),i=1,\ldots,10\}=}\\
& = & \{(0,0),(1,0),(0,1),(1,1),(0,\epsilon'),(0,\epsilon),
(\epsilon',0),(\epsilon,0),(\epsilon',\epsilon),(\epsilon,\epsilon')\}\\ &\subset &
{\cal O}_K\times{\cal O}_K.\end{eqnarray*}
On a~:
\begin{eqnarray}
\Theta(\und{z}) &= &\prod_{i=1}^{10}\theta\left(0,\left(\frac{z_1\epsilon}{\sqrt{5}},-\frac{z_2\epsilon'}{\sqrt{5}}\right);
\alpha_i,\frac{\beta_i}{\sqrt{5}}\right)\nonumber\\
& = & \prod_{i=1}^{10}\sum_{\nu\in{\cal O}_K}(-1)^{\tr(\nu\beta_i/\sqrt{5})}\exp\left\{\pi{\rm
i}\tr\left(\left(\nu+\frac{\alpha_i}{2}\right)^2\frac{\epsilon\und{z}}{\sqrt{5}}\right)\right\}.\label{eq:produit_10}
\end{eqnarray}

Une vŽrification directe que nous laissons au lecteur entra"ne (\footnote{Utiliser le fait que ${\cal O}_K'={\cal O}_K$.})~:
\[\Theta(z,z')=-\Theta(z',z).\] On dit que $\Theta$ est une forme
modulaire de Hilbert {\em antisymŽtrique}. Elle s'annule dans $\Xi_2=\{(z,z')\mbox{ avec }z=z'\}$. Nous montrons
maintenant que ceci dŽtermine bien le lieu d'annulation de $\Theta$, et que la multiplicitŽ est $1$.

\subsection{Un domaine fondamental ŽloignŽ du bord de ${\cal H}^2$.\label{section:goetzky}}

Nous commenons par montrer qu'il existe un domaine fondamental pour l'action de $\Upsilon$ sur ${\cal H}^2$
plus convenable que celui que nous avons introduit dans le paragraphe \ref{section:description}. Nous voulons
travailler avec un domaine fondamental qui soit le plus ŽloignŽ que possible du bord de ${\cal H}^2$,
et celui que nous avons construit au paragraphe \ref{section:sans_demonstration}, ne l'est pas assez. 
La construction que
nous  donnons est celle de \cite{Goetzky:Anwendung}, pp. 416-422.
On considre les parties ${\mathfrak A},{\mathfrak B},{\mathfrak C}$ de ${\cal H}^2$ suivantes.

\medskip

\noindent On pose ${\mathfrak A}=\{(z,z')\in{\cal H}^2\mbox{ tel que }|zz'|\geq 1\}$. C'est un domaine fondamental
pour l'action du sous-groupe $\{1,T\}\subset\Upsilon$.

\medskip

\noindent On pose ${\mathfrak B}$ le sous-ensemble de tous les couples $(z,z')\in{\cal H}^2$ ayant la propriŽtŽ
que $\epsilon'{}^2\leq \Im(z)\Im(z')\leq\epsilon^2$. C'est un domaine fondamental pour l'action du sous-groupe
$U^\ZZ=\{U^{n}\mbox{ tel que }n\in\ZZ\}\subset\Upsilon$, o $U$ est dŽfini dans (\ref{eq:relations_goetzky}).

\medskip

\noindent Pour dŽfinir ${\mathfrak C}$, nous allons d'abord dŽcrire l'ensemble ${\mathfrak C}_{s,s'}$ dont les ŽlŽments sont les
$(z,z')\in{\mathfrak C}$ tels que $\Im(z)=s$ et $\Im(z')=s'$, de telle sorte que (union disjointe)~:
\[{\mathfrak C}=\bigcup_{s,s'\in\RR_{>0}}{\mathfrak C}_{s,s'}.\]
Soit ${\cal Q}_{s,s'}$ un domaine fondamental 
de $\RR^2$ pour l'action des translations de $\Sigma({\cal O}_K)$ ayant la propriŽtŽ que pour tout $(r,r')\in{\cal Q}_{s,s'}$
la quantitŽ~:
\[(r^2+s^2)(r'{}^2+s'{}^2)=|(r+{\rm i}s)(r'+{\rm i}s')|^2\] est la plus petite possible.
Alors ${\mathfrak C}_{s,s'}$ est le translatŽ de ${\cal Q}_{s,s'}$ par $({\rm i}s,{\rm i}s')$ dans $\CC^2$. 
Avec cette construction,
${\mathfrak C}$ est un domaine fondamental pour l'action des translations de ${\cal O}_K$.

\medskip

\begin{Lemme} L'ensemble ${\cal G}={\mathfrak A}\cap{\mathfrak B}\cap {\mathfrak C}\subset{\cal H}^2$ 
contient un domaine fondamental pour l'action
de $\Upsilon$ sur ${\cal H}^2$. De plus, si $(z,z')\in{\cal G}$, alors $\Im(z)\Im(z')>0.54146$.
\label{lemme:domaine_de_goetzky}
\end{Lemme}

\noindent {\bf DŽmonstration.} 
Ici on Žcrira $\und{t}_i=(t_i,t_i')$ et $\und{z}_i=(z_i,z_i')$.
Il est clair
que ${\mathfrak B}\cap {\mathfrak C}$ est un domaine fondamental pour l'action de $\Upsilon_\infty$ sur ${\cal H}^2$.
Montrons maintenant que pour tout $\und{z}\in{\cal H}^2$, il existe $\gamma\in\Upsilon$ tel que $\gamma(\und{z})\in{\cal G}$.

Soit donc $\und{z}_1\in{\cal H}^2$. Il existe $\gamma_1\in\Upsilon_\infty$ tel que $\und{t}_1=
\gamma_1(\und{z}_1)\in{\mathfrak B}\cap {\mathfrak C}$. Si $\Im(t_{1})\Im(t_{1}')\geq 1$ alors $|t_1t_1'|\geq 1$ et 
$\und{t}_1\in{\cal
G}$. Sinon $|t_1t_1'|<1$ et $\und{z}_2=(z_{2},z_{2}')=T(\und{t}_1)\in{\mathfrak A}$. Observons que~:
\begin{eqnarray*}
\Im(z_{2})\Im(z_{2}') & = & \Im(z_{1})\Im(z_{1}')|t_{1}t_{1}'|^{-2}\\
& > &\Im(z_{1})\Im(z_{1}').
\end{eqnarray*}
Il existe un ŽlŽment $\gamma_2\in\Upsilon_\infty$ tel que $\und{t}_2=\gamma_2(\und{z}_2)\in{\mathfrak B}\cap {\mathfrak C}$.
Si $\Im(t_2)\Im(t_2')\geq 1$ alors $\und{t}_2\in{\cal G}$, et
nous avons dŽmontrŽ un cas particulier du lemme pour $\und{z}_1$. Sinon,
$|t_2t_2'|< 1$, et nous pouvons continuer en construisant successivement $\und{z}_3,\und{t}_3,\und{z}_4,\ldots$.

Nous avons une suite $(\und{z}_i)_{i=1,2,\ldots}$ d'ŽlŽments de ${\cal H}^2$ Žquivalents modulo l'action de $\Upsilon$, telle que~:
\[\Im(z_{1})\Im(z_{1}')< \Im(z_{2})\Im(z_{2}')<\Im(z_{3})\Im(z_{3}')<\cdots.\]
Si aucun de ces points n'est dans ${\cal G}$, alors ils sont tous dans le domaine~:
\[\{(z,z')\in{\cal H}\mbox{ tel que }|zz'|\leq
1,\epsilon'{}^2\leq\frac{\Im(z)}{\Im(z')}\leq\epsilon^2,\Im(z)\Im(z')\geq\Im(z_{1})\Im(z_{1}')\}.\] Ce domaine est compact, et
l'action de $\Upsilon$ est propre~: nous avons une contradiction, car un sous-ensemble compact de ${\cal H}^2$ ne contient
qu'une finitude d'ŽlŽments Žquivalents modulo $\Upsilon$. Donc, pour tout $\und{z}\in{\cal H}^2$, il existe
$\gamma\in\Upsilon$ tel que $\gamma(\und{z})\in{\cal G}$, et ${\cal G}$ est un ensemble fondamental.

Ensuite, nous allons prouver que si $\und{z}\in{\mathfrak A}\cap{\mathfrak C}$, alors $\Im(z)\Im(z')>0.54$. 
Soit $\aa$ un nombre rŽel tel que $|\aa|\leq 1/\sqrt{5}$, 
soit $P(\aa)$ le parallŽlogramme de $\RR^2$ ayant pour sommets~:
\begin{eqnarray*}
& \left(\frac{1}{4}(2+(1-\aa)\sqrt{5}),\frac{1}{4}(2-(1+\aa)\sqrt{5})\right) & \\
& \left(\frac{1}{4}(2-(1-\aa)\sqrt{5}),\frac{1}{4}(2+(1+\aa)\sqrt{5})\right) & \\
& \left(-\frac{1}{4}(2+(1-\aa)\sqrt{5}),-\frac{1}{4}(2-(1+\aa)\sqrt{5})\right) & \\
& \left(-\frac{1}{4}(2-(1-\aa)\sqrt{5}),-\frac{1}{4}(2+(1+\aa)\sqrt{5})\right).&
\end{eqnarray*}
On voit que $P(\aa)$ est un domaine fondamental pour l'action des translations de $\Sigma({\cal O}_K)$ sur $\RR^2$.
Posons $f_{s,s'}(r,r')=(r^2+s^2)(r'{}^2+s'{}^2)$. 

Supposons pour commencer que 
\begin{equation}
\displaystyle{\left(\frac{s'}{s}\right)^2=\frac{1+\aa}{1-\aa}}\mbox{ avec }|\aa|\leq\frac{1}{\sqrt{5}}.
\label{eq:egalite}\end{equation}
Pour $(r,r')\in P(\aa)$, l'inŽgalitŽ suivante est un exercice ŽlŽmentaire que nous laissons au lecteur~:
\begin{equation}f_{s,s'}(r,r')\leq \left(\frac{5}{16}\right)^2+\frac{9}{8}ss'+s^2s'{}^2.\label{eq:516}\end{equation}
Les inŽgalitŽs concernant $\aa$ sont restrictives. Mais pour tout $(s,s')\in\RR_{>0}^2$, il existe $n\in\ZZ$ unique,
tel que~:
\[\frac{\epsilon^ns'}{\epsilon^{-n}s}=\epsilon^{2n}\frac{s'}{s}=\left(\frac{1+\aa}{1-\aa}\right),\] avec
$|\aa|\leq 1/\sqrt{5}$~: il suffit de remarquer que la fonction $g(\aa)=\displaystyle{\frac{1+\aa}{1-\aa}}$ est strictement
croissante dans $]-\infty,1[$, et que $g(-1/\sqrt{5})=\epsilon'{}^2$, $g(1/\sqrt{5})=\epsilon^2$.

Soit $(z,z')\in{\mathfrak C}$, $z=x+{\rm i}s,z'=x'+{\rm i}s'$. Soient $n,\aa$ comme ci-dessus. Soit $(r,r')\in P(\aa)$ tel que
$(x-r,x'-r')\in\Sigma({\cal O}_K)$. On a~:
\begin{eqnarray*}
|zz'| & = & f_{s,s'}(x,x')\\
& \leq & f_{s,s'}(\epsilon^nr,\epsilon'{}^nr')\\
& \leq & f_{\epsilon^{-n},\epsilon^n}(r,r')\\
& \leq & \left(\frac{5}{16}\right)^2+\frac{9}{8}(\epsilon^{-n}s)(\epsilon^ns')+(\epsilon^{-2n}s^2)(\epsilon^{2n}s'{}^2)\\
& \leq & \left(\frac{5}{16}\right)^2+\frac{9}{8}ss'+(ss')^2.
\end{eqnarray*}
Or, si $(z,z')\in{\mathfrak A}$, alors $|zz'|\geq 1$. Donc 
\[\left(\frac{5}{16}\right)^2+\frac{9}{8}ss'+s^2s'{}^2\geq 1,\] ce qui implique~:
\[ss'\geq\frac{-9+\sqrt{312}}{16}\geq 0.54146,\] d'o l'inŽgalitŽ annoncŽe, et le lemme est dŽmontrŽ. 
On peut dŽmontrer que ${\cal G}$ est un
domaine fondamental pour l'action de $\Upsilon$ sur ${\cal H}^2$, mais nous ne le ferons pas ici car nous n'aurons
pas besoin de ceci~: voir 
les dŽtails dans \cite{Goetzky:Anwendung}.

\subsection{Le diviseur de $\Theta$.}

\begin{Lemme} Le lieu d'annulation de $\Theta$
est l'image de $\Xi_2$ dans $X_\Upsilon$, avec multiplicitŽ $1$.\label{lemme:diviseur_theta}\end{Lemme}

\noindent {\bf DŽmonstration.} Notre technique d'Žtude (due ˆ Gštzky~: \cite{Goetzky:Anwendung}) 
est compl\-tement
ŽlŽmentaire. Plus bas, nous mentionnerons des techniques plus Žlabo\-rŽes, menant ˆ des
rŽsultats plus gŽnŽraux.

Nous considŽrons les facteurs du produit (\ref{eq:produit_10}), et leur comportement dans le domaine
${\cal G}$ du lemme \ref{lemme:domaine_de_goetzky}. Plus prŽcisŽment, nous Žtablissons des minorations de valeurs 
absolues, qui garantissent la non nullitŽ de ces facteurs sur ${\cal G}-\Xi_2$. La modularitŽ de $\Theta$ complŽtera le lemme.

Soit ${\cal F}$ un domaine fondamental pour l'action de $\Upsilon$ sur ${\cal H}^2$ contenu
dans l'ensemble ${\cal G}$ du lemme \ref{lemme:domaine_de_goetzky}.

On Žtudie les $10$ fonctions thtas sur ${\cal F}$ ci-dessus. Posons~:
\[\vartheta_{\alpha,\beta}(\und{z})=\sum_{\nu\in{\cal O}_K}(-1)^{\tr(\nu\beta/\sqrt{5})}\exp\left\{\pi{\rm
i}\tr\left(\left(\nu+\frac{\alpha}{2}\right)^2\frac{\und{z}}{\sqrt{5}}\right)\right\},\] avec $(\alpha,\beta)\in{\cal E}$,
de telle sorte que
\[\Theta(z,z')=\prod_{i=1}^{10}\vartheta_{\alpha_i,\beta_i}(\epsilon z,\epsilon' z').\]
Ces fonctions sont dŽfinies sur ${\cal H}\times{\cal H}^-$, o ${\cal H}^-=\{z\in\CC\mbox{ tel que }\Im(z)<0\}$.
Le groupe $\Upsilon$ agit sur ${\cal H}\times{\cal H}^-$ de manire Žvidente. L'ensemble~:
\[{\cal F}^\flat=\{(z,\bar{z}')\mbox{ avec }(z,z')\in{\cal F}\}\] est un domaine fondamental pour cette action,
et pour tout $(z,z')\in{\cal F}^\flat$ on a $-\Im(z)\Im(z')>0.54146$.
Comme $\epsilon'{}^2\leq -\Im(z)/\Im(z')\leq\epsilon^2$, on a $\Im(z),-\Im(z')>0.454777$ dans ${\cal F}^\flat$. On en dŽduit
des majorations pour $|e_\mu|$, o $\mu\in K_+$, et~:
\[e_\mu=\exp\left\{\pi{\rm i}\tr\left(\mu \frac{\und{z}}{\sqrt{5}}\right)\right\}.\]
Soit $\und{z}\in{\cal F}^\flat$. Si $\mu=1$, comme \[\Im(z)-\Im(z')\geq 2(\Im(z)|\Im(z')|)^{1/2}\geq
2\cdot (0.54146)^{1/2}\approx 1.47167,\] on trouve~:
\begin{equation}
|e_1|=\exp\{-\pi/\sqrt{5}(\Im(z)-\Im(z'))\}<0.126482\mbox{ dans ${\cal F}^\flat$.}\label{eq:un_septieme}
\end{equation}
Naturellement, pour tout $\und{z}\in{\mathfrak B}^\flat\supset{\cal F}^\flat$~:
\begin{equation}|e_\epsilon|,|e_{\epsilon'}|\leq 1,\label{eq:lequ_uno}\end{equation} o ${\mathfrak B}^\flat$ est 
dŽfini de manire Žvidente.
Pour tous $p,q\in\RR$ on a~:
\[(|e_\epsilon|^p-|e_{\epsilon}|^q)(|e_{\epsilon'}|^p-|e_{\epsilon'}|^q)\geq 0.\]
Comme $e_\epsilon e_{\epsilon'}=e_1$, on en dŽduit, pour $\und{z}\in{\cal F}^\flat$~:
\begin{eqnarray}
\lefteqn{|e_{\epsilon}|^p|e_{\epsilon'}|^q+|e_{\epsilon}|^q|e_{\epsilon'}|^p=}\nonumber\\
& = & |e_1|^p+|e_1|^q-(|e_\epsilon|^p-|e_{\epsilon}|^q)(|e_{\epsilon'}|^p-|e_{\epsilon'}|^q)\nonumber\\
& \leq & |e_1|^p+|e_1|^q.\label{eq:utile1}
\end{eqnarray}
Nous Žtudions les fonctions $\vartheta_{\alpha,\beta}(\und{z})$ avec $(\alpha,\beta)\in{\cal E}$.
Observons que~:
\[\exp\left\{-\frac{\pi}{4}{\rm i}\tr\left(\alpha^2\frac{\und{z}}{\sqrt{5}}\right)\right\}\vartheta_{\alpha,\beta}(\und{z})=
\sum_{\nu\in{\cal
O}_K}(-1)^{\tr(\beta\nu/\sqrt{5})} \exp\left\{\pi{\rm i}\tr\left(\nu(\nu+\alpha)\frac{\und{z}}{\sqrt{5}}\right)\right\}.\]
Nous dŽtaillons nos estimations dans le cas $\alpha=0,1$. Si $\alpha=0$ alors $\beta\in\{0,1,\epsilon,\epsilon'\}$.
Si $\alpha=1$ alors $\beta\in\{0,1\}$.
On trouve~:
\begin{eqnarray*}
\lefteqn{e_{\alpha^2}^{-1/4}\vartheta_{\alpha,\beta}(\und{z})=
\displaystyle{\sum_{(n_1,n_2)\in\ZZ^2}\exp\left\{\pi{\rm
i}\tr\left(\frac{\beta(n_1\epsilon+n_2\epsilon')}{\sqrt{5}}\right)\right\}\times}}\\  & &\times
e_1^{(n_1-n_2)^2}e_\epsilon^{n_1(n_1+\alpha)}e_{\epsilon'}^{n_2(n_2+\alpha)},
\end{eqnarray*}
car $(\epsilon,\epsilon')$ est une base de ${\cal O}_K$. Donc
\begin{eqnarray*}
\lefteqn{|e_{\alpha^2}^{-1/4}\vartheta_{\alpha,\beta}(\und{z})-(1+\alpha)|\leq}\\
&\leq &  -(1+\alpha) +\frac{1}{2}\sum|e_1|^{(n_1-n_2)^2}\times\left(|e_\epsilon|^{n_1(n_1+\alpha)}|e_{\epsilon'}|^{n_2(n_2+\alpha)}+\right.\\
& & \left.|e_\epsilon|^{n_2(n_2+\alpha)}|e_{\epsilon'}|^{n_1(n_1+\alpha)}\right)\\
& \leq &  -(1+\alpha) +\frac{1}{2}\sum|e_1|^{(n_1-n_2)^2}(|e_1|^{n_1(n_1+\alpha)}+|e_1|^{n_2(n_2+\alpha)})\\
& \leq & -(1+\alpha) +\sum_{m=-\infty}^\infty|e_1|^{m^2}\sum_{n=-\infty}^\infty|e_1|^{n(n+\alpha)}, 
\end{eqnarray*} en utilisant les inŽgalitŽs (\ref{eq:utile1}).
L'inŽgalitŽ (\ref{eq:un_septieme}) implique~:
\begin{eqnarray*}
|e_{\alpha^2}^{-1/4}\vartheta_{\alpha,\beta}(\und{z})|& \geq &
2(1+\alpha)-\sum_{m=-\infty}^\infty|e_1|^{m^2}\sum_{n=-\infty}^\infty|e_1|^{n(n+\alpha)}\\ & > & 0.4288.
\end{eqnarray*}
Donc, si $\alpha\in\{0,1\}$, $\vartheta_{\alpha,\beta}$ ne s'annule pas dans ${\cal F}^\flat$, ce qui donne la
non nullitŽ dans ${\cal F}^\flat$ de six facteurs du produit dŽfinissant $\Theta$.

Si $\alpha\in\{\epsilon,\epsilon'\}$ (quatre cas), on voit que $\vartheta$ s'annule dans $\{(z,z)\}\subset{\cal H}^2$.
La technique expliquŽe ci-dessus, opportunŽment modifiŽe, implique que pour $\und{z}\in{\cal F}^\flat$~:
\begin{eqnarray*}
\left|\frac{\vartheta_{\alpha,\beta}(\und{z})}{2e_{\alpha^2}^{1/4}(1+(-1)^{\tr(\alpha'\beta/\sqrt{5})}e_{\alpha'})}\right| & > & \frac{1}{100},
\end{eqnarray*}
donc ne s'annule pas. Nous prŽfŽrons ne pas dŽtailler ce deuxime cas, car les
techniques sont les mmes~: pour des dŽtails, voir \cite{Goetzky:Anwendung}. Ceci 
implique que la restriction de $\Theta$ ˆ ${\cal F}$ s'annule seulement dans $\Xi_2$ avec multiplicitŽ $1$~: le lemme
est dŽmontrŽ.

La thŽorie des produits de Borcherds (\cite{Borcherds:Automorphic1}, \cite{Borcherds:Automorphic2})
permet de retrouver ces
rŽsultats, et bien d'autres propriŽtŽs que nous ne pouvons pas  dŽcrire ici. 
Dans \cite{Bruinier:Borcherds}
on montre que~:
\[\Theta(\und{z})=64\exp\left\{2\pi{\rm i}\left(\frac{\epsilon z}{\sqrt{5}}-\frac{\epsilon'z'}{\sqrt{5}}
\right)\right\}\prod_{{\tiny \begin{array}{c}\nu\in{\cal O}_K^*\\
\epsilon\nu'-\epsilon'\nu> 0
\end{array}}}(1-e^{2\pi{\rm i}\tr(\nu \und{z})})^{s(5\nu\nu')a(5\nu\nu')},\]
o $\sum_{n=0}^\infty a(n)q^n$ est la sŽrie de Fourier d'une certaine forme modulaire d'une variable complexe,
et $s:\NN\rightarrow\ZZ$ est une certaine fonction. Tout ceci est compltement explicite~: voir les
dŽtails dans \cite{Bruinier:Borcherds}.

\subsection{Structure d'un certain anneau de formes modulaires.}

\noindent {\bf DŽfinition.} Soit $f$ une forme modulaire de Hilbert de poids parallle $r$ pour un certain
groupe modulaire de Hilbert
$\Gamma$. On dit que $f$ est {\em symŽtrique} (resp. {\em antisymŽtrique}) si $f(z,z')=f(z',z)$ (resp. $f(z,z')=-f(z',z)$).

\medskip

L'existence de la fonction $\Theta$ permet de faire une description explicite de l'anneau de certaines 
formes modulaires pour $\Upsilon$.

\medskip

Avant de continuer, nous faisons une remarque. Nous voulons dŽcrire encore plus en dŽtail la structure des espaces vectoriels de 
formes modulaires de Hilbert dans le cas particulier de $K=\QQ(\sqrt{5})$. Nous allons utiliser des notations
particulires ˆ ce cas. Ainsi, lorsqu'il s'agit de formes modulaires pour le groupe modulaire de Hilbert associŽ ˆ
ce corps de nombres, et {\em seulement dans ce cas}, nous posons 
(\footnote{Ces notations sont pratiquement celles utilisŽes par Gundlach et Resnikoff, dans leurs travaux sur ces formes modulaires, ˆ
part le fait qu'ils Žcrivent $\chi_5^-$, et nous $\chi_5$.
On remarque que $\varphi_2,\chi_5,\chi_6$ ont leur sŽries de Fourier ˆ coefficients entiers rationnels.})~:
\begin{eqnarray*}
\varphi_2 & = & \zeta_K(2)^{-1}E_{2,{\cal O}_K},\\
\chi_5 & = & 2^{-5}\Theta,\\
\chi_6 & = & \frac{67}{21600}(\varphi_2^3-\zeta_K(6)^{-1}E_{6,{\cal O}_K}).
\end{eqnarray*}
\begin{Theoreme} 
Toute forme modulaire symŽtrique de poids parallle $2r$ avec $r\in\NN$ pour $\Upsilon$ est un polyn™me isobare 
en les trois formes modulaires
$\varphi_2,\chi_5^2,\chi_6$.\label{lemme:structure1}\end{Theoreme}

\medskip

Nous verrons que les formes modulaires $\varphi_2,\chi_5,\chi_6$ sont algŽbriquement 
indŽ\-pendantes sur $\CC$ (corollaire \ref{lemme:indep_alg}).

\medskip

\noindent {\bf DŽmonstration.} On commence par observer que $\chi_6$ est une forme parabo\-lique de poids $(6,6)$, 
et donc $\chi_6(z,z)$ est
une forme modulaire pour $\SL_2(\ZZ)$ de poids $12$, Žgale ˆ $c\Delta(z)$, pour $c\in\CC$. 
Comme $E_{2,{\cal O}_K},E_{6,{\cal O}_K}$ sont
symŽtriques, $\chi_6$ est symŽtrique.

On dŽtermine $c$ en calculant les premiers coefficients de Fourier de $\chi_6$. 
Ceci est possible
car les coefficients de Fourier des sŽries d'Eisenstein sont explicites (cf. sous-paragraphe 
\ref{section:recette}). On trouve~:
\[\chi_6(z,z)=\Delta(z).\]
On dit que $\chi_6$ est un {\em relevŽ} de $\Delta$.

Nous avons dŽjˆ vu que $\varphi_2(z,z)=E_4(z)$. Comme l'anneau des formes modulaires de poids divisible par $4$
pour le groupe $\SL_2(\ZZ)$ est Žgal ˆ~:
\[\CC[E_4(z),\Delta(z)],\]
on voit que toute forme modulaire de poids divisible par $4$ peut tre relevŽe en une forme modulaire symŽtrique
de poids $(2r,2r)$ pour $\Upsilon$. Nous avons donc un isomorphisme d'anneaux~:
\[\mu:\CC[E_4(z),\Delta(z)]\rightarrow\CC[\varphi_2(\und{z}),\chi_6(\und{z})],\]
inverse du morphisme de restriction ˆ $\Xi_2=\{(z,z)\}$, ce qui prouve entre autres 
que les formes modulaires $\varphi_2$ et $\chi_6$ sont algŽbriquement indŽpendantes.

Soit maintenant $F$ une forme modulaire symŽtrique de poids $(2r,2r)$ pour $\Upsilon$. Si 
$F\in\CC[\varphi_2,\chi_6]$
nous n'avons rien ˆ dŽmontrer. Supposons que
\[H:=F-\mu(F(z,z))\not=0.\]
La forme modulaire $H$ s'annule dans $\Xi_2$ par construction. L'Žtude de la fonction 
$\chi_5=2^{-5}\Theta$ que nous avons faite
implique que \[H_1:=\frac{H}{\chi_5^2}\] est une forme modulaire symŽtrique de poids 
$(2r-10,2r-10)$. La dŽmons\-tration se complte par rŽcurrence sur $r$. 

\section{PropriŽtŽs diffŽrentielles de formes modulaires de Hilbert.}

Dans le paragraphe \ref{section:prelude}, nous avons vu comment on dŽtermine, en utilisant certains crochets de Rankin-Cohen, 
la structure de 
l'anneau engendrŽ par toutes les formes modulaires elliptiques. 

On peut Žtendre la dŽfinition de crochet de Rankin-Cohen aux formes modulaires de Hilbert. Par
exemple, ceci est  Žcrit dans \cite{Lee:Hilbert}. Ainsi on associe, ˆ deux formes modulaires de Hilbert $F,G$ 
de poids $\und{f},\und{g}\in\NN^n$ (pour un corps de
nombres totalement rŽel de degrŽ $n$ sur $\QQ$), et ˆ un
$n$-uplet de nombres entiers positifs $\und{s}$, une forme parabolique $[F,G]_{\und{s}}$ de poids $\und{f}+\und{g}+2\und{s}$~:
\begin{eqnarray*}
[F,G]_{\und{s}} & = & \frac{1}{(2\pi{\rm i})^{|\und{s}|}}\sum_{r_1=0}^{s_1}\cdots\sum_{r_n=0}^{s_n}
(-1)^{|\und{r}|}\frac{\partial^{|\und{r}|}F}{\partial z_1^{r_1}\cdots
\partial z_n^{r_n}}\frac{\partial^{|\und{s}-\und{r}|}G}{\partial z_1^{s_1-r_1}\cdots
\partial z_n^{s_n-r_n}}\times\\
& &\prod_{i=1}^n\binomial{f_i+s_i-1}{s_i-r_i}\binomial{g_i+s_i-1}{r_i},
\end{eqnarray*}
o $|(a_1,\ldots,a_n)|=a_1+\cdots+a_n$. 

Pour $x\in\ZZ$, nous Žcrivons $x_i={}^t(0,\ldots,0,x,0,\ldots,0)\in\ZZ^n$ avec le coefficient $x$ ˆ la $i$-me place.
Voici l'exemple le plus simple de crochet de Rankin-Cohen. La fonction~:
\[
[F,G]_{1_i} = \frac{1}{(2\pi {\rm i})}\left(g_iG\frac{\partial F}{\partial z_i}-f_iF\frac{\partial G}{\partial z_i}\right)\]
est une forme parabolique de poids $\und{f}+\und{g}+2_i$ (crochet de Rankin).

Posons~:
\[\Lambda_{i,k}F=\frac{1}{2f_if_k}[F,F]_{1_i+1_k},\quad\Pi_{i}F=\frac{1}{(f_i+1)}[F,F]_{2_i}.\]
Ces opŽrateurs diffŽrentiels s'Žcrivent explicitement de la manire suivante~:
\begin{eqnarray*} 
(2\pi {\rm i})^2\Pi_iF & := & f_iF\frac{\partial^2 F}{\partial z_i^2}-(f_i+1)\left(\frac{\partial F}{\partial z_i}\right)^2,\\
(2\pi {\rm i})^2\Lambda_{i,k}F & := & F\frac{\partial^2 F}{\partial z_i\partial z_k}-\frac{\partial F}{\partial z_i}
\frac{\partial F}{\partial z_k}.
\end{eqnarray*}
On vŽrifie que~:
\begin{eqnarray*}
\Pi_i:M_{\und{f}}(\Gamma) & \rightarrow & S_{2\und{f}+4_i}(\Gamma),\\
\Lambda_{i,k}:M_{\und{f}}(\Gamma) & \rightarrow & S_{2\und{f}+2_i+2_k}(\Gamma).
\end{eqnarray*}
Voici maintenant un lemme ŽlŽmentaire qui dŽcrit les propriŽtŽs d'annulation des images des opŽrateurs diffŽrentiels $[\cdot,\cdot]_{1_i},\Pi_i,\Lambda_{i,k}$.
\begin{Lemme}
Soient $F,G$ deux formes modulaires de Hilbert de $n$ variables complexes, non constantes, de poids $\und{f}$ et $\und{g}$. Alors~:
\begin{enumerate}
\item $\Lambda_{i,k}F\not=0$.
\item $\Pi_iF\not=0$.
\item $[F,G]_{1_i}=0$ si et seulement si la fonction $F^{g_i}/G^{f_i}$ est constante.
\end{enumerate}
\label{lemme:annulation_operateurs}
\end{Lemme}
\noindent {\bf DŽmonstration.} (1). On a \[\Lambda_{i,k}F=\frac{1}{(2\pi{\rm i})^2}F^2
\displaystyle{\frac{\partial^2}{\partial z_i\partial z_k}\log F}.\]
Nous pouvons supposer que $i=1,k=2$. Si $\Lambda_{i,k}F=0$, alors 
la fonction mŽromorphe $H=F^{-1}\partial F/\partial z_2$ ne dŽpend pas de la variable $z_1$. On a de plus
\[H(z_2+\sigma_2(\nu),\ldots,z_n+\sigma_n(\nu))=H(z_2,\ldots,z_n)\quad\mbox{ pour }\nu\in{\cal O}_K.\]
Soit $\pi_1$ la projection $\RR^n\rightarrow\RR^{n-1}$ obtenue en
effaant le premier coefficient. L'ensemble $\pi_1(\Sigma({\cal O}_K))$ est dense dans $\RR^{n-1}$. En particulier $H$
est constante dans $E=\pi_{1}\und{z}+\RR^{n-1}$. Le plus petit sous-ensemble analytique de ${\cal H}^{n-1}$ contenant
$E$ est ${\cal H}^{n-1}$~: donc $H$ est constante dans ${\cal H}^{n-1}$, non nulle car $F$ est non constante.
Donc~:
\[\frac{\partial F}{\partial z_2}=\lambda F,\] avec $\lambda\not=0$. Mais on voit clairement que ceci est incompatible
avec la modularitŽ de $F$.

\medskip 

\noindent (2). Si $\Pi_{i}F=0$ et $F$ est non constante, alors
on peut Žcrire $F(\und{z})=(gz_i+h)^{-f_i}$, o $g,h$ sont deux fonctions ne dŽpendant pas de $z_i$. On vŽrifie directement qu'une 
telle fonction ne peut pas tre une forme modulaire. On peut aussi utiliser les sŽries de Fourier du lemme \ref{lemme:injectivite}.

\medskip

\noindent (3). On a~:
\[[F,G]_{1_i}=-\frac{1}{2\pi{\rm i}}FG\frac{\partial}{\partial z_i}\log\left(\frac{G^{f_i}}{F^{g_i}}\right).\] Si $G^{f_i}/F^{g_i}\in\CC^\times$,
alors $[F,G]_{1_i}=0$. Supposons maintenant que $[F,G]_{1_i}=0$~: on a alors que $H:=G^{f_i}/F^{g_i}$ ne dŽpend pas de 
la variable $z_i$. D'autre part, c'est un quotient de formes modulaires, donc il satisfait~:
\begin{eqnarray*}
H(z_2,\ldots,z_n) & = & H(z_2+\sigma_2(\nu),\ldots,z_n+\sigma_n(\nu)),
\end{eqnarray*}
pour tout $\nu\in{\cal H}^n$ et $z_2,\ldots,z_n\in{\cal H}$, lorsque ces valeurs sont dŽfinies. Comme au point (1),
$H$ est constante dans ${\cal H}^{n-1}$~: la dŽmonstration du lemme \ref{lemme:annulation_operateurs} est complte. 
On remarque en particulier que $H$ est une fonction modulaire, et que les poids de $F$ et $G$ satisfont
$\und{f}\in\QQ^\times\und{g}$. 

\medskip

\noindent {\bf Exemple 1.}
Comme il existe toujours une forme modulaire non nulle de poids $(2,\ldots,2)$ pour $\Gamma$ si $n>1$, on peut construire
dans certains cas des formes modulaires non nulles de poids non parallle \[(t_1,\ldots,t_n)\in(2\ZZ)^n,\] avec
$t_i\geq 2$ pour tout $i=1,\ldots,n$. Par exemple, si $K$ est quadratique, $\Pi_1E_{2,{\cal O}_K}$ est une forme
parabolique de poids $(8,4)$ qui est non nulle (d'aprs le lemme \ref{lemme:annulation_operateurs}), 
et constitue un premier specimen 
de forme parabolique de poids non parallle. 

\medskip

\noindent {\bf Exemple 2.} Si $[K:\QQ]=2$, alors les opŽrateurs diffŽrentiels $\Lambda=\Lambda_{1,1}$ et~:
\begin{eqnarray*}
\Pi :F & \mapsto & (\Pi_1F)(\Pi_2F)
\end{eqnarray*}
agissent sur l'anneau graduŽ engendrŽ par les formes modulaires de Hilbert de poids parallle.
Par de simples arguments d'algbre linŽaire (basŽs sur des calculs des premiers
coefficients de Fourier de formes modulaires), 
on peut dŽterminer des relations diffŽrentielles entre gŽnŽrateurs, induites
par $\Lambda$ et $\Pi$. Par exemple, si $K=\QQ(\sqrt{5})$, on trouve~:
\begin{eqnarray}
\Lambda\varphi_2 & = & 24 \chi_6,\nonumber\\
\Lambda\chi_6 & = & \frac{1}{20}\varphi_2(\varphi_2\chi_5^2 - \chi_6^2),\label{eq:systeme2}\\  
\Lambda\chi_5^2 & = & \frac{1}{10}\chi_5^2\left(\chi_6^2-\varphi_2\chi_5^2\right).\nonumber
\end{eqnarray}

\medskip

\noindent {\bf \'Equations diffŽrentielles.}
On remarque que les formes modulaires $\chi_6$ et $\chi_5^2$ peuvent tre construites par
application d'opŽrateurs diffŽrentiels, ˆ partir de $\varphi_2$. Gr‰ce ˆ ces identitŽs, on peut aussi expliciter des
Žquations aux dŽrivŽes partielles ayant des formes modulaires de Hilbert comme solutions~: en voici un exemple.

On a le systme diffŽrentiel (\ref{eq:systeme2}), et 
on sait que $\Pi \varphi_2\in M_{\und{12}}(\Gamma)$. D'aprs le thŽorme \ref{lemme:structure}, $\Pi \varphi_2$
est un polyn™me isobare en $\varphi_2,\chi_6$ et $\chi_5^2$. La rŽsolution d'un systme linŽaire
explicite (dans ce cas, sept Žquations
et sept inconnues) engageant les coefficients de Fourier de ces formes modulaires implique~:
\begin{equation}\Pi \varphi_2=576(9\chi_6^2-5\varphi_2\chi_5^2).\label{eq:deri2}\end{equation}
En Žliminant $\chi_5^2$ dans (\ref{eq:deri2}) et la troisime Žquation de (\ref{eq:systeme2}), on trouve~:
\begin{equation}
100(\Lambda\circ \Lambda)\varphi_2+\varphi_2\Pi \varphi_2-4\varphi_2(\Lambda\varphi_2)^2=0.\label{eq:equadiff}
\end{equation}
Cette Žquation diffŽrentielle a la propriŽtŽ intŽressante suivante (cf. \cite{Resnikoff:Automorphic})~: 
$\varphi_2$ est la seule solution de cette Žquation
qui soit une forme modulaire de Hilbert pour $\QQ(\sqrt{5})$ dont la sŽrie de Fourier a les termes indŽxŽs par les plus petites traces~:
\[1+120\exp\left\{2\pi{\rm i}\tr\left(\frac{\epsilon}{\sqrt{5}}z\right)\right\}+\cdots\]

\subsection{OpŽrateurs multilinŽaires.} 

On peut dŽfinir des
opŽrateurs multilinŽaires agissant sur certains espaces de formes modulaires
de Hilbert (et mme sur des espaces de polyn™mes non isobares en des formes modulaires), de
la manire suivante. Soit
${\cal I}$ un sous-ensemble non vide de
$\{1,\ldots,n\}$~: pour simplifier, choisissons ${\cal I}=\{1,\ldots,m\}$.
Soit $K$ un corps de nombres totalement rŽel de degrŽ $n$, soit $\Gamma$ le groupe modulaire de Hilbert associŽ. 

\medskip

\noindent {\bf DŽfinition.} On dit qu'une forme modulaire de Hilbert $F$ pour $\Gamma$ est de {\em poids
${\cal I}$-parallle}, si son poids $\und{f}=(f_1,\ldots,f_n)$ est tel que $f_i=f$ pour tout $i\in{\cal I}$.
On dit dans ce cas que $f$ est le ${\cal I}$-poids de $F$.

\medskip

Ainsi, une forme modulaire de Hilbert de poids parallle $f$ pour $\Gamma$ est une forme modulaire de
Hilbert de poids ${\cal I}$-parallle, avec ${\cal I}=\{1,\ldots,n\}$. 

\begin{Lemme}
Soient $F_1,\ldots,F_{m+1}$ des formes modulaires de poids ${\cal I}$-parallles $r_1,\ldots,r_{m+1}$. La fonction~:
\[
\langle F_1,\ldots,F_{m+1}\rangle_{\und{1}}:=\displaystyle{\det\left(\begin{array}{cccc}
r_1F_1 & r_2F_2 & \cdots & r_{m+1}F_{m+1}\\ & & & \\
\frac{\partial F_1}{\partial z_1} & \frac{\partial F_2}{\partial z_1} & \cdots & \frac{\partial F_{m+1}}{\partial z_{1}} \\
& & & \\
\frac{\partial F_1}{\partial z_2} & \frac{\partial F_2}{\partial z_2} & \cdots & \frac{\partial F_{m+1}}{\partial z_{2}} \\
& & & \\
\vdots & \vdots & & \vdots \\
& & & \\
\frac{\partial F_1}{\partial z_m} & \frac{\partial F_2}{\partial z_m} & \cdots & \frac{\partial F_{m+1}}{\partial z_m}
\end{array}\right)}
\]
est une forme modulaire de Hilbert de poids ${\cal I}$-parallle $r_1+\cdots+r_{m+1}+2$.\label{lemme:multi}\end{Lemme}

\noindent {\bf DŽmonstration}. 
On note
${\mathfrak S}_m$ le groupe des permutations de ${\cal I}$, $\epsilon(\sigma)$ la signature d'une permutation $\sigma\in{\mathfrak S}_m$,
on pose~:
\[[F_1,\ldots,F_{m+1}]_{\und{1}}:=\sum_{\sigma\in{\mathfrak
S}_{m}}\epsilon(\sigma)[F_1,[F_2,\ldots,[F_m,F_{m+1}]_{1_{\sigma(m)}}\ldots]_{1_{\sigma(2)}}]_{1_{\sigma(1)}}.\]
La preuve est une consŽquence immŽdiate de la formule suivante~:
\begin{equation}[F_1,\ldots,F_{m+1}]_{\und{1}}=(2\pi{\rm i})^{-m}\Phi_m\langle F_1,\ldots,F_{m+1}\rangle_{\und{1}},\label{eq:formulotta}\end{equation}
avec
\[\Phi_m=(r_2+\cdots+r_{m+1})(r_3+\cdots+r_{m+1})\cdots(r_m+r_{m+1}),\]
(produit ayant $m-1$ facteurs)
car le terme de gauche dans (\ref{eq:formulotta}) est clairement une forme modulaire de Hilbert de poids ${\cal I}$-parallle $r_1+\cdots+r_{m+1}+2$,
d'aprs les propriŽtŽs de base du crochet de Rankin. 

Par exemple, si $F_1,F_2$ sont des formes modulaires elliptiques
de poids $r_1,r_2$, la formule (\ref{eq:formulotta}) nous donne tout simplement un crochet de Rankin~:
\begin{eqnarray*}
[F_1,F_2]_{\und{1}}
&=&(2\pi{\rm i})^{-1}\det\left(\begin{array}{cc}r_1F_1 & r_2F_2\\ & \\ \displaystyle{\frac{\partial F_1}{\partial z_1}}&
\displaystyle{\frac{\partial F_2}{\partial z_1}}\end{array}\right)
\end{eqnarray*}
Si $G_1,G_2,G_3$ sont des formes modulaires de Hilbert de deux variables complexes, 
de poids parallles $s_1,s_2,s_3$, alors la formule (\ref{eq:formulotta})
implique~:
\begin{eqnarray}
{[}G_1,G_2,G_3{]}_{\und{1}}&:=&[G_1,[G_2,G_3]_{1_2}]_{1_1}-[G_1,[G_2,G_3]_{1_1}]_{1_2}\nonumber\\
& = & (2\pi{\rm i})^{-2}(s_2+s_3)\det\left(\begin{array}{ccc}s_1G_1 & s_2G_2& s_3G_3\\ & \\ \displaystyle{\frac{\partial G_1}{\partial z_1}}&
\displaystyle{\frac{\partial G_2}{\partial z_1}}&\displaystyle{\frac{\partial G_3}{\partial z_1}}\label{eq:m3}\\
& & \\ \displaystyle{\frac{\partial G_1}{\partial z_2}}&
\displaystyle{\frac{\partial G_2}{\partial z_2}}&\displaystyle{\frac{\partial G_3}{\partial z_2}}\end{array}\right)
\end{eqnarray}
La dŽmonstration de la formule (\ref{eq:formulotta}) repose 
sur l'ŽgalitŽ suivante, dont la dŽmonstration est ŽlŽmentaire et laissŽe au lecteur~:
\begin{eqnarray}
\lefteqn{\sum_{s=1}^m(-1)^s\frac{\partial}{\partial z_s}\langle F_2,\ldots,F_{m+1}\rangle_{\und{1}-1_s}=}
\label{eq:formulacchia}\\
& = & (r_2+\cdots+r_{m+1})\det\left(\frac{\partial
F_i}{\partial z_j}\right)_{{\tiny\begin{array}{c} i=2,\ldots,m+1\\ j=1,\ldots,m\end{array}}}
\nonumber,
\end{eqnarray} o $\und{1}-1_s=(\underbrace{1,\ldots,1}_{s-1\mbox{ \tiny termes }},0,1,\ldots,1)$
et o
\[
\langle F_2,\ldots,F_{m+1}\rangle_{\und{1}-1_s}:=\displaystyle{\det\left(\begin{array}{cccc}
r_2F_2 & r_3F_3 & \cdots & r_{m+1}F_{m+1}\\ & & & \\
\frac{\partial F_2}{\partial z_1} & \frac{\partial F_3}{\partial z_1} & \cdots & \frac{\partial F_{m+1}}{\partial z_{1}} \\
& & & \\
\vdots & \vdots & & \vdots \\
& & & \\
\frac{\partial F_2}{\partial z_{s-1}} & \frac{\partial F_3}{\partial z_{s-1}} & \cdots & \frac{\partial F_{m+1}}{\partial z_{s-1}} \\
& & & \\
\frac{\partial F_2}{\partial z_{s+1}} & \frac{\partial F_3}{\partial z_{s+1}} & \cdots & \frac{\partial F_{m+1}}{\partial z_{s+1}} \\
& & & \\
\vdots & \vdots & & \vdots \\
& & & \\
\frac{\partial F_2}{\partial z_m} & \frac{\partial F_3}{\partial z_m} & \cdots & \frac{\partial F_{m+1}}{\partial z_m}
\end{array}\right)}
\]
On dŽmontre (\ref{eq:formulotta}) par rŽcurrence
sur
$m$. Pour $m=1$, l'identitŽ est triviale, comme nous l'avons dŽjˆ remarquŽ. Supposons que $m$ soit plus grand~; 
comme il existe une bijection entre ${\mathfrak S}_m$ et $\cup_{s=1}^m{\mathfrak S}_s'$
o ${\mathfrak S}_s'$ est l'ensemble des fonctions bijectives $\{2,\ldots,m\}\rightarrow\{1,\ldots,s-1,s+1,\ldots,m\}$,
on a~:
\begin{eqnarray*}
[F_1,\ldots,F_m]_{\und{1}} & = & \sum_{s=1}^m(-1)^{s+1}[F_1,\sum_{\tilde{\sigma}\in{\mathfrak
S}_s'}\epsilon(\tilde{\sigma})[F_2,[\cdots[F_m,F_{m+1}]_{1_{\tilde{\sigma}(m)}}\cdots]_{1_{\tilde{\sigma}(2)}}]_{1_s},
\end{eqnarray*}
($\epsilon$
dŽsigne le prolongement de l'application signature ˆ ${\mathfrak
S}_s'$). Donc~:
\begin{eqnarray}
[F_1,\ldots,F_m]_{\und{1}} & = & [F_1,\sum_{s=1}^m(-1)^s[F_2,\ldots,F_{m+1}]_{\und{1}-1_s}]_{1_s}\nonumber\\
& = & (2\pi{\rm i})^{m-1}\Phi'\sum_{s=1}^m(-1)^s[F_1,\langle F_2,\ldots,F_{m+1}\rangle_{\und{1}-1_s}]_{1_s}\label{a}\\
& = & (2\pi{\rm i})^{m}\Phi'\sum_{s=1}^m(-1)^s\left(r_1F_1\frac{\partial}{\partial z_s}\langle F_2,\ldots,F_{m+1}\rangle_{\und{1}-1_s}\right.\nonumber\\
& &\left.-(r_2+\cdots+r_{m+1})\langle F_2,\ldots,F_{m+1}\rangle_{\und{1}-1_s}\frac{\partial F_1}{\partial z_s}\right)\nonumber\\
& = & (2\pi{\rm i})^{m}\Phi_m\left(r_1F_1\det\left(\frac{\partial F_i}{\partial z_j}\right)_{i,j}\right.\label{b}\\
& &\left.-\frac{\partial F_1}{\partial z_s}\langle F_2,\ldots,F_{m+1}\rangle_{\und{1}-1_s}\right)\nonumber\\
& = & (2\pi{\rm i})^{m}\Phi_m\langle F_1,\ldots, F_{m+1}\rangle_{\und{1}},\nonumber
\end{eqnarray} o le symbole $[F_2,\ldots,F_{m+1}]_{\und{1}-1_s}$ a une signification Žvidente,
et $\Phi'=(r_3+\cdots+r_{m+1})\cdots(r_m+r_{m+1})$
(en (\ref{a}) nous avons appliquŽ l'hypothse de rŽcurrence, et en (\ref{b}) nous avons appliquŽ l'ŽgalitŽ
(\ref{eq:formulacchia})). La preuve du lemme est terminŽe.

\medskip

\noindent {\bf Remarque.} 
Les opŽrateurs $[\cdots]_{\und{1}}$ sont dŽfinis en combinant entre eux
des crochets de Rankin, et ˆ premire vue, on est temptŽs de dire que ce sont des polyn™mes diffŽrentiels (voir le paragraphe \ref{section:differentiels})
d'ordre
$m$ dont  toutes leurs propriŽtŽs algŽbriques sont consŽquence de la
thŽorie des crochets de Rankin-Cohen. Cependant, le lemme \ref{lemme:multi} implique qu'il s'agit en rŽalitŽ
d'opŽrateurs diffŽrentiels d'ordre $1$, 
en particulier multilinŽaires. Nous verrons que cette pŽculiaritŽ fait des ces opŽrateurs des 
objets intŽressants, dont l'Žtude sera poursuivie dans des autres travaux.

\medskip

Dans la suite de ce texte, nous Žtudions ces opŽrateurs uniquement dans le cas o ${\cal I}=\{1\}$ et dans
le cas o ${\cal I}=\{1,2\}$ et le nombre de variables complexes est deux, 
(dans ce dernier cas, les opŽrateurs sont bien dŽfinis uniquement sur des formes
modulaires de Hilbert de poids parallle). Voici maintenant une proposition qui dŽcrit les propriŽtŽs de base dans le cas ${\cal I}=\{1,2\}$.

\begin{Proposition}
Soit $K$ un corps de nombres quadratique rŽel, et $\Gamma$ le groupe modulaire de Hilbert associŽ.
Soient $F,G,H$ des formes modulaires de Hilbert pour $\Gamma$ de poids parallle $f,g,h$.
La forme modulaire $[F,G,H]:=[F,G,H]_{\und{1}}$ est de poids parallle $f+g+h+2$.
Si $F$ est antisymŽtrique et $G,H$ sont 
symŽtriques, ou si $F,G,H$ sont antisymŽtriques, alors $[F,G,H]$ est symŽtri\-que. Si $F$ est symŽtrique et 
$G,H$ sont anti\-symŽ\-triques, ou si $F,G,H$ sont symŽtri\-ques, alors $[F,G,H]$ est antisymŽtrique.

De plus, les deux conditions suivantes sont Žquivalentes.
\begin{enumerate}
\item La forme modulaire $[F,G,H]$ est non nulle.
\item Les fonctions $F,G,H$ sont algŽbriquement indŽpendantes sur $\CC$.
\end{enumerate}
\label{lemme:diff}\end{Proposition}
\noindent {\bf DŽmonstration.} D'aprs le lemme \ref{lemme:multi}, $[F,G,H]$ est une forme modulaire de Hilbert de poids
$f+g+h+2$, et on a
\[[F,G,H]=\frac{1}{(2\pi{\rm i})^2}(g+h)\displaystyle{\det\left(\begin{array}{ccc}
fF & gG & hH \\ & & \\ \displaystyle{\frac{\partial F}{\partial z}} & \displaystyle{\frac{\partial G}{\partial z}} & 
 \displaystyle{\frac{\partial H}{\partial z}} \\  & & \\
\displaystyle{\frac{\partial F}{\partial z'}} & \displaystyle{\frac{\partial G}{\partial z'}} & 
 \displaystyle{\frac{\partial H}{\partial z'}}
\end{array}\right)},\] ou de manire Žquivalente~:
\[[F,G,H]=\frac{1}{(2\pi{\rm i})^2}(g+h)\left(fF(G\wedge H)+hH(F\wedge G)-gG(F\wedge H)
\right),\] avec
\[\alpha\wedge\beta=\det\sqm{\displaystyle{\frac{\partial \alpha}
{\partial z}}}{\displaystyle{\frac{\partial \alpha}{\partial z'}}}
{\displaystyle{\frac{\partial \beta}{\partial z}}}
{\displaystyle{\frac{\partial \beta}{\partial z'}}}.\]
Soit ${\cal C}:\CC^2\rightarrow\CC^2$ la symŽtrie ${\cal C}(z,z')=(z',z)$. Pour toute fonction
analytique $A$, on a~:
\[\frac{\partial}{\partial z}(A\circ{\cal C})=\left(\frac{\partial A}{\partial z'}\right)\circ{\cal C}\mbox{ et }
\frac{\partial}{\partial z'}(A\circ{\cal C})=\left(\frac{\partial A}{\partial z}\right)\circ{\cal C}.\]
On en dŽduit~:
\begin{eqnarray*}
(\alpha\wedge\beta)\circ{\cal C} & = & \left(\frac{\partial\alpha}{\partial z}\frac{\partial\beta}{\partial z'}
-\frac{\partial\alpha}{\partial z'}\frac{\partial\beta}{\partial z}\right)\circ{\cal C}\\
& = & \left(\frac{\partial\alpha}{\partial z}\circ{\cal C}\right)\left(\frac{\partial\beta}{\partial z'}\circ{\cal C}\right)
-\left(\frac{\partial\alpha}{\partial z'}\circ{\cal C}\right)\left(\frac{\partial\beta}{\partial z}\circ{\cal C}\right)\\
& = & \frac{\partial}{\partial z'}(\alpha\circ{\cal C})\frac{\partial}{\partial z}(\beta\circ{\cal C})-
\frac{\partial}{\partial z}(\alpha\circ{\cal C})\frac{\partial}{\partial z'}(\beta\circ{\cal C}).
\end{eqnarray*} 
Cette dernire quantitŽ est Žgale ˆ $-(\alpha\wedge\beta)$ si $\alpha,\beta$ symŽtriques ou antisymŽtriques,
et Žgale ˆ $\alpha\wedge\beta$ si $\alpha$ symŽtrique (resp. antisymŽtrique) et $\beta$
antisymŽtrique (resp. symŽtrique).

DŽcrivons maintenant les conditions d'annulation de $[F,G,H]$.
On vŽrifie directement que si $F,G,H$ sont multiplicativement dŽpendantes modulo $\CC^\times$, alors $[F,G,H]=0$.
Montrons que la condition (1) implique la condition (2). 
Si $M$ est non nulle, alors $F,G,H$ sont multiplicativement indŽpendantes modulo $\CC^\times$, donc $H^f/F^h$ et $H^g/G^h$ sont aussi 
multiplicativement indŽpen\-dantes modulo $\CC^\times$. On peut appliquer le thŽorme 3 p. 253 de \cite{Ax:Shanuel}, et on
trouve que  le corps de fonctions
\[{\mathfrak K}:=\CC\left(\log\left(\frac{H^f}{F^h}\right),\log\left(\frac{H^g}{G^h}\right),\frac{H^f}{F^h},\frac{H^g}{G^h}\right)\]
a degrŽ de transcendance minorŽ par $2+\mu$, o $\mu$ est le rang de la matrice jacobienne logarithmique
\[\displaystyle{J=\frac{1}{(2\pi{\rm i})^2}\left(\begin{array}{cc}\displaystyle{
\frac{\partial}{\partial z}\log\left(\frac{H^f}{F^h}\right)} & \displaystyle{\frac{\partial}{\partial z}\log\left(\frac{H^g}{G^h}\right)}\\
& \\
\displaystyle{\frac{\partial}{\partial z'}\log\left(\frac{H^f}{F^h}\right)} & \displaystyle{\frac{\partial}{\partial z'}\log\left(\frac{H^g}{G^h}\right)}
\end{array}\right)}.\]
Le crochet $[F,G,H]$ est liŽ au dŽterminant de $J$~:
\begin{equation}[F,G,H]=\frac{g+h}{h}\det(J)F^{1+h}G^{1+h}H^{1-f-g}.\label{eq:formule}\end{equation} Donc $\mu=2$, et le corps ${\mathfrak K}$ est de
degrŽ de transcendance $4$. Ainsi, le corps \[{\mathfrak T}:=\CC\left(\frac{H^f}{F^h},\frac{H^g}{G^h}\right)\]
est de degrŽ de transcendance $2$ sur $\CC$. Comme une relation de dŽpendance algŽbri\-que pour $F,G,H$ est toujours dŽterminŽe
par un polyn™me isobare, ceci implique que $F,G,H$ sont algŽbriquement indŽpendantes.

Montrons ensuite que la condition (2) implique la condition (1). 
Si $F,G,H$ sont algŽbriquement indŽpendantes, alors 
le degrŽ de transcendance de ${\mathfrak T}$ est $2$.
On peut choisir deux paramtres complexes analytiquement indŽpendants $t_1,t_2$ tels que~:
\[{\mathfrak K}\cong\CC(t_1,t_2,e^{t_1},e^{t_2}),\] donc ${\mathfrak K}$ a degrŽ de transcendance $4$ sur $\CC$.

Supposons par l'absurde que $M=0$. La formule (\ref{eq:formule})
implique que les fonctions $\log(H^f/F^h)$ et $\log(H^g/G^h)$ sont analytiquement dŽpendantes. Donc les fonctions $H^f/F^h,H^g/G^h$
sont aussi analytiquement dŽpendantes. Mais deux fonctions modulaires sur ${\cal H}^2$ analytiquement dŽpendantes sont aussi algŽbri\-quement
dŽpendantes, ce qui donne au degrŽ de
transcendance
$\delta$ de
${\mathfrak K}$ la majoration
$\delta\leq 3$, d'o une contradiction. La proposition \ref{lemme:diff} est maintenant entirement dŽmontrŽe.

\medskip

\noindent {\bf Remarque.} Gr‰ce au lemme \ref{lemme:multi}, la proposition \ref{lemme:diff} se gŽnŽralise au cas de $n$ variables complexes avec $n\geq 3$,
mais nous n'en donnons pas les dŽtails ici.

\subsection{ThŽorme de structure pour $K=\QQ(\sqrt{5})$.\label{section:structure}}
On peut utiliser les opŽrateurs diffŽrentiels $\Lambda,\Pi,[\cdot,\cdot,\cdot]$ pour dŽterminer explicitement la structure de certains
anneaux de formes modulaires. Nous dŽcrivons explicitement toutes les formes modulaires de Hilbert 
dans le cas $K=\QQ(\sqrt{5})$.

\medskip

\begin{Lemme}
\`A multiplication par un nombre complexe non nul prs, il existe une unique forme modulaire de Hilbert pour $\Upsilon=\SL_2({\cal O}_K)$,
symŽtrique de poids $(15,15)$.
\label{lemme:chi15}\end{Lemme}
\noindent {\bf DŽmonstration.} Nous construisons la forme modulaire, puis nous dŽmon\-trons son unicitŽ.
En appliquant la proposition \ref{lemme:diff} ˆ $F=\chi_6,G=\varphi_2,H=\chi_5$, on voit que
\[\tilde{\chi}=[\chi_6,\varphi_2,\chi_5]\] est
symŽtrique de poids $(15,15)$ et non nulle. La non nullitŽ de $\tilde{\chi}$ peut
aussi tre vŽrifiŽe directement ˆ partir de la dŽfinition de $M$~: le coefficient de $\exp\{2\pi{\rm i}\tr(\und{z})\}$
dans la sŽrie de Fourier de $\tilde{\chi}$ est non nul, ce qui implique $\tilde{\chi}\not=0$. 
Ceci termine la preuve de l'existence.

Montrons maintenant l'unicitŽ d'une telle forme modulaire, ˆ une constante multiplicative prs.
Soit $\chi$ une forme modulaire de Hilbert, symŽtrique de poids $(15,15)$.  

Calculons $\chi|_{\Xi_2}$. 
Nous regardons les premiers coefficients de Fourier
de $\chi$ (ordonnŽs par trace croissante)~:
\[\chi(\und{z})=c_0+c_\nu\exp\{2\pi{\rm i}\tr(\nu\und{z})\}+c_{\nu'}\exp\{2\pi{\rm i}\tr(\nu'\und{z})\}+c_1\exp\{2\pi{\rm
i}\tr(\und{z})\}+\cdots\] avec $\nu=\epsilon/\sqrt{5}$.
La forme $\chi$ Žtant modulaire de poids impair (ici $(15,15)$), on a $\chi(U(\und{z}))=\no(\epsilon')^{15}\chi(\und{z})=-\chi(\und{z})$.
Ainsi $c_0=-c_0$ et $\chi$ est parabolique (\footnote{Une autre manire de remarquer ceci est 
d'Žcrire $M_{\und{15}}(\Upsilon)=\langle E_{15}\rangle\oplus S_{\und{15}}(\Upsilon)$, car $\Upsilon$ n'a qu'une
seule classe d'Žquivalence de pointe. De plus, puisque $K$ a une unitŽ de norme nŽgative, $E_{15}=0$.}).

Donc $\chi|_{\Xi_2}$
est une forme parabolique elliptique de poids $30$, donc
une combinaison linŽaire $x\Delta^2E_6+y\Delta E_6^3$ (noter que $p=15$ est le plus petit
entier impair tel que l'espace des formes paraboliques elliptiques de dimension $2p$ 
est de dimension $\geq 2$).

De mme, $c_\nu=-c_{\epsilon^2\nu}$, ce qui donne $c_\nu=-c_{\nu'}$. D'autre part, $\chi$ est symŽtrique, et donc $c_\nu=c_{\nu'}$.
Ainsi $c_\nu=c_{\nu'}=0$. Ceci implique $\chi|_{\Xi_2}=xq^2+\cdots=x\Delta^2E_6$~: naturellement, un
argument similaire s'applique aux formes modulaires symŽtriques de poids parallle impair (cf. plus bas)~:
on en dŽduit que si $h$ est une forme modulaire symŽtrique de poids impairs, alors $(h|_{\Xi_2})\Delta^{-2}E_6^{-1}$
est une forme modulaire elliptique.

Il existe une combinaison linŽaire $\phi=a\chi+b\tilde{\chi}$ ayant sa restriction ˆ $\Xi_2$ nulle. 
$\phi/\chi_5^2$ est une forme modulaire symŽtrique de poids $(5,5)$, nulle d'aprs le thŽorme \ref{lemme:structure1}~:
le lemme est dŽmontrŽ. En particulier, $\tilde{\chi}$ ne s'annule pas dans $\Xi_2$.
Nous dŽmontrons maintenant~:

\begin{Theoreme} L'anneau des formes modulaires de Hilbert pour le groupe $\Upsilon=\SL_2({\cal O}_K)$, avec
$K=\QQ(\sqrt{5})$ est Žgal ˆ l'anneau (graduŽ par les poids) des polyn™mes en $\varphi_2,\chi_5,\chi_6,\tilde{\chi}$,
quotientŽ par la relation~:
\begin{eqnarray}
\frac{5}{49}\tilde{\chi}^2 & = & 50000\chi_5^6-1000\varphi_2^2\chi_6
\chi_5^4+\varphi_2^5\chi_5^4-2\varphi_2^4\chi_6^2\chi_5^2+\nonumber\\
& &+1800\varphi_2\chi_6^3\chi_5^2+\varphi_2^3\chi_6^4-864\chi_6^5.\label{eq:klein}
\end{eqnarray}\label{lemme:structure2}\end{Theoreme}
\noindent {\bf DŽmonstration.} Comme $\tilde{\chi}^2$ est une forme modulaire symŽtrique de poids parallle pair $30$,
c'est un polyn™me en $\varphi_2,\chi_5^2,\chi_6$ d'aprs le thŽorme
\ref{lemme:structure1}. La relation (\ref{eq:klein}) est uniquement
dŽterminŽe. 

Nous remarquons qu'il suffit de dŽcrire le sous-anneau des formes modulaires symŽtriques. En effet, soit
$\phi$ une forme modulaire de Hilbert de poids parallle. On peut Žcrire~:
\begin{eqnarray*}
\phi(z,z') & = & \frac{1}{2}(f(z,z')+f(z',z))+\frac{1}{2}(f(z,z')-f(z',z))\\
& = & f_1+\chi_5f_2,
\end{eqnarray*}
avec $f_1,f_2$ deux formes modulaires symŽtriques.

Soit maintenant $f$ une forme modulaire symŽtrique de poids $(r,r)$, avec $r$ impair. On voit que $r\geq 15$.
La forme modulaire
elliptique $g=f|_{\Xi_2}$ est de poids $2r$, et s'annule ˆ l'infini avec une multiplicitŽ $\geq 2$ (reprendre l'idŽe
de la dŽmonstration du lemme \ref{lemme:chi15}). Ainsi~:
\[g=E_6\Delta^2P(E_4,\Delta),\] pour un certain polyn™me isobare $P$.
Il existe un polyn™me isobare $Q$ tel que la forme modulaire de Hilbert symŽtrique~:
\[h=f-\tilde{\chi}Q(\varphi_2,\chi_6)\]
s'annule sur $\Xi_2$. Or, $h/\chi_5^2$ est une forme modulaire de poids $(r-10,r-10)$~: on applique une hypothse
de rŽcurrence. On trouve que $f=\tilde{\chi}k$, avec $k$ forme modulaire symŽtrique de poids pair $(r-15,r-15)$.
Le thŽorme  est dŽmontrŽ, mais reste la question du calcul de la relation (\ref{eq:klein}).

\medskip

Sans passer par la gŽomŽtrie des surfaces (ce qui oblige ˆ calculer explicitement une dŽsingularisation
de la surface $X_\Upsilon$, comme le fait Hirzebruch dans
\cite{Hirzebruch:Klein}), la technique la plus avantageuse est encore une fois de passer par
des relations diffŽrentielles, et de rŽsoudre explicitement des systmes linŽaires avec peu d'Žquations
et  d'inconnues. 

En calculant explicitement les coefficients de Fourier de 
\[[\tilde{\chi}^2,\varphi_2,\chi_5],[\tilde{\chi}^2,\varphi_2,\chi_6],[\tilde{\chi}^2,\chi_5,\chi_6],\]
on explicite la relation (\ref{eq:klein}).

\medskip

\noindent {\bf Remarque.} On montre que \[\tilde{\chi}=\frac{28}{\sqrt{5}}\chi_{15},\] suivant les notations de Resnikoff.
On peut comparer notre construction ˆ celle de Gundlach dans \cite{Gundlach:Bestimmung1} pp. 241-247. Gundlach doit d'abord 
construire des sŽries d'Eisenstein tordues par des caractres, de poids $(1,1)$ pour des sous-groupes de congruence, ce qui
est difficile en gŽnŽral.
Notre construction est considŽrablement plus simple, et plus explicite.

\subsection{Polyn™mes diffŽrentiels\label{section:differentiels}.}
Soit $\und{p}\in\NN^n$. Nous posons~:
\[D_{\und{p}}:=\frac{\partial^{|\und{p}|}}{(\partial z_1)^{p_1}\cdots(\partial z_n)^{p_n}},\] o $|\und{p}|:=p_1+\cdots+p_n$.

\medskip

\noindent {\bf DŽfinition.} Un {\em polyn™me diffŽrentiel} $D$ d'ordre $\leq \sigma$, agissant sur 
l'espace des fonctions holomorphes sur ${\cal H}^n$, est par dŽfinition un opŽrateur de la forme~:
\begin{equation}
DX=\sum_{({\cal I},h)}a_{{\cal I},h}\prod_{\und{p}\in{\cal I}}(D_{\und{p}}(X))^{h(\und{p})},\label{eq:ecriture}
\end{equation}
o toutes les sommes et produits sont finis, ${\cal I}$ est un sous-ensemble fini de $\NN^n$, $h$ une fonction ${\cal
I}\rightarrow\NN$, la somme sur les $({\cal I},h)$ porte sur des ensembles ${\cal I}$ tels
que pour tout
$\und{p}\in{\cal I}$ on a $|\und{p}|\leq \sigma$, les $a_{{\cal I},h}$ sont des nombres complexes.

\medskip

Soit $\Gamma$ un groupe modulaire de Hilbert.
Nous nous intŽressons ici aux polyn™mes diffŽrentiels $D$ ayant la propriŽtŽ que $D(M_{\und{r}}(\Gamma))\subset M_{\und{s}}(\Gamma)$,
pour $\und{r},\und{s}\in\NN^n$. Nous commenons par un lemme ŽlŽmentaire.

\medskip

\begin{Lemme} Soit $F\in M_{\und{f}}(\Gamma)$, soit $\und{p}$ un ŽlŽment de $\NN^n$
tel que $|\und{p}|>1$. 
Alors il existe un entier $s\in\ZZ$ tel que~:
\[D_{\und{p}}F=F^s({\cal D}F+\Phi F),\] o ${\cal D}$ est un polyn™me diffŽrentiel qui est une composition de multiples
d'opŽ\-rateurs
diffŽrentiels $\Pi_i,\Lambda_{i,k},(X,Y)\mapsto [X,Y]_{1_i}$ avec $i,k\in\{1,\ldots,n\}$ agissant sur $F$, 
et $\Phi F$ est l'image de $F$ par un polyn™me diffŽrentiel d'ordre $<|\und{p}|$.
\label{lemme:differentiel_basique}\end{Lemme}

\medskip

\noindent {\bf DŽmonstration.} 
Nous commenons par Žcrire~:
\begin{eqnarray*}
\Pi_iF & = & \frac{1}{(2\pi{\rm i})^2}f_iF\frac{\partial^2F}{\partial z_i^2}+\cdots\\
\Lambda_{i,k}F & = & \frac{1}{(2\pi{\rm i})^2}F\frac{\partial^2F}{\partial z_i\partial z_j}+\cdots\\
{[}F,{[}F,\ldots {[}F,\Lambda_{i,k}F{]}_{1_{k_s}}\ldots{]}_{1_{k_2}}{]}_{1_{k_1}} & = & 
\frac{1}{(2\pi{\rm i})^{s+2}}F^s\frac{\partial^{s+2} F}{\partial z_{k_1}\cdots
\partial z_{k_s}\partial z_i\partial z_j}+\cdots\\
{[}F,{[}F,\ldots {[}F,\Pi_iF{]}_{1_{i}}\ldots{]}_{1_{i}}{]}_{1_{i}} & = & 
\frac{1}{(2\pi{\rm i})^{s+2}}f_iF^s\frac{\partial^{s+2} F}{\partial z_i^{s+2}}+\cdots
\end{eqnarray*}
o les symboles $+\cdots$ dŽsignent la prŽsence de termes qui sont des polyn™mes diffŽrentiels d'ordre infŽrieur, ŽvaluŽs en $F$.
Le lemme \ref{lemme:differentiel_basique} est dŽmontrŽ en explicitant l'expression de $D_{\und{p}}F$
dans les formules ci-dessus.

\medskip

\noindent {\bf Remarque.} Si $n=1$ et $F\in M_k(\SL_2(\ZZ))$, alors en posant~:
\[G_1F:=[F,[F,F]_2]_1\mbox{ et }G_s(F):=[F,G_{s-1}F]_1,\] on a que $G_s(F)$ est une forme modulaire telle que
\[G_sF=cF^{s+2}\frac{d^{s+2}F}{dz^{s+2}}+\cdots,\] o $c$ est une constante dŽpendant de $k$ et $s$. 

\medskip

Rankin a montrŽ que tout polyn™me diffŽrentiel $D$ tel que \[D(M_r(\SL_2(\ZZ)))\subset M_s(\SL_2(\ZZ)),\]
a la propriŽtŽ que $DF$ appartient ˆ $\CC(F)[[F,F]_4,G_1f,G_2f,\ldots]$. Nous gŽnŽra\-lisons son rŽsultat,
suivant Resnikoff \cite{Resnikoff:Automorphic}.

\medskip

\begin{Lemme} Soit 
$F\in M_{\und{f}}(\Gamma)$ une forme modulaire non nulle. Soit
$D$ un polyn™me diffŽrentiel non nul, ayant la propriŽtŽ que $DF\in M_{\und{t}}(\Gamma)$,
pour $\und{t}\in\NN^n$. Alors $DF$ est un polyn™me ˆ coefficients dans $\CC(F)$,
en des compositions
de multiples d'opŽ\-rateurs diffŽrentiels $X\mapsto[X,X]_{4_i},X\mapsto\langle X\rangle_{2_i+2_j},(X,Y)\mapsto[X,Y]_{2_i}$ agissant sur $F$.
\label{lemme:operateurs}\end{Lemme}
 
\medskip

\noindent {\bf DŽmonstration.} On a l'Žcriture (\ref{eq:ecriture})~: on peut appliquer le lemme
\ref{lemme:differentiel_basique}. Notons ${\cal D}_{\und{\sigma}}$ un opŽrateur diffŽrentiel ${\cal D}$ obtenu par
application de ce lemme, associŽ au mon™me diffŽrentiel $D_{\und{\sigma}}$.  On a~:
\begin{eqnarray*}
DF & = & \sum_{({\cal I},h)}a_{{\cal I},h}\prod_{\und{p}\in{\cal I}}(D_{\und{p}}F)^{h(\und{p})}\\
& = & \sum_{({\cal I},h)}a_{{\cal I},h}\prod_{\und{p}\in{\cal I}}F^{s(\und{p})}({\cal D}_{\und{p}}F+\Phi_{\und{p}}F)^{h(\und{p})}\\ 
& = & \sum_{s\in\ZZ}\sum_{({\cal J},k)}b_{{\cal
J},k,s}F^s{\prod_{\und{p}\in{\cal J}}}^\sharp({\cal D}_{\und{p}}F)^{k(\und{p})}
{\prod}^\flat_{\und{p}\in{\cal J}}(D_{\und{p}}F)^{k(\und{p})},
\end{eqnarray*}
o toutes les sommes et produits sont finies, une somme est indexŽe sur des couples $({\cal J},k:{\cal J}\rightarrow\NN)$, les symboles $s$ dŽsignent
des fonctions ${\cal I}\rightarrow\ZZ$, 
et le produit $\prod^\sharp$ est indexŽ
par les
$\und{p}\in{\cal J}$ tels que
$|\und{p}|>1$,  le produit $\prod^\flat$ est indexŽ par les $\und{p}\in{\cal J}$ tels que $|\und{p}|=1$.
Noter que si $|\und{p}|=1$, alors il existe $j$ tel que~:
\[D_{\und{p}}F=\frac{\partial F}{\partial z_j}.\]
Les termes mon™miaux donnŽs par les produits $F^s\prod_{{\cal J}}^\sharp$ 
sont tous des formes modulaires, disons de poids $\und{d}({\cal J},s)$.
On a, pour $\gamma=\sqm{a}{b}{c}{d}\in\Gamma$~:
\[\left(\frac{\partial F}{\partial z_j}\right)(\gamma(\und{z}))=\left(
\prod_{i=1}^n(c_iz_i+d_i)^{f_i}\right)(c_jz_j+d_j)^2\left(\frac{f_jc_j}{c_jz_j+d_j}F(\und{z})+
\left(\frac{\partial F}{\partial z_j}\right)(\und{z})\right).\]
Nous pouvons supposer qu'au moins un parmi les coefficients $b_{{\cal J},k,s}$ soit non nul.
D'aprs cette identitŽ, on voit que pour ${\cal J},k,s$ tels que le coefficient $b_{{\cal J},k,s}$ est non nul~:
\[\und{t}=\und{d}({\cal J},k,s)+\sum^\flat_{\und{p}\in{\cal I}}(k(\und{p})\und{f}+2_j).\]
On en dŽduit que
pour $\und{z}\in{\cal H}^n$ fixŽ, les fonctions $\Gamma\rightarrow\CC$~:
\[\prod^\flat_{\und{p}\in{\cal
J}}\left(\frac{f_{j}c_{j}}{c_{j}z_{j}+d_{j}}
F(\und{z})+D_{\und{p}}F(\und{z})\right)^{k(\und{p})}\] et $1$ sont
$\CC$-linŽairement dŽpendantes.
En d'autres termes, l'image de la fonction $R:\Gamma\rightarrow\CC^n$ dŽfinie par~:
\[\gamma=\sqm{a}{b}{c}{d}\mapsto\left(\frac{f_1c_1}{c_1z_1+d_1},\ldots,\frac{f_nc_n}{c_nz_n+d_n}\right)\]
est contenue dans une sous-variŽtŽ algŽbrique propre de $\CC^n$.

Nous montrons maintenant que cela n'est pas possible~: ceci nous don\-ne\-ra que $k(\und{p})=0$ pour tout $\und{p}$, et la
fin de la dŽmonstration du lemme.

L'image de $R$ est fermŽe (au sens de Zariski) dans $\CC^n$ si et seulement si
\[\displaystyle{\left\{\left(\frac{1}{z_1+\frac{c_1}{d_1}},\ldots,\frac{1}{z_n+\frac{c_n}{d_n}}\right)\right\}_{\gamma\in\Gamma}}\]
l'est (nous avons vu que $\prod_if_i\not=0$~: lemme \ref{lemme:M0r=CC}), c'est-ˆ-dire si et seulement si
\[\left\{\left(\frac{c_1}{d_1},\ldots,\frac{c_n}{d_n}\right)\right\}_{\gamma\in\Gamma}\] est fermŽe.
Or, on voit sans difficultŽ que cet ensemble est dense dans $\RR^n$ pour la topologie euclidienne.

\begin{Corollaire} Supposons que $n=2$. Soit $D$ un polyn™me diffŽrentiel non nul d'ordre $\leq 2$, 
et de degrŽ minimal avec la propriŽtŽ que pour une forme modulaire de Hilbert non constante de poids parallle
$F$ on ait $DF=0$.
alors il existe un opŽrateur diffŽrentiel non nul $EX\in\CC[X,\Lambda X,\Pi X]$, tel que $EF=0$.
\label{lemme:symetrique}\end{Corollaire}

\medskip

\noindent {\bf DŽmonstration.} D'aprs le lemme \ref{lemme:operateurs}, $DF\in\CC(F)[\Lambda F,\Pi_1F,\Pi_2F]$, car
l'ordre de $D$ Žtant $\leq 2$, il ne peut pas y avoir de terme faisant intervenir les crochets $(X,Y)\mapsto[X,Y]_{1_i}$ dans
toute expression de $D$. \'Ecrivons donc~:
\[DF=\sum_{(x,y,z,t)\in\ZZ\times\NN^3}c_{x,y,z,t}F^x(\Lambda F)^y(\Pi_1F)^z(\Pi_2F)^t,\]
et supposons $F$ de poids parallle $f$. Le poids de $F^x(\Lambda F)^y(\Pi_1F)^z(\Pi_2F)^t$ est 
$(fx+(2f+2)y+2f(z+t))\und{1}+4_1z+4_2t$. D'aprs le lemme \ref{lemme:annulation_operateurs}, tous les mon™mes 
$F^x(\Lambda F)^y(\Pi_1F)^z(\Pi_2F)^t$ sont non nuls. Il est alors clair que s'il existe $(x,y,z,t)\in\ZZ\times\NN^3$ avec 
$c_{x,y,z,t}\not=0$ et $z\not=t$, alors une puissance non nulle et positive de $\Pi_1F$ (ou de $\Pi_2F$)
divise $DF$ de telle sorte qu'on puisse Žcrire~:
\[DF=(\Pi_iF)^s\sum_{(x,y,z)\in\ZZ\times\NN^2}c_{x,y,z}'F^x(\Lambda F)^y(\Pi F)^t=(\Pi_iF)^sF^wEF,\] avec
$EF\in\CC[F,\Lambda F,\Pi F]$ non nul. Le corollaire est dŽmontrŽ.


\subsection{\'Equations diffŽrentielles ayant des solutions mo\-dulaires.}

Dans ce sous-paragraphe, nous supposons $[K:\QQ]=2$.
Nous commenons par un lemme technique ŽlŽmentaire, qui sera utilisŽ dans la suite.

\begin{Lemme} 
Les opŽrateurs diffŽrentiels \begin{small}\[X\mapsto c\Pi X -(\Lambda X)^2,\quad X\mapsto [X,\Pi X]_{1_i},\quad X\mapsto [X,\Lambda
X]_{1_i},\quad X\mapsto [X,\Pi_iX]_{1_i}\]\end{small}  avec $c\in\CC$, n'annulent aucune forme modulaire non constante.
\label{lemme:injectivite}\end{Lemme}

\noindent {\bf Esquisse de dŽmonstration.}
Nous faisons seulement une partie des dŽ\-monstrations.
 Nous commenons par donner les formules explicites dŽcri\-vant l'action de ces opŽrateurs sur les sŽries de Fourier.
Soient $F,G$ deux formes modulaires de poids $\und{f}$ et $\und{g}$, de sŽries de Fourier~:
\[F(\und{z})=\sum_{\nu\in{\cal O}_{K,+}^*\cup\{0\}}a_\nu \exp\{2\pi{\rm i}\tr(\nu\und{z})\},\quad
G(\und{z})=\sum_{\nu\in{\cal O}_{K,+}^*\cup\{0\}}b_\nu \exp\{2\pi{\rm i}\tr(\nu\und{z})\}.\]
On a~:
\begin{eqnarray}
\Pi_iF & = & \sum_{\tau\in{\cal O}_{K,+}^*}\exp\{2\pi{\rm i}\tr(\tau\und{z})\}{\sum_{\nu+\mu=\tau}}^\sharp a_\nu
a_\mu(f_i\sigma_i(\mu)^2-(f_i+1)\sigma_i(\nu\mu)),\nonumber\\ 
\Lambda F & = & \sum_{\tau\in{\cal O}_{K,+}^*}\exp\{2\pi{\rm
i}\tr(\tau\und{z})\}{\sum_{\nu+\mu=\tau}}^\sharp a_\nu a_\mu(\no(\mu)-\nu\mu'),\nonumber\\ 
\Pi f & = & \sum_{\tau\in{\cal O}_{K,+}^*}\exp\{2\pi{\rm
i}\tr(\tau\und{z})\}\times\nonumber\\ & & \times{\sum_{\alpha+\beta+\gamma+\delta=\tau}}^\sharp 
a_\alpha a_\beta a_\gamma a_\delta(f_1\alpha^2-(f_1+1)\beta\alpha)(f_2\gamma^2-(f_2+1)\delta\gamma),\nonumber\\
{[}F,G{]}_{1_i} & = & \sum_{\tau\in{\cal O}_K^*{}_+}\exp\{2\pi{\rm i}\tr(\tau\und{z})\}{\sum_{\nu+\mu=\tau}}^\sharp
(g_i\sigma_i(\nu)-f_i\sigma_i(\mu))a_{\nu}b_{\mu},
\label{eq:efg}
\end{eqnarray}
o les sommes $\sum^\sharp$ sont indexŽes par des ŽlŽments de ${\cal O}_{K,+}^*\cup\{0\}$.
Certaines parties du lemme se dŽmontrent bien en utilisant ces formules. Pour d'autres parties, c'est mieux d'appliquer 
l'hypothse de modularitŽ.

\medskip

\noindent (1). Pour simplifier, nous supposons que $F$ est une forme parabolique. Si $F$ est non nulle,
alors il existe un plus petit entier $m>0$ tel que pour quelques $\nu\in{\cal O}_{K,+}^*$ avec $\tr(\nu)=m$ on ait
$a_\nu\not=0$. On peut supposer que le nombre rŽel positif $\sigma_1(\nu)$ soit le plus grand possible avec cette
propriŽtŽ. On trouve alors que le coefficient de $\exp\{2\pi{\rm i}\tr(4\nu\und{z})\}$ dans la sŽrie de Fourier de $\Pi F$
est Žgal ˆ $\no(\nu)^2a_\nu^4\not=0$, et pour tout $\mu\in{\cal O}_{K,+}^*$ de trace $<\tr(4\nu)$, 
le coefficient de $\exp\{2\pi{\rm i}\tr(\mu\und{z})\}$ dans la sŽrie de Fourier de $\Pi F$ est nul.

Mais le choix de $\nu$ implique que le coefficient de $\exp\{2\pi{\rm i}\tr(2\nu\und{z})\}$ dans la sŽrie de Fourier de $\Lambda F$,
est nul (car Žgal ˆ $a_\nu^2(\no(\nu)-\nu\nu')$), et pour tout $\mu\in{\cal O}_{K,+}^*$ de trace $<\tr(2\nu)$, 
le coefficient de $\exp\{2\pi{\rm i}\tr(\mu\und{z})\}$ dans la sŽrie de Fourier de $\Lambda F$ est nul.
Si $c\not=0$ on en dŽduit que $a_\nu=0$, une contradiction. Donc $c=0$, et nous sommes ramenŽs au cas (1) du lemme
\ref{lemme:annulation_operateurs}.

\medskip

\noindent (2). Commenons par considŽrer $F,G$ non nulles, telles que $[F,G]_{1_i}=0$.
Supposons pour simplifier qu'elles soient toutes deux paraboliques. 
Il existe
$\nu,\mu\in{\cal O}_K^*{}_+$, de trace minimale, $\sigma_1(\nu),\sigma_1(\mu)$ plus grands possibles, 
avec la propriŽtŽ que $a_\nu$ et $b_\mu$ soient non nuls.
Pour $\tau=\nu+\mu$, on trouve que le coefficient de Fourier de $[F,G]_{1_i}$ associŽ ˆ $\tau$ est Žgal ˆ 
$(g_i\sigma_i(\nu)-f_i\sigma_i(\mu))a_{\nu}b_{\mu}$. Il faut donc que $g_i\sigma_i(\nu)-f_i\sigma_i(\mu)=0$.

Posons $G=\Pi F$~: nous savons que $G\not=0$. Alors $b_{4\nu}=a_{\nu}^4\nu^4$ est non nul, et 
$\mu=4\nu$ est de trace minimale avec cette propriŽtŽ.
Donc $g_i=4f_i$, ce qui est impossible, car $G$ est une forme modulaire de poids $\und{g}=4\und{f}+\und{4}$.

\medskip

\noindent (3). Pour montrer que $[F,\Lambda F]_{1_i}$ ne peut pas s'annuler sur 
$M_{\und{f}}(\Gamma)$,
il suffit de montrer que $\displaystyle{\frac{\partial}{\partial z_i}
\frac{(\Lambda F)^{f_i}}{F^{2 f_i+2}}}\not=0$~: on utilise la modularitŽ de $F$
pour montrer que c'est le cas.

\medskip

\noindent (4). Les formes modulaires $F$ et $G=\Pi_iF$ sont non nulles et leur poids $\und{f}$ et $\und{g}$ satisfont
$\und{f}\not\in\QQ^\times\und{g}$. D'aprs le lemme \ref{lemme:annulation_operateurs}, (3), la forme modulaire $[F,G]_{1_i}$
ne peut pas tre nulle. Nous laissons au lecteur le soin de complter les autres parties de la dŽmonstration du lemme
\ref{lemme:injectivite}.

\begin{Proposition}
Soit $F$ une forme modulaire de Hilbert symŽtrique de poids $r$, annulant un poly\-n™me diffŽrentiel non nul d'ordre $\leq 2$.
Alors la forme modulaire $TF$~:
\begin{equation}
TF=[F,\Lambda F]_{1_1}[F,\Pi F]_{1_2}-[F,\Lambda F]_{1_2}[F,\Pi F]_{1_1}\label{eq:finipio}\end{equation} est nulle.
\label{lemme:ind_alg_diff}\end{Proposition}

\noindent {\bf DŽmonstration.} Soit $F$ une forme modulaire symŽtrique~:
son poids $\und{r}$ est parallle. Soit $D$ un opŽrateur diffŽrentiel d'ordre
$\leq 2$ s'annulant en $F$. 
Le corollaire \ref{lemme:symetrique} implique qu'il existe un opŽrateur diffŽrentiel $E$ tel que
$EX\in\CC[X,\Lambda X,\Pi X]$, et
\begin{equation}
EF=\sum_{\und{i}}c_{\und{i}}F^{i_1}(\Lambda F)^{i_2}(\Pi F)^{i_3}=0,\label{eq:principio}
\end{equation} o la somme est finie, et porte sur des triplets $\und{i}$ d'entiers
positifs ou nuls.

Naturellement, le polyn™me $EF$ est isobare, et on a $i_1r+(2i_2+4i_3)(r+1)=s$ pour tout $\und{i}$.
On voit aussi que $EX$ doit avoir tous ses degrŽs partiels non nuls. En effet, si $EX$ n'a qu'un seul degrŽ partiel non nul,
on rencontre une contradiction avec le lemme \ref{lemme:injectivite}.

Si $EX$ a seulement deux degrŽs partiels non nuls, on procde de la mme faon.
Si par exemple $EX$ ne dŽpend pas de $X$, alors~:
\[\sum_{i}c_{i}\left(\frac{\Pi F}{(\Lambda F)^2}\right)^i=0\] et on est ramenŽ ˆ rŽsoudre l'Žquation diffŽrentielle
$\Pi F=c(\Lambda F)^2$, qui n'a pas de solutions modulaires autres que la solution nulle, d'aprs le lemme \ref{lemme:injectivite}.

Donc $EX$ a tous ses degrŽs partiels non nuls. On peut reformuler (\ref{eq:principio}) ainsi~: les fonctions
\[\alpha=\frac{\Lambda F}{F^{(2r+2)/r}},\beta=\frac{\Pi F}{F^{(4r+4)/r}}\] sont algŽbriquement dŽpendantes sur $\CC$.

Ceci implique que $\alpha,\beta$ sont analytiquement dŽpendantes~: $\alpha\wedge\beta=0$. 
Mais cette condition n'est qu'une rŽŽcriture de (\ref{eq:finipio})~:
on voit que $\alpha\wedge\beta=(rF^{(6r+6)/r})^{-1}H$ et $TF=rF^2H$, pour une 
forme modulaire $H$, d'o $TF=0$ si et seulement si $H=0$ si et seulement si $\alpha\wedge\beta=0$.
La proposition
\ref{lemme:ind_alg_diff} est dŽmontrŽe.

On remarque que $TF$ est antisymŽtrique et divisible par $F^2$. Pour voir ceci on utilise
les techniques de preuve de la premire partie de la proposition \ref{lemme:diff}.
Pour $K=\QQ(\sqrt{5})$, le thŽorme \ref{lemme:structure2}
implique que $(T\varphi_2)\varphi_2^{-2}$ est
un multiple non nul de $\chi_5\tilde{\chi}$~: on trouve $(T\varphi_2)\varphi_2^{-2}=2^{11}\cdot3^3\cdot 5\cdot
7^{-1}\chi_5\tilde{\chi}$.

\medskip

\begin{Corollaire}
Si $f$ est une forme modulaire symŽtri\-que qui n'annule pas l'opŽrateur diffŽrentiel $T$ dŽfini par 
(\ref{eq:finipio}), alors $F,\Lambda F,\Pi_iF$ et $\Pi_2F$ sont 
algŽ\-briquement indŽpen\-dantes sur $\CC$, et donc
$F,\Lambda F$ et $\Pi F$ sont aussi algŽbri\-quement indŽpendantes.\label{lemme:indep_alg}\end{Corollaire}

Dans le cas o $K=\QQ(\sqrt{5})$, un calcul numŽrique explicite sur la sŽrie de Fourier de $\varphi_2$ permet de vŽrifier que $\varphi_2$ 
ne satisfait pas (\ref{eq:finipio}). Donc $\varphi_2$ ne satisfait aucune Žquation diffŽrentielle d'ordre $\leq 2$, et
$\varphi_2,\chi_6,\chi_5$ sont algŽbriquement indŽpendantes. Dans le cas gŽnŽral, nous pouvons utiliser
la proposition suivante.

\begin{Proposition} Soit $F$ une forme modulaire symŽtrique, soit
\[F(\und{z})=a_0+\sum_{\nu\in{\cal O}^*_{K,+}}a_\nu \exp\{2\pi{\rm i}\tr(\nu\und{z})\}\] son dŽveloppement en sŽrie de
Fourier ˆ l'infini. Supposons qu'il existe au moins deux ŽlŽments distincts $\nu,\mu\in{\cal O}^*_{K,+}$, de trace
minimale avec la propriŽtŽ que $a_\nu,a_\mu\not=0$. Alors $TF$, dŽfinie par (\ref{eq:finipio}), est non nulle.
\label{lemme:resnikoff}\end{Proposition}

\noindent {\bf DŽmonstration.} 
Supposons pour simplifier que $F$ n'est pas une forme para\-bolique, et n'est pas constante.
On calcule explicitement le dŽveloppement en sŽrie de Fourier de $TF$.
Posons $G=\Lambda F,H=\Pi F$~: ce sont des formes modulaires de poids parallle 
$g$ et $h$. \'Ecrivons les sŽries de
Fourier de
$G,H$ comme suit~:
\[G(\und{z})=\sum_{\nu\in{\cal O}^*_{K,+}}b_\nu \exp\{2\pi{\rm i}\tr(\nu\und{z})\},\quad H(\und{z})=
\sum_{\nu\in{\cal O}^*_{K,+}}c_\nu \exp\{2\pi{\rm i}\tr(\nu\und{z})\}.\]
Si on Žcrit~:
\[(TF)(\und{z})=\sum_{\nu\in{\cal O}^*_{K,+}}s_\nu \exp\{2\pi{\rm i}\tr(\nu\und{z})\},\]
en utilisant (\ref{eq:efg}) on trouve que
\[s_\nu=\sum_{\alpha+\beta+\gamma+\delta=\nu}t_{\alpha,\beta,\gamma,\delta}a_\alpha a_\gamma b_\beta c_\delta,\]
o
\begin{eqnarray*}
t_{\alpha,\beta,\gamma,\delta} & = &g\det\sqm{\sigma_1(\delta)}{\sigma_2(\delta)}{\sigma_1(\alpha)}{\sigma_2(\gamma)}+
h\det\sqm{\sigma_1(\alpha)}{\sigma_2(\gamma)}{\sigma_1(\beta)}{\sigma_2(\beta)}+\\ & &+
f\det\sqm{\sigma_1(\beta)}{\sigma_2(\beta)}{\sigma_1(\delta)}{\sigma_2(\delta)}.\end{eqnarray*}

Soient $\mu,\nu$ comme dans les hypothses. Soit $l=\tr(\mu)=\tr(\nu)>0$.
Les hypothses de la proposition impliquent qu'on peut choisir $l$ minimal, $\mu,\nu$ 
de telle sorte que $\Sigma(\mu),\Sigma(\nu)$ soient 
dans le segment \[\{\Sigma(\tau)\mbox{ avec }\tr(\tau)=l\},\]
la distance entre $\Sigma(\mu)$ et $\Sigma(\nu)$ soit la plus petite possible, et $\sigma_1(\mu)$ soit le plus grand possible~:
on voit alors que $\mu,\nu$ sont $\QQ$-linŽairement indŽpendants.

Cette condition de minimalitŽ impose une Žcriture plus simple pour certains coefficients $s_\nu$. 
En effet, si $\beta,\delta\not=0$, on vŽrifie que $b_\beta\not=0$ implique $\tr(\beta)\geq l$ et $c_\delta\not=0$
implique $\tr(\delta)\geq 2l$ (se souvenir que $a_0\not=0$).

Si $\tau=2\mu+\nu$
on a $\tau=\alpha+\beta+\gamma+\delta$ et $a_\alpha a_\gamma b_\beta c_\delta\not=0$ si et  seulement si $\alpha=\gamma=0,\beta=\nu,\delta=2\mu$, ou
$\alpha=\gamma=0,\beta=\mu,\delta=\nu+\mu$, et~:
\begin{eqnarray*}
s_\tau & = & r_f^2a_0^2(\sigma_1(\nu)\sigma_2(\mu)-\sigma_1(\mu)\sigma_2(\nu))(2b_{\nu}c_{2\mu}-b_{\mu}c_{\mu+\nu}).
\end{eqnarray*}
Les coefficients $b_{\nu},c_{2\mu},b_{\mu},c_{\mu+\nu}$ peuvent se calculer explicitement en fonction de $a_\mu,a_\nu$ gr‰ce ˆ la 
condition de minimalitŽ de $l$. On trouve 
\[s_\nu=r_f^4\no(\mu)a_0^5(\sigma_1(\nu)\sigma_2(\mu)-\sigma_1(\mu)\sigma_2(\nu))^3a_\mu^2a_\nu.\]
Cette quantitŽ ne peut pas s'annuler si $r_f\not=0$, car $\nu,\mu$ sont $\QQ$-linŽairement indŽpendants.
En particulier si non nulles, les sŽries $E_{r,{\mathfrak A}}$ satisfont les conditions de la proposition.

La dŽmonstration dans le cas o $f$ est une forme parabolique est similaire, nous la laissons au lecteur.

\medskip

\noindent {\bf Question}. {\em Y a-t-il une forme modulaire non nulle annulant (\ref{eq:finipio})?}

\medskip

On voit que toute sŽrie d'Eisenstein non nulle satisfait les hypothses de la proposition \ref{lemme:resnikoff}.
On en dŽduit~:

\begin{Theoreme}
L'anneau $\CC[E_2,\Lambda E_2,\Pi E_2]$
a un degrŽ de transcendance $3$ sur $\CC$. 
\label{lemme:structure}\end{Theoreme}

Comme le degrŽ de transcendance du corps $F_{\und{0}}(\Gamma)$ des fonctions mo\-dulaires de poids $\und{0}$ est $2$, on trouve 
que la cl™ture algŽbrique du corps de fractions de l'anneau \[{\cal T}(\Gamma)=\bigoplus_rM_{r}(\Gamma)\] des formes modulaires de poids parallle 
est de degrŽ de transcendance $3$. Ainsi, quatre formes modulaires de Hilbert de poids parallles sont toujours algŽbriquement
dŽpendantes.

\begin{Corollaire}
Pour tout corps quadratique rŽel $K$, il existe trois formes modulaires de Hilbert pour $\Gamma_K$, symŽtriques et algŽbriquement indŽpendantes, 
de poids parallle $2,6,10$, et une forme modulaire antisymŽtrique de poids $20$. 
\label{lemme:hammond}\end{Corollaire}

\noindent {\bf DŽmonstration.} Si $F$ est une forme modulaire de Hilbert de poids $\und{r}=(r_1,r_2)$,
\[\Phi F:=(\Pi F-(r_1+1)(r_2+1)(\Lambda F)^2)F^{-1}\] est une forme modulaire de Hilbert de poids $3\und{r}+\und{4}$ qui est non nulle
d'aprs le lemme \ref{lemme:injectivite} (dans l'expression explicite de 
$\Pi F-(r_1+1)(r_2+1)(\Lambda F)^2$, tous les termes sont des produits de $F$ et de fonctions holomorphes sur ${\cal H}^2$). 
En particulier, $\Phi E_2$ est une
forme modulaire non nulle de poids $10$.

D'aprs le corollaire \ref{lemme:indep_alg}, les formes modulaires $F=E_2,G=\Lambda E_2$ et $\Pi E_2$ sont algŽbriquement
indŽpendantes, de poids parallles $2,6,12$, et sont aussi clairement symŽtriques. Donc la forme modulaire $H=\Phi E_2$
est symŽtrique de poids $10$ et $F,G,H$ sont algŽbriquement indŽpendantes.

La proposition \ref{lemme:diff} implique
l'existence d'une forme modulaire non nulle $M$, de poids $20$. Cette forme modulaire est antisymŽtrique,
car $F,G,H$ sont symŽtriques. 

\medskip

\noindent {\bf Remarque.}
Soit ${\cal T}_0={\cal T}_0(\Gamma)$ l'anneau engendrŽ par les formes modulaires de Hilbert symŽtriques de poids pair
associŽes ˆ $\Gamma$,
soit  ${\mathfrak P}$ l'idŽal de ${\cal T}_0$ engendrŽ par toutes les formes modulaires de Hilbert $F$ telles que
$F|_{\Xi_2}=0$.
Si le discriminant de $K$ est une somme de deux carrŽs, et si ${\cal P}$ est principal,
alors ${\cal T}_0$ est engendrŽ par $F,G$ et un gŽnŽrateur de ${\cal P}$ (cf. \cite{Gundlach:Bestimmung1}).
Ici, nous ne pouvons pas dŽmontrer que ${\cal T}_0=\CC[F,G,H]$, mais seulement que
$H\in{\cal P}$. En effet, $(\Lambda E_2)|_{\Xi_2}$ est parabolique de poids $12$, donc proportionnelle ˆ $\Delta$.
De mme, $\Pi_1 E_2|_{\Xi_2},\Pi_2 E_2|_{\Xi_2}$ sont proportionnelles ˆ $\Delta$, ce qui implique que $\Phi^\sharp E_2:=\Pi E_2-9(\Lambda E_2)^2$
satisfait $(\Phi^\sharp E_2)|_{\Xi_2}=\lambda\Delta^2$ avec $\lambda\in\CC$. La constante de proportionnalitŽ $\lambda$ est nulle, car $\Phi^\sharp E_2=E_2\Phi
E_2$, donc $E_{4,\ZZ}(\Phi E_2)|_{\Xi_2}=\lambda\Delta^2$, mais il n'existe pas de forme parabolique non nulle de poids $10$ 
pour $\SL_2(\ZZ)$. Donc $H=\Phi E_2\in{\mathfrak P}$.

Le corollaire \ref{lemme:hammond} gŽnŽralise le thŽorme 4.1 p. 507 de \cite{Hammond:Modular}~:
de plus, la mŽthode
utilisŽe dans \cite{Hammond:Modular}
(Žtude de {\em plongements modulaires} d'Igusa) ne detecte pas de forme modulaire antisymŽtrique non nulle de poids $20$.
Noter de plus qu'en gŽnŽral, $M^2\not\in\CC[F,G,H]$. 

\subsection{L'anneau des formes modulaires de tout poids pour $\Gamma$, et sa cl™ture algŽbrique.}

Nous avons explicitŽ un lien entre les propriŽtŽs diffŽrentielles de formes modulaires et
la structure des anneaux des formes modulaires de poids parallle. Nous avons un anneau~:
\[{\cal L}=\bigoplus_{\und{r}}M_{\und{r}}(\Gamma),\]
et nous nous demandons si les propriŽtŽs diffŽrentielles des formes modulaires de poids parallle
permettent de dŽcrire sa structure.

\begin{Lemme}
L'anneau ${\cal L}$ n'est pas de type fini.
\label{lemme:non_finitude}\end{Lemme}
\noindent {\bf Esquisse de dŽmonstration.} On reprŽsente les espaces vectoriels non-nuls de formes modulaires de poids
$(a,b)\in\ZZ^2$ comme des points de $\RR^2$~: nous construisons ainsi un ensemble discret ${\cal E}$ contenu dans le premier
quadrant.

Si l'anneau ${\cal L}$ est de type fini, alors ${\cal E}$ est contenu dans un c™ne d'angle au sommet $<\pi/2$,
car il n'y a pas de formes modulaires non nulles, de poids $(0,a)$ ou $(b,0)$. L'action 
des opŽrateurs $\Pi_i$ induit des applications ${\cal E}\rightarrow\NN^2$. Aucune de ces applications
n'a son image contenue dans un c™ne d'angle au sommet $<\pi/2$~: ${\cal L}$ n'est donc pas un anneau de type fini.

\medskip

On revient aux processus de construction de formes modulaires. Nous en avons pour construire des formes modulaires 
de poids parallle {\em ex-novo}~: sŽries d'Eisenstein, fonctions thta. 

Pour construire des formes modulaires
de poids non parallle, nous ne connaissons jusqu'ˆ ici que des mŽthodes diffŽrentielles~: application des opŽrateurs
$\Pi_i,\Lambda_{i,k},(X,Y)\mapsto[X,Y]_{1_i},\ldots$ ˆ des formes modulaires de poids parallle. Les formules de traces
de \cite{Saito:Operator} peuvent tre aussi utilisŽes.

Voici une question liŽe ˆ ce problme~: posons $K=\QQ(\sqrt{5})$. 
On voit qu'il n'y a pas de formes modulaires de Hilbert de poids non parallle $(a,b)$ avec
$a+b<12$. Nous avons $M_{\und{6}}=\langle \varphi_2^3,\chi_6\rangle$. Nous voulons construire explicitement
des bases de $M_{(a,b)}$ avec $a+b=12$ et $a\not= b$.

Par exemple, $M_{(8,4)}$ contient $\Pi_1E_2$ qui est non nul, et la dimension de cet espace est $1$~:
nous en avons une base.

Les opŽrateurs diffŽrentiels que nous avons introduit
ne permettent pas de calculer des formes modulaires non nulles dans $M_{(a,b)}$
avec \[(a,b)\not\in\{(8,4),(6,6),(4,8)\},\] ˆ partir de formes modulaires de poids
parallle~: {\em a priori} la seule chose que l'on puisse dire est que ces espaces sont de dimension au plus $1$.

Nous savons que si $f$ est une forme modulaire non nulle dans un de ces espaces, elle doit satisfaire
une certaine Žquation diffŽrentielle. Dans chaque cas, on peut construire son dŽveloppement en sŽrie
de Fourier, montrer que sa restriction ˆ $\Xi_2$ n'est pas un multiple 
de $\Delta$, et finalement dŽmontrer que pour $a+b=12$ et $(a,b)\not\in\{(8,4),(6,6),(4,8)\}$,
alors $M_{(a,b)}=\{0\}$.

L'anneau ${\cal L}$ est-il engendrŽ par l'image de ${\cal T}$ par application itŽrŽe de tous les opŽra\-teurs
diffŽrentiels $\Pi_i,\Lambda_{i,k},(X,Y)\mapsto[X,Y]_{1_i},\ldots$?

La rŽponse ˆ cette question est {\em non}. 
Les formules des traces d'opŽrateurs de Hecke permettent de construire des
algorithmes de calcul de la dimension d'espaces de formes modulaires de poids non parallle (cf. \cite{Saito:Operator},
et les algorithmes de calcul dŽveloppŽs dans \cite{Okada:Trace}). 
Par exemple, en utilisant la proposition 2.7 de \cite{Okada:Trace}, on peut montrer que $S_{(14,2)}(\Upsilon)$ est de dimension
$\geq 1$. Il est facile de montrer que $S_{(14,2)}$ ne peut pas contenir l'image d'un opŽrateur diffŽrentiel 
dŽfini sur un espace vectoriel de formes modulaires. 

Des bases d'espaces de formes modulaires
$S_{(a,b)}$ peuvent tre calculŽes explicitement par les formules de traces (cf. \cite{Okada:Trace}),
mais pas par application d'opŽrateurs diffŽrentiels sur des formes modulaires de poids parallle. 

\medskip

Posons~:
\[{\cal M}_0=\CC\left[E_2,\frac{1}{2\pi{\rm i}}\frac{\partial E_2}{\partial z_1},\frac{1}{2\pi{\rm i}}\frac{\partial E_2}{\partial z_2},
\frac{1}{(2\pi{\rm i})^2}\frac{\partial^2 E_2}{\partial z_1 \partial z_2},\frac{1}{(2\pi{\rm i})^2}\frac{\partial^2 E_2}{\partial z_1^2},
\frac{1}{(2\pi{\rm i})^2}\frac{\partial^2 E_2}{\partial z_2^2}\right].\]
Cet anneau n'est pas stable pour les opŽrateurs de dŽrivation 
$(2\pi{\rm i})^{-1}(\partial/\partial z)$ et $(2\pi{\rm i})^{-1}(\partial/\partial z')$, mais
les propriŽtŽs diffŽrentielles des formes modulaires permettent d'obtenir~:
\begin{Theoreme}
La clotžre algŽbrique ${\cal B}$ du corps des fractions de ${\cal M}_0$, munie des opŽrateurs 
de dŽrivation $(2\pi{\rm i})^{-1}(\partial/\partial z)$ et $(2\pi{\rm i})^{-1}(\partial/\partial z')$, est un corps diffŽrentiel
qui est de degrŽ de transcendance $6$ sur $\CC$.
Si $K=\QQ(\sqrt{5})$, alors l'anneau
\[{\cal M}={\cal M}_0[\varphi_2^{-1},\chi_5^{-1},[\varphi_2,\chi_5]_{2_1}^{-1},[\varphi_2,\chi_5]_{2_2}^{-1}],\] muni des opŽrateurs de dŽrivation 
$(2\pi{\rm i})^{-1}(\partial/\partial z)$ et $(2\pi{\rm i})^{-1}(\partial/\partial z')$,
est un anneau differentiel qui est isomorphe ˆ un quotient de l'anneau de polyn™mes $\CC[X_1,\ldots,X_{10}]$ 
par un idŽal premier de hauteur algŽbrique
$4$.\label{lemme:degre6}\end{Theoreme}

\noindent {\bf DŽmonstration}. Les propositions \ref{lemme:ind_alg_diff} et \ref{lemme:resnikoff} impliquent que
les fonctions~:
\begin{eqnarray*}
& & X=E_2,Y_1=\frac{1}{2\pi{\rm i}}\frac{\partial E_2}{\partial z_1},Y_2=\frac{1}{2\pi{\rm i}}\frac{\partial E_2}{\partial z_2},\\ & & Z=
\frac{1}{(2\pi{\rm i})^2}\frac{\partial^2 E_2}{\partial z_1
\partial z_2}, T_1=\frac{1}{(2\pi{\rm i})^2}\frac{\partial^2 E_2}{\partial z_1^2},T_2=\frac{1}{(2\pi{\rm i})^2}
\frac{\partial^2 E_2}{\partial z_2^2},
\end{eqnarray*}
 sont algŽbriquement indŽpendantes sur
$\CC$. Donc le degrŽ de transcendance de ${\cal B}$ est au moins $6$.

Posons~:
\[A_1=\frac{1}{(2\pi{\rm i})^3}\frac{\partial^3 E_2}{\partial z_1^3},A_2=\frac{1}{(2\pi{\rm i})^3}\frac{\partial^3 E_2}{\partial z_2^3},
B_1=\frac{1}{(2\pi{\rm i})^3}\frac{\partial^3 E_2}{\partial z_1^2 \partial z_2},B_2=\frac{1}{(2\pi{\rm i})^3}\frac{\partial^3 E_2}{\partial z_1 \partial z_2^2}.\]
On calcule les formes modulaires $[E_2,\Pi_iE_2]_{1_j}$ pour $i,j=1,2$, qui sont non nulles d'aprs le lemme \ref{lemme:injectivite}. On trouve~:
\begin{eqnarray*}
{[}E_2,\Pi_2E_2{]}_{1_1} & = & 4(B_2X^2-T_2Y_1X+3Y_2(Y_1Y_2-XZ))\\
& = & 4B_2X^2+{\cal B}_2\\
{[}E_2,\Pi_1E_2{]}_{1_2} & = & 4(B_1X^2-T_1Y_2X+3Y_1(Y_2Y_1-XZ))\\
& = & 4B_1X^2+{\cal B}_1\\
{[}E_2,\Pi_1E_2{]}_{1_1} & = & 4(6Y_1^3-6T_1XY_1+A_1X^2)\\
& = & 4A_1X^2+{\cal A}_1\\
{[}E_2,\Pi_2E_2{]}_{1_2} & = & 4(6Y_2^3-6T_2XY_2+A_2X^2)\\
& = & 4A_2X^2+{\cal A}_2,
\end{eqnarray*}
avec ${\cal A}_i,{\cal B}_i\in\CC[X,Y_1,\ldots,T_2]$. Ces formes modulaires 
sont de poids $(8,10),\\ (10,8),(12,6)$ et $(6,12)$ respectivement.
De ces formes modulaires on construit des formes modulaires de poids parallle non nulles~:
\begin{eqnarray*}
{[}E_2,\Pi_1E_2{]}_{1_2}^2\Pi_2E_2,{[}E_2,\Pi_2E_2{]}_{1_2}{[}E_2,\Pi_1E_2{]}_{1_2}
\Pi_1E_2 & \in & S_{\und{24}}\\
{[}E_2,\Pi_1E_2{]}_{1_1}^2(\Pi_2E_2)^3,{[}E_2,\Pi_2E_2{]}_{1_2}^2(\Pi_1E_2)^3 & \in & S_{\und{36}}
\end{eqnarray*}
On dŽduit que les fonctions
$A_1,A_2,B_1,B_2$ se trouvent chacune dans une extension quadratique de
\[{\cal T}(\Gamma)[X,X^{-1},Y_1,\ldots,T_2,(\Pi_1E_2)^{-1},(\Pi_2E_2)^{-1}].\] Mais la
cl™ture algŽbrique du corps de  fractions de cet l'anneau est Žgale ˆ la cl™ture algŽbrique de $\CC(X,Y_1,\ldots,T_2)$ (\footnote{Noter
qu'il y a bien plus de rŽlations. Par exemple~:
\begin{eqnarray*}
{[}E_2,\Pi_2E_2{]}_{1_1}{[}E_2,\Pi_1E_2{]}_{1_2} &\in & S_{\und{18}}\\
{[}E_2,\Pi_1E_2{]}_{1_1}{[}E_2,\Pi_2E_2{]}_{1_2} &\in & S_{\und{18}}\\
{[}E_2,\Pi_2E_2{]}_{1_1}^2\Pi_1E_2,{[}E_2,\Pi_1E_2{]}_{1_1}{[}E_2,\Pi_2E_2{]}_{1_1}
{[}E_2,E_2{]}_{4_2}
& \in & S_{\und{24}}\\
{[}E_2,\Pi_1E_2{]}_{1_1}{[}E_2,\Pi_1E_2{]}_{1_2}(\Pi_2E_2)^2 & \in & S_{\und{30}}\\
{[}E_2,\Pi_2E_2{]}_{1_2}{[}E_2,\Pi_2E_2{]}_{2_1}(\Pi_1E_2)^2 & \in & S_{\und{30}}\\
\ldots & \ldots & \ldots
\end{eqnarray*}
Toutes ces relations ci-dessus sont redondantes~: il y a donc des relations qui sont satisfaites par les formes modulaires
de poids parallle dans $S_{6r}(\Gamma)$, indŽpendantes de $K$.}).
Ceci prouve que $A_1,A_2,B_1,B_2$ sont algŽbri\-ques sur $\CC(X,Y_1,\ldots,T_2)$.

Supposons maintenant que $K=\QQ(\sqrt{5})$, et dŽmontrons les propriŽtŽs de l'anneau ${\cal M}$.
L'anneau ${\cal M}_1$ engendrŽ par toutes les derivŽes partielles de la fonction
$\varphi_2$ est diffŽrentiellement stable. Cet anneau contient $\chi_6$ et $\chi_5^2$, 
car il contient ${\cal M}_0$ qui contient ˆ son tour
$\chi_6$ et $\chi_5^2$ d'aprs le thŽorme \ref{lemme:structure1} (ce sont des images 
de $\varphi_2$ par des opŽrateurs diffŽrentiels d'ordre $2$). Donc toutes les derivŽes partielles de toutes les formes 
modulaires symŽtriques de poids pairs appartiennent ˆ ${\cal M}_1$.

Comme $\chi_{15}$ est proportionnelle ˆ $\chi_5^{-1}[\varphi_2,\chi_5^2,\chi_6]$, on trouve aussi
que toutes les derivŽes partielles de toutes les formes modulaires de poids parallle appartiennent ˆ ${\cal M}_1[\chi_5^{-1}]$
(il suffit d'appliquer le thŽorme \ref{lemme:structure2}). Donc toutes les derivŽes partielles de $\chi_5^{-1}$ appartiennent ˆ ${\cal M}_1[\chi_5^{-1}]$
et il en rŽsulte que cet anneau est diffŽrentiellement stable.

Passons ˆ l'Žtude de l'anneau ${\cal M}$.
On a que ${\cal M}\supset{\cal M}_0[\chi_5^{-1}]$, et ce dernier anneau contient toutes
les formes modulaires de poids parallle, comme l'on a dŽjˆ remarquŽ.

Montrons maintenant que ${\cal M}_1[\chi_5^{-1}]\subset{\cal M}$. Il suffit de dŽmontrer
que $A_1,A_2,B_1,B_2\in{\cal M}$~: nous montrons ici seulement que $A_1,B_1\in{\cal M}$.
On a que $[\varphi_2,\Pi_1\varphi_2]_{1_1}[\varphi_2,\chi_5]_{1_2}^3$ est une forme parabolique 
non nulle de poids parallle $33$, et que $[\varphi_2,\Pi_1\varphi_2]_{1_1}[\varphi_2,\chi_5]_{1_2}$
est une forme parabolique 
non nulle de poids parallle $17$. Comme \[[\varphi_2,\Pi_1\varphi_2]_{1_1}=4\varphi_2^2A_1+P,\quad
[\varphi_2,\Pi_1\varphi_2]_{2_1}=4\varphi_2^2B_1+Q\]
avec $P,Q\in{\cal M}_0$, et comme ${\cal M}$ contient toutes les formes modulaires de poids parallle,
nous obtenons ${\cal M}_1[\chi_5^{-1}]\subset{\cal M}$. Comme ${\cal M}$ contient aussi
toutes les derivŽes partielles de $\varphi_2^{-1},\chi_5^{-1},[\varphi_2,\chi_5]_{1_1}^{-1}$ et $[\varphi_2,\chi_5]_{1_2}^{-1}$, 
${\cal M}$ est diffŽrentiellement
stable.

L'anneau ${\cal M}$ est clairement un quotient de $\CC[X_1,\ldots,X_{10}]$ par un idŽal de rŽlations
${\cal I}$. De plus, comme ${\cal M}$ est un anneau de fonctions mŽromorphes, il n'a pas de diviseurs de zŽro,
et donc ${\cal I}$ est un idŽal premier.
Comme le degrŽ de transcendance de ${\cal B}$ est $6$, l'idŽal ${\cal I}$ a hauteur $4$ (\footnote{Pour des applications
envisagŽables de ce thŽorme, il convient de remarquer que ${\cal M}$ est un anneau de Cohen-Macaulay, puisque
${\cal I}$ est clairement un idŽal d'intersection complte. Les rŽlations dŽfinissant ${\cal I}$ peuvent tre
explicitŽes, mais elles sont assez compliquŽes.}).

La preuve du thŽorme 
\ref{lemme:degre6} est terminŽe.
On montre aussi que la cl™ture algŽbri\-que du corps engendrŽ par toutes les dŽrivŽes partielles de fonctions modulaires
de poids $\und{0}$ est de degrŽ de transcendance $6$.

\medskip

\noindent {\bf Remarque.} 
Connaissant la structure de l'anneau de formes modulaires de poids parallle pour un certain
groupe modulaire de Hilbert $\Gamma$, on peut dŽcrire les relations algŽbriques liant les fonctions
$A_i,B_j$. Naturellement, ces relations sont en gŽnŽral trs compliquŽes. Attention~: le corps ${\cal B}$ {\em n'est
pas} la cl™ture algŽbrique du corps des fractions de l'anneau ${\cal L}$ des formes modulaires de tout poids.
On a la proposition qui suit.

\begin{Proposition}
\label{lemme:quatre_formes}
Cinq formes modulaires de Hilbert sont toujours algŽbri\-quement dŽpendantes. Il existe
quatre formes modulaires de Hilbert algŽbri\-quement indŽpendantes.
\end{Proposition}
{\bf DŽmonstration.} Le thŽorme \ref{lemme:structure}
implique qu'il existe trois formes modulaires de poids parallles et algŽbriquement
indŽpendantes. De plus, on peut facilement vŽrifier que si $X$ est une forme modulaire de Hilbert de poids non parallle non nulle, et si $F$
est une forme modulaire de Hilbert de poids parallle non constante, alors $F,X$ sont algŽbriquement indŽpendantes.
Dans ce texte nous avons construit plusieurs formes modulaires de Hilbert non constantes, de poids non parallle.
Ainsi, on peut toujours construire quatre formes modulaires de Hilbert algŽbriquement indŽpendantes.

Soient maintenant $F,G,H,X,Y$ cinq formes modulaires non constantes, avec $F,G,H$ de poids parallle $f,g,h$ et $X,Y$ de poids non parallles
$\und{x},\und{y}$. 
Supposons que $F,G,H$ soient algŽbriquement 
indŽpendantes, et que $X,Y$ soient aussi algŽbriquement indŽpendantes~: on a aussi que $F,G,H,X$ sont algŽbriquement indŽpendantes,
et que $F,G,H,Y$ sont algŽbriquement indŽ\-pendantes. 
Si $\und{x}\in\QQ^\times\und{y}$, alors il existe deux entiers $a,b$ tels que $X^a/Y^b$ soit une fonction modulaire $Q$ de poids parallle $0$,
non constante. Mais $Q$ est un quotient de deux formes modulaires de poids parallles, donc $Q\in\CC(F,G,H)$, et $F,G,H,X,Y$ sont
algŽbriquement dŽpendantes. Supposons maintenant que $\und{x}\not\in\QQ^\times\und{y}$~: donc les
vecteurs $\und{x}$ et $\und{y}$ sont $\QQ$-linŽairement indŽpendants. Il existe une combinaison linŽaire
$a\und{x}+b\und{y}\in(\ZZ-\{0\})\und{1}$, avec $a,b\in\ZZ$, et la forme modulaire $X^aY^b$ est une forme modulaire non constante $P$ de poids parallle.
Nous avons dŽjˆ remarquŽ que quatre formes modulaires de Hilbert sont toujours algŽbriquement dŽpendantes, donc 
$F,G,H,P$ sont algŽbriquement dŽpendantes et on a aussi dans ce cas que
$F,G,H,X,Y$ sont algŽbriquement dŽpendantes.

Nous laissons au lecteur le soin de complŽter la dŽmonstration de cette proposition dans le cas gŽnŽral o
$F,G,H$ n'ont pas nŽcessairement de poids parallles, en utilisant ce qui precde.

\end{document}